 \documentclass[final,3p,times]{elsarticle}

\usepackage{amssymb}

\usepackage{amsthm}

\theoremstyle{definition}

\usepackage{amsmath,bm}

\usepackage{lineno}
\usepackage{algorithm}
\usepackage[noend]{algpseudocode}
\usepackage{booktabs} 
\usepackage{standalone}
\usepackage{tikz}
\usepackage{pgfplots}
\pgfplotsset{compat=1.13}
\usetikzlibrary{calc}
\usetikzlibrary{arrows}
\usetikzlibrary{shapes,shapes.multipart, shapes.geometric, arrows}
\usetikzlibrary{positioning}
\usetikzlibrary{decorations.markings}
\usetikzlibrary{decorations.pathreplacing}
\usepackage{placeins} 
\usepackage{subcaption}
\usepackage{graphicx}
\usepackage{xargs}                      
\usepackage{verbatim}
\usepackage{multicol}
\usepackage{hyperref}

\tikzstyle{startstop} = [rectangle, rounded corners, minimum width=3cm,minimum height=1cm,text centered,text width=12cm, draw=black, fill=red!30]
\tikzstyle{io} = [trapezium, trapezium left angle=70, trapezium right angle=110, minimum width=1cm, minimum height=1cm, text width=6cm,text centered, draw=black, fill=blue!30]
\tikzstyle{process} = [rectangle, minimum width=3cm, minimum height=1cm, text centered, text width=5cm, draw=black, fill=orange!30]
\tikzstyle{decision} = [ellipse, minimum width=3cm, minimum height=1cm, text centered, draw=black, fill=green!30]
\tikzstyle{blackbox} = [rectangle, minimum width=3cm, minimum height=1cm, text centered, text width=5cm, draw=black, fill=gray!30]
\tikzstyle{arrow} = [thick,->,>=stealth]

\usepackage[colorinlistoftodos,prependcaption,textsize=tiny,textwidth=2cm,]{todonotes}
\usepackage{xcolor}
\newcommandx{\unsure}[2][1=]{\todo[linecolor=red,backgroundcolor=red!25,bordercolor=red,#1]{#2}}
\newcommandx{\change}[2][1=]{\todo[linecolor=blue,backgroundcolor=blue!25,bordercolor=blue,#1]{#2}}
\newcommandx{\info}[2][1=]{\todo[linecolor=green,backgroundcolor=green!25,bordercolor=green,#1]{#2}}
\newcommandx{\improvement}[2][1=]{\todo[linecolor=purple,backgroundcolor=purple!25,bordercolor=purple,#1]{#2}}
\newcommandx{\thiswillnotshow}[2][1=]{\todo[disable,#1]{#2}}

\usepackage{threeparttable}

\newcommand{\myset}[3]{\{#1_{#2}\}^{#3}_{#2=1}}

\newcommand{\x}{\mathbf{x}}
\newcommand{\y}{\mathbf{y}}
\newcommand{\X}{\mathbf{X}}

\newcommand{\kernelM}[2]{M\left(\mathbf{#1},\mathbf{#2} \right)}

\makeatletter
\def\BState{\State\hskip-\ALG@thistlm}
\makeatother

\begin{document}

\begin{frontmatter}



\title{An integral equation based numerical method for the forced heat equation on
complex domains}


%

\author[KTH]{Fredrik Fryklund\corref{cor1}}
\ead{ffry@kth.se}
\cortext[cor1]{Corresponding author}
\author[SFU]{Mary Catherine A. Kropinski}
\ead{mkropins@sfu.ca}
\author[KTH]{Anna-Karin Tornberg}
\ead{akto@kth.se}

\address[KTH]{Department of Mathematics, KTH Royal Institute of Technology, Stockholm, Sweden}
\address[SFU]{Department of Mathematics, Simon Fraser University, Burnaby, Canada}

\begin{abstract}

Integral equation based numerical methods are directly applicable to homogenous elliptic PDEs, and offer the ability to solve these with high accuracy and speed on complex domains. In this paper, extensions to problems with inhomogeneous source terms and time dependent PDEs, such as the heat equation, have been introduced. One such approach for the heat equation is to first discretize in time, and in each time-step solve a so-called modified Helmholtz equation with a parameter depending on the time step size. The modified Helmholtz equation is then split into two parts: a homogenous part solved with a boundary integral method and a particular part, where the solution is obtained by evaluating a volume potential over the inhomogeneous source term over a simple domain. In this work, we introduce two components which are critical for the success of this approach: a method to efficiently compute a high-regularity extension of a function outside the domain where it is defined, and a special quadrature method to accurately evaluate singular and nearly singular integrals in the integral formulation of the modified Helmholtz equation for all time step sizes.

\end{abstract}

\end{frontmatter}



\section{Introduction}
\label{s:introduction}
In this paper we present a highly accurate numerical method for solving the forced isotropic heat equation with Dirichlet data on complex multiple  connected domains in two dimensions. We adapt the solution methodology introduced by \citeauthor{Kropinski2011Heat} in \cite{Kropinski2011Heat}, but extend and generalise their work to allow for solution of a wider class of problems with improved discretisation in time and uniform accuracy all the way up to the boundary. First, the heat equation is discretized in time with an implicit treatment of the diffusion term, an approach that is sometimes referred to as Rothe's method \cite{chapko1997,chapko2001} or elliptic marching. This results in a sequence of modified Helmholtz equations, also known as the linearised Poisson-Boltzmann equation, to be solved at each time step. Doing so advances the solution to the parabolic heat equation in time. A relaxed definition of the modified Helmholtz equation reads $\alpha^{2}u-\Delta u = f$, with
$\alpha^{2}$ inversely proportional to the time step. Utilising the linearity, this equation is further split into two: one that finds a particular solution for the specific right hand side without enforcing the boundary conditions, and a homogeneous problem that ensures that the sum of the two solutions solves the original problem.  The homogeneous problem is solved with a boundary integral method with a panel-based Nystr\"{o}m quadrature scheme, as introduced in \cite{Kropinski2011modHelm} by Kropinski and Quaife.
The particular solution is written as a volume potential with the free space Green's function for the modified Helmholtz equation, also known as the Yukawa-or screened Columb potential. To avoid constructing quadrature methods for the evaluation of this volume potential over complex domains, an extension of the right hand side $f$ is introduced, allowing for integration over a simple rectangular domain.

In \cite{Kropinski2011Heat}  the authors \citeauthor{Kropinski2011Heat} demonstrated the potential of developing an efficient and accurate general boundary integral solver for the heat equation on complex domains. Moreover, they list the major remaining issues that require further investigation. At that time only examples for which a continuous extension of $f$ could be constructed by hand was considered, thus excluding complex geometries and general data. Another impediment was the loss of accuracy for evaluating layer potentials close to their sources. Their solution was to over-resolve the boundary, but the loss of accuracy is still significant as a target point approaches the boundary. In this paper we introduce the following developments:
\begin{itemize}
  \item High order adaptive methods for time evolution.
  \item A method to efficiently compute a high-regularity extension of a function $f$  to an enclosing and geometrically simple domain, given only its values at discrete locations in $\Omega$.
  \item
A special purpose quadrature method to avoid loss of accuracy when evaluating layer potentials close to the boundary and the kernel becomes nearly singular.
\end{itemize}

Two main groups of semi-implicit time stepping methods are Runge-Kutta methods \cite{KENNEDY2003139} and spectral deferred correction methods \cite{Dutt2000,JIA20081739,minion2003}. We use the former to obtain an adaptive scheme, but the approach we propose is general with respect to the choice of semi-implicit time stepper.

It is not a simple problem to construct a high regularity extension of a function, for which only its values are known in discrete points inside the original domain $\Omega$. In \cite{Kropinski2011Heat}, Kropinski and Quaife considered only examples for which a continuous extension could be constructed by hand. We use a partition of unity extension technique (PUX) by \citeauthor*{FRYKLUNDPUX} in \cite{FRYKLUNDPUX}. They solve the Poisson equation with the above-mentioned split into a particular and an homogeneous problem. We now use this method for function extension in the context of the modified Helmholtz equation with excellent results and can hence increase the class of solvable problems as compared to \cite{Kropinski2011Heat}. An alternative approach for function extension is given in \cite{ASKHAM20171}, where the function to be extended outside of
$\Omega$ sets the boundary Dirichlet data on $\partial \Omega$ for the Laplace equation in $\mathbb{R}^{2} \setminus \Omega$. The solution to this problem is computed with an integral equation based method, and defines a continuous function extension.  See \cite{FRYKLUNDPUX} and the references therein for other extension techniques, such as Fourier continuation methods or extending the unknown solution or solution from previous time step \cite{BRUNO20102009,STEIN2017155,Shirokoff2015}.

When evaluating a layer potential close to a boundary, the kernel becomes nearly singular.
A well known challenge with boundary integral based methods is accurate numerical integration of singular (for evaluation on the boundary) and nearly singular kernels.
The comparative study \cite{Hao2014} complemented with \cite{Helsing2015} give an overview of state of the art methods. The latter includes panel-based explicit kernel-split schemes with product integration, pioneered by \citeauthor{specialQuad} \cite{specialQuad} for the Laplace equation. This methodology is applicable to a large class of linear elliptic PDEs, and achieves excellent results also for e.g. the Helmholtz  \cite{Helsing2015} and Stokes equations \cite{OJALA2015145}.  However, for the modified Helmholtz equation product integration may fail altogether for sufficiently large $\alpha$, i.e. for small time steps in our setting. The quadrature rule will in this case require an unfeasibly high resolution of the boundary, which is not motivated by the geometry nor the resolution requirement for the layer density.
This spurred the development of a quadrature scheme to solve this problem.
The Yukawa potential decays as $\exp \left(-\alpha\right)$, and the kernel becomes more localised as $\alpha$ increases. In this process, the product integration requires an increasing amount of upsampling, but only over a decreasing interval, and hence only local upsampling is needed. In a separate paper \cite{2019arXiv190607713K}, we present an adaptive quadrature scheme in the spirit of \cite{HELSING20098892} that lifts the previous restriction on $\alpha$.

A parallel development of a boundary integral based solver for the heat equation is based on direct approximation of the heat kernel, thus avoiding discretisation of the differential operator with respect to time. In the initial work \cite{doi:10.1137/080732389} it was observed that to achieve the desired accuracy for domains with high curvature the time step must be considerably smaller than the formal rate of convergence would suggest. The authors refer to this as \emph{geometrically induced stiffness}. In recent work towards solving the heat equation with said method \citeauthor{Wang2019} has developed a \emph{hybrid method} that allows for evaluation of the boundary and volume potentials including the space-time heat kernel without the constraints from geometric stiffness \cite{Wang2019}.

Efforts to solve the heat equation with boundary integral equation based techniques are not only motivated by that specific task. Surely, there are other methods to solve the heat equation on a complex domain, such as finite element methods.  However, the algorithmic development in these efforts is essential to increase the applicability of integral equation based numerical methods which sport several attractive features, including that complex geometry naturally enters the problem and generation of an unstructured mesh is redundant, ill-conditioning associated with discretising the operators is avoided, high accuracy can be attained, and boundary data and far field conditions are simple to incorporate.  Developments for the heat equation are also related to extension from Stokes to Navier-Stokes equations.

The focus of this paper is on the heat equation. However, fast integral equations for the modified Helmholtz equation are of interest for the many applications that equation applies to. These include, but are not limited to: electrostatic interactions in protein and related biological functions, macroscopic electrostatics, Debye–Huckel theory for dilute electrolytes, water wave problems, in the linearisation of the Poisson–Boltzmann equation and approximation of surfaces \cite{proteinstructurePossinBoltz,JUFFER1991144,CHEN2006551,Vorobjev2019,LIANG19971830,KOUIBIA2019262}. Consequently, there is active research on solution methods and analysis thereof for the modified Helmholtz equation. In \cite{chen_jiang_chen_yao_2015} the method of fundamental solution is used, while in \cite{LI200666} it is solved by plane wave functions.

\subsection{Overview of the paper}
The mathematical problem is formulated in Section \ref{s:formulation}, both for the heat equation and the modified Helmholtz equation. Section \ref{s:discretisation} contains the numerical methods for solving the homogeneous problem and the particular problem for the modified Helmholtz equation, including an introduction to PUX.  It is assumed that the heat equation as been appropriately discretised in time. Thereafter we present the numerical results in Section \ref{s:results}, for the modified Helmholtz equation, the heat equation and a reaction-diffusion type problem. Finally we present our conclusions and an outlook in Section \ref{s:conclusions}.
See \ref{ss:discretisation_temporal} for instructions on how IMEX Runge-Kutta methods reduce the heat equation to a sequence of modified Helmholtz equations. There are simple examples, Butcher tableaus and a note on adaptivity. In \ref{appendix:flowchart} we present a graphical overview of the solution procedure for the modified Helmholtz equation.

\section{Formulation}
\label{s:formulation}
\begin{figure}
  \centering
  \includegraphics[width=0.4\textwidth]{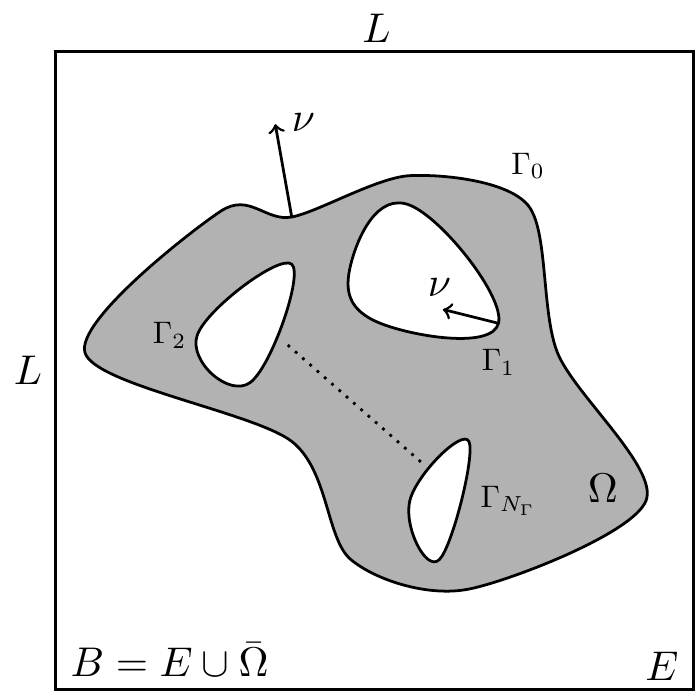}
\caption{The heat equation \eqref{eq:HeatEq}--\eqref{eq:HeatEqBC} is defined in $\Omega$. It is enclosed in a box $B=[L,L]^2$. The boundaries are denoted $\Gamma_{n}$,  $n = 0,\ldots,N_{\Gamma}$. The outer boundary is $\Gamma_{0}$ and the outward directed normal is denoted by $\nu$.}
\label{fig:Omega}
\end{figure}
Consider the forced isotropic heat equation
\begin{align}
  \frac{\partial U(t,\mathbf{x})}{\partial t} - \Delta U(t,\mathbf{x}) &= F(t,\mathbf{x}),\quad t_{0}<t,\quad \mathbf{x}\in \Omega\subset\mathbb{R}^{2},\label{eq:HeatEq}\\
  U(t_{0},\mathbf{x}) &= U_{0}(\mathbf{x}),\quad \mathbf{x}\in  \Omega,\label{eq:HeatEqIC}\\
  U(t,\mathbf{x}) &= g(t,\mathbf{x}),\quad \mathbf{x}\in \Gamma,\label{eq:HeatEqBC}
\end{align}
subject to initial- and Dirichlet boundary data  $U_{0}$ and $g$, respectively. To fix notation let $\Omega$ be a time
independent, compact $(N_{\Gamma}  + 1)$-ply connected region in $\mathbb{R}^{2}$ with a boundary
$\Gamma$ consisting  of $(N_{\Gamma} + 1)$ closed curves. These form the set
$\Gamma = \{\Gamma_{n}\}_{n=0}^{N_{\Gamma}}$, where $\Gamma_{0}$ is the outer boundary of the
region $\Omega$, see Fig. \ref{fig:Omega}. The component curves are individually smooth and parametrisation each is assumed to be known.
The outward directed normal at $\y \in \Gamma$ is denoted $\bm{\nu}(\y) = \bm{\nu}_{\y}$ and $\kappa(\y)$ denotes the curvature at $\y \in \Gamma$.
\subsection{Discretising in time and the modified Helmholtz equation}
\label{ss:formulation_discintime}
The heat equation \eqref{eq:HeatEq} is first discretised in time, an approach known as elliptic marching or Rothe's method. To prevent severe time step  restrictions an implicit-explicit (IMEX) scheme is used. It consists of using an implicit discretisation of the stiff terms and an explicit one for the nonstiff terms \cite{AscherIMEX}. Regardless of the specifics of the IMEX scheme, to advance the solution $U$ in discrete time a sequence of modified Helmholtz equations are solved. The modified Helmholtz equation is stated as
\begin{align}
 \label{eq:ModHelmEq}
 \alpha^{2}u(\mathbf{x}) - \Delta u(\mathbf{x}) = f(\mathbf{x}),\quad \mathbf{x} \in \Omega,\\
 \label{eq:ModHelmEqBC}
 u(\mathbf{x}) = g(\mathbf{x}),\quad \mathbf{x} \in \Gamma,
\end{align}
with $u$ unknown in $\Omega$. The scalar parameter $\alpha^{2}$ is inversely proportional to the time step $\delta t$; its explicit form along with $f$ and $g$ are given by the specific IMEX scheme. We use an adaptive IMEX Runge-Kutta method of fourth order in this paper, see \ref{ss:discretisation_temporal}. However, what follows holds for any IMEX scheme.

  Using the linearity of the differential operator $\alpha^{2} - \Delta$, the solution $u$ to \eqref{eq:ModHelmEq}--\eqref{eq:ModHelmEqBC} is decomposed into a homogeneous solution $u^{H}$ and a particular solution $u^{P}$, such that $u(\mathbf{x}) = u^{H}(\mathbf{x}) + u^{P}(\mathbf{x})$ for $\mathbf{x} \in \Omega$. First the particular solution is acquired by solving a free space problem
  \begin{align}
    \label{eq:ModHelmEq_part}
      \alpha^{2} u^{P}(\x) - \Delta u^{P}(\x) = f^{e}(\x),& \quad \mathbf{x} \in \mathbb{R}^{2},\\
      u(\x)\rightarrow 0,&\quad |\x|\rightarrow \infty,
      \label{eq:ModHelmEq_partBC}
  \end{align}
assuming the existence of an extension $f^{e}\in C^{k}(\mathbb{R}^{2})$, for some $k\geq0$, of the right hand side $f$ from \eqref{eq:ModHelmEq}, such that
  \begin{align}
    \label{eq:ModHelmEq_partfe}
    f^{e}(\mathbf{x}) &= f(\mathbf{x}),\, \forall \mathbf{x}\in \Omega,\\
    \label{eq:ModHelmEq_partfe2}
    \text{supp}(f^{e})&\subset B = [-L,L]^{2},
  \end{align}
  for some finite $L$. The boundary condition, given by the Dirichlet data $g$ in \eqref{eq:ModHelmEqBC}, is satisfied by $u$ if $u^{H}$ is a solution to
  \begin{align}
    \alpha^{2}u^{H} - \Delta u^{H} = 0,\quad \mathbf{x} \in \Omega,\label{eq:ModHelmEq_homo}\\
    \label{eq:ModHelmEqBC_homo}
    u^{H} = \tilde{g}(\mathbf{x})=g(\mathbf{x})-u^{P}(\mathbf{x})\vert_{\Gamma},\quad \mathbf{x} \in \Gamma.
  \end{align}

In short, first solve the free space problem \eqref{eq:ModHelmEq_part} to obtain the boundary data for the homogeneous problem \eqref{eq:ModHelmEq_homo}--\eqref{eq:ModHelmEqBC_homo}. The solution to the modified Helmholtz equation is the sum of the two solutions, $u(\x)=u^{H}(\x)+u^{P}(\x)$ for $\x \in \Omega$. See the flowchart in \ref{appendix:flowchart} for a graphical overview.

\subsubsection{The inhomogeneous modified Helmholtz equation}
\label{sss:fomulation_modhelm_inhomo}
The free-space modified Helmholtz equation \eqref{eq:ModHelmEq_part}--\eqref{eq:ModHelmEq_partBC} can be solved with Fourier transforms.  Let $\hat{u}^{P} = \hat{u}^{P}(\bm{\xi})$ and $\hat{f}^{e}=\hat{f}^{e}(\bm{\xi})$ denote the Fourier transforms for $u^{P}$ and $f^{e}$, respectively. Here $\bm{\xi} = [\xi_{1},\, \xi_{2}]\in \mathbb{R}^{2}$ with $\xi = |\bm{\xi}|$. Then under the Fourier transform \eqref{eq:ModHelmEq_part} is
\begin{equation}
  \alpha^{2}\hat{u}^{P}(\bm{\xi}) + \xi^{2}\hat{u}^{p}(\bm{\xi}) = \hat{f}^{e}(\bm{\xi}), \quad \bm{\xi} \in \mathbb{R}^{2}
\end{equation}
and we obtain
\begin{equation}
  \label{eq:solveFourierCoeff}
  \hat{u}^{P}(\bm{\xi}) = \frac{\hat{f}^{e}(\bm{\xi})}{\alpha^{2} + \xi^{2}}, \quad \bm{\xi} \in \mathbb{R}^{2}.
\end{equation}
Note that the above expression is free of singularities, since $\alpha^{2}\neq 0$. The solution is given by the inverse Fourier transform
\begin{equation}
  \label{eq:inverseFourier}
  u^{P}(\mathbf{x}) = \frac{1}{(2\pi)^{2}}\int\limits_{\mathbb{R}^{2}} \hat{u}^{P}(\bm{\xi})e^{i\bm{\xi}\cdot \mathbf{x}}\,\text{d}\bm{\xi}.
\end{equation}
For this solution to be well-defined the extension $f^{e}$ must be in $L^{1}(\mathbb{R}^{2})$. How to construct said extension and compute an approximation of $u^{P}$ is presented in \ref{sss:discretisation_adaptive_modhelm_inhomo}.

\subsubsection{The homogeneous modified Helmholtz equation}
\label{sss:fomulation_modhelm_homo}
Consider the homogeneous modified Helmholtz equation \eqref{eq:ModHelmEq_homo}--\eqref{eq:ModHelmEqBC_homo}. The free-space Green's function $G(\x,\y)$ for the operator $\alpha^{2}-\Delta$ is
\begin{equation}
  \label{eq:greenFunc}
  G(\x,\y) = \frac{\alpha^{2}}{2 \pi} K_{0}\left(\alpha\|\mathbf{y}-\mathbf{x}\|\right),
\end{equation}
where $K_{0}$ denotes the zeroth-order modified Bessel function of the second kind.  In other contexts the kernel $G(\x,\y)$ is also referred to as the Yukawa or screened Coulomb potential. As in \cite{Kropinski2011Heat,Kropinski2011modHelm}, we seek the solution $u^{H}(\mathbf{x})$ for $\bf{x}\in \Omega$ in the form of a double layer potential:
\begin{equation}
\label{eq:dblPot}
u^{H}(\mathbf{x}) = \frac{\alpha^{2}}{2\pi} \int\limits_{\Gamma}\kernelM{x}{y}\mu(\mathbf{y})\,ds_{\mathbf{y}},\quad \forall \mathbf{x}\in \Omega,
\end{equation}
 with the  kernel
\begin{equation}
\label{eq:kerM}
\kernelM{x}{y} = -\frac{\partial}{\partial \nu_{\mathbf{y}}}  K_{0}\left(\alpha\|\mathbf{y}-\mathbf{x}\|\right)=  \alpha K_{1}\left(\alpha\|\mathbf{y}-\mathbf{x}\|\right)\frac{\mathbf{y}-\mathbf{x}}{\|\mathbf{y}-\mathbf{x}\|}\cdot \nu_{\mathbf{y}},
\end{equation}
where $K_{1}$ denotes the first-order modified Bessel function of the second kind. The limiting value as $\x$ goes to $\y$ along a  boundary segment $\Gamma_{n}$ is well defined:
\begin{equation}
\label{eq:Mlimit}
\lim_{\mathbf{x}\rightarrow \mathbf{y}} M(\mathbf{x},\mathbf{y}) = -\frac{1}{2\pi}\kappa(\mathbf{y}),\quad\x,\,\mathbf{y}\in\Gamma_{n},
\end{equation}
where $\kappa(\y)$ is the curvature of $\Gamma_{n}$ at $\y\in\Gamma_{n}$, $n=1,\ldots,N_{\Gamma}$. The density $\mu:\Gamma \rightarrow \mathbb{R}$ is not known a priori; it is found through the solution of a boundary integral equation. Such an equation of the second kind for $\mu$ can be formulated as
\begin{equation}
\label{eq:Density}
\mu(\mathbf{x}) + \frac{1}{\pi}\int\limits_{\Gamma}\kernelM{x}{y}\mu(\mathbf{y})\,ds_{\mathbf{y}} = -\frac{2}{\alpha^{2}}\tilde{g}(\mathbf{x}),\quad \forall \mathbf{x}\in \Gamma.
\end{equation}
For a derivation see e.g.  \cite{QUAIFETHESIS}. For $\tilde{g}\equiv 0$ only the trivial solution $\mu\equiv 0$ along $\Gamma$ satisfies \eqref{eq:Density}. Thus by the Fredholm alternative the solution $\mu$ exists and is unique for any integrable $\tilde{g}$, for both simply and multiply connected domains \cite{atkinson1997numerical}. This property is inherited by the corresponding discretised systems as well, introduced in Section \ref{sss:discretisation_adaptive_modhelm_homo}.

Each contour $\Gamma_{n}$ is split into $N_{P,n}$ intervals, referred to as panels, where $\Gamma_{n,k}$ is the $k$th panel on the $n$th contour and $N_{P}$ the total number of panels over $\Gamma$. A panel $\Gamma_{n,k}$ is represented by a known parametrisation $\bm{\gamma}_{n,k}$, such that
\begin{equation}
   \label{eq:gammaparam}
   \Gamma_{n,k} = \{\bm{\gamma}_{n,k}(t)\,|\,t\in[-1,1]\}.
\end{equation}
By introducing a speed function $s_{n,k}(t) = |\bm{\gamma}^{\prime}_{n,k}(t)|$ and $\mu_{n,k}(t) = \mu(\bm{\gamma}_{n,k}(t))$ the layer potential \eqref{eq:dblPot} can be written as
\begin{equation}
\label{eq:dblPotpan}
u^{H}(\mathbf{x}) = \frac{\alpha^{2}}{2\pi} \sum\limits_{n=0}^{N_{\Gamma}}\sum\limits_{k=1}^{N_{P,n}}\int\limits_{-1}^{1}M\left(\bm{x},\bm{\gamma}_{n,k}(t)\right)\mu_{n,k}(t)s_{n,k}(t)\,dt,\quad \forall \mathbf{x}\in \Omega
\end{equation}
and analogously for the boundary integral equation \eqref{eq:Density}
\begin{equation}
 \label{eq:semidischomo}
 \mu(\mathbf{x}) + \frac{1}{\pi}\sum\limits_{n=0}^{N_{\Gamma}}\sum\limits_{k=1}^{N_{P,n}}\int\limits_{-1}^{1}M\left(\bm{x},\bm{\gamma}_{n,k}(t)\right)\mu_{n,k}(t)s_{n,k}(t)\,dt = -\frac{2}{\alpha^{2}}\tilde{g}(\mathbf{x}),\quad \forall \mathbf{x}\in \Gamma.
\end{equation}

\section{Discretisation}
\label{s:discretisation}

 This section covers the numerical treatment of the modified Helmholtz equation. Note that two different methods are needed, one for the inhomogeneous problem and one for the homogeneous problem.  We assume some suitable IMEX scheme has been chosen for temporal discretisation of the heat equation \eqref{eq:HeatEq}--\eqref{eq:HeatEqBC}, e.g. the Runge-Kutta methods presented in \ref{ss:discretisation_temporal}.

Consider a box $B = [-L,L]^{2}$ in $\mathbb{R}^{2}$ that contains $\bar{\Omega}$. The complement of $\bar{\Omega}$ relative to $B$ is denoted by $E$. Denote the \emph{grid} by $\mathbf{X}$, which is a set of $N_{u}^{2}$ elements $\x$, referred to as nodes or points. They are uniformly distributed with spacing $\delta x$ over $B$. Let subscripts indicate subsets of $\mathbf{X}$, such as $\mathbf{X}_{\Omega} = \{\mathbf{x}\in \mathbf{X}|\x \in \Omega\}$ and $\mathbf{X}_{E} = \{\mathbf{x}\in \mathbf{X}|\x \in E\}$.

 $\mathbf{X}_{E} = \{\mathbf{x}\in \mathbf{X}\,|\,\x \in E\}$ such that it satisfies \eqref{eq:ModHelmEq_partfe}--\eqref{eq:ModHelmEq_partfe2}.  Thereafter we consider the homogeneous problem \eqref{eq:ModHelmEq_homo}--\eqref{eq:ModHelmEqBC_homo}, formulated as a boundary integral equation on $\Gamma$. The solution is computed at the locations $\mathbf{X}_{\Omega} = \{\mathbf{x}\in \mathbf{X}\,|\,\x \in \Omega\}$ in a post-processing step.
 The solution to the modified Helmholtz equation is computed at all grid points that fall inside $\Omega$, i.e. the elements of $\mathbf{X}_{\Omega}$. First we present how to find this solution for the free space problem \eqref{eq:ModHelmEq_part}--\eqref{eq:ModHelmEq_partBC}. This involves extending the function $f$, based on the data at $\mathbf{X}_{\Omega} = \{\mathbf{x}\in \mathbf{X}\,|\,\x \in \Omega\}$ to




\label{ss:discretisation_adaptive_modhelm}

\subsection{The inhomogeneous problem and function extension}
\label{sss:discretisation_adaptive_modhelm_inhomo}
An approximate solution to the free-space problem \eqref{eq:ModHelmEq_part}--\eqref{eq:ModHelmEq_partBC}  is computed by discretising the integral in \eqref{eq:inverseFourier} with the trapezoidal rule. It is evaluated efficiently with FFTs on the regular grid $\mathbf{X}$ in $B$, thus in $\mathbf{X}_{\Omega}$ as well, and on the boundary $\Gamma$ with a non-uniform inverse FFT. The latter is used to modify the given Dirichlet boundary data  \eqref{eq:ModHelmEqBC} for the homogeneous modified Helmholtz equation.

If the compactly supported $f^e$ in \eqref{eq:ModHelmEq_part} is smooth, then the coefficients in the Fourier series expansion decay exponentially fast with the wave number, and this procedure would be specially accurate. With limited regularity, the Fourier coefficients instead decay algebraically, with one additional order for each continuous derivative. This approach requires an extension $f^{e}$ of $f$ defined on $\mathbf{X}$, preferably with high global regularity and compact support. It is constructed with PUX, which is briefly reviewed in this subsection. The basic concept is to blend local extensions by a partition of unity into a global extension with compact support, enforced by weight functions. The global regularity of the extension is directly related to the construction of said partition of unity. This is achieved by distributing overlapping partitions along the boundary $\Gamma$ of $\Omega$. In each partition the local values of $f$ are used to extend it to the points in the partition that fall outside $\Omega$. For a more extensive treatment see the original work \cite{FRYKLUNDPUX}.
\label{sss:discretisation_aptive_modhelm_pux}

\subsubsection{Partition of unity}
\label{sss:PUX}
\begin{figure}[ht]
    \begin{center}

\begin{tikzpicture}
\node (starfish) {\includegraphics[trim=7.5cm 10cm 7cm 10cm, clip, width=0.4\textwidth]{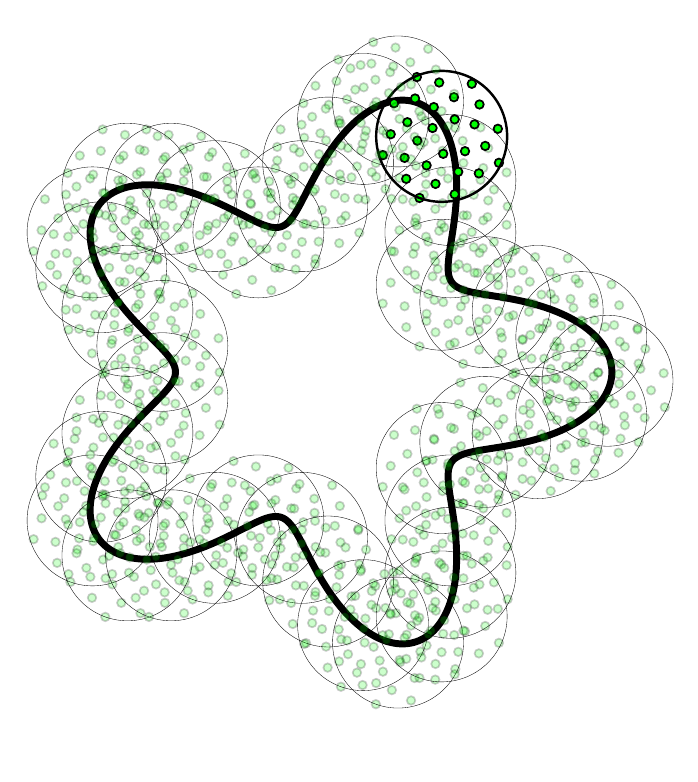}};
\node (partition) [right = of starfish] {\includegraphics[trim=5.5cm 9cm 5cm 8cm, clip, width=0.4\textwidth]{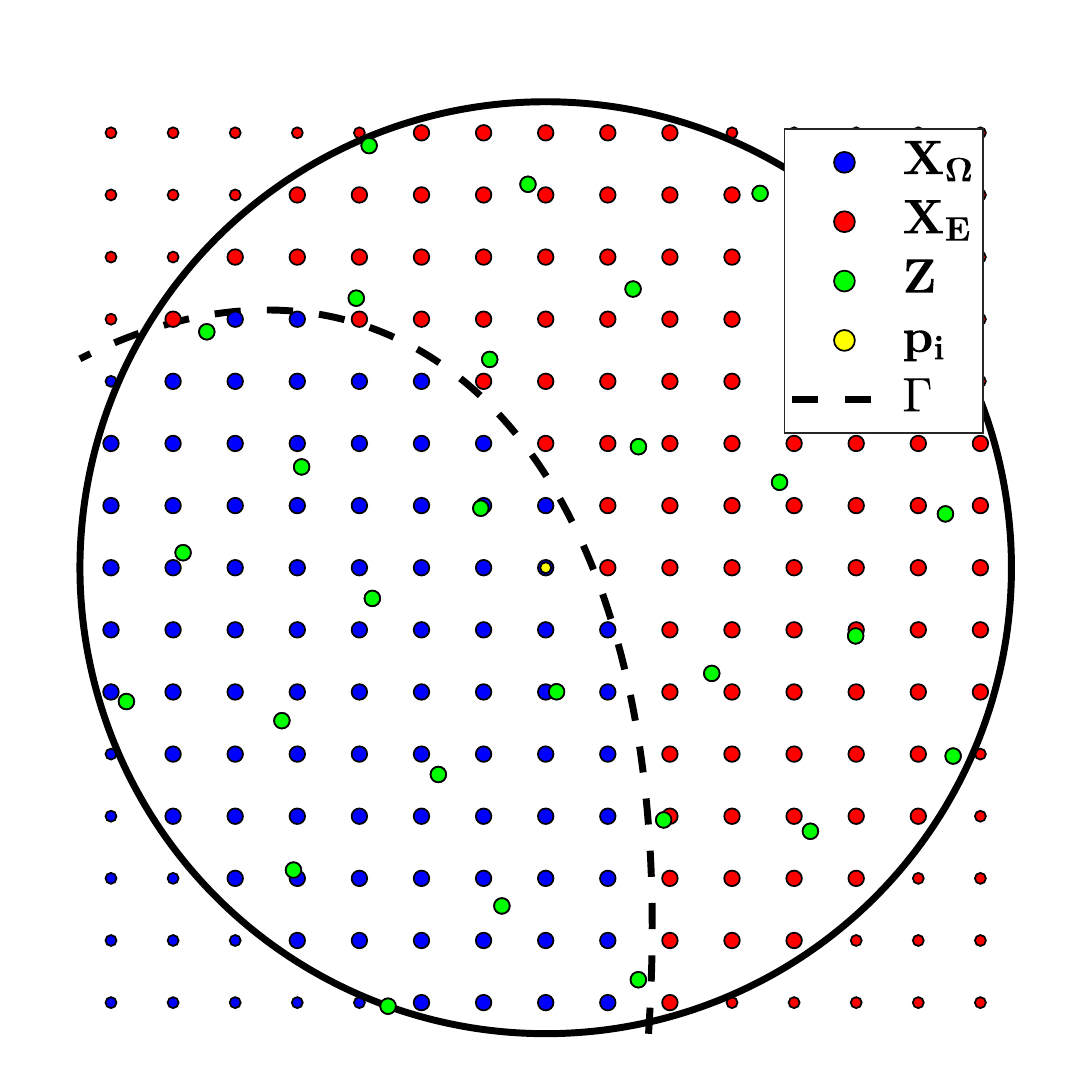}};
\draw [ultra thick,black,->] (2.0,3.5) to[bend left] (7.0,3.5);
\end{tikzpicture}
    \end{center}
    \caption{Left: Schematic figure of distribution of extension partitions along $\Gamma$ for a complex domain. The green markers correspond to RBF centres $\mathbf{Z}$, generated by \eqref{eq:Vogel},  and the distribution is repeated for every partition. Right: classification of points in $\mathbf{X}$ as inside or outside $\Omega$. The larger markers denote points in $\mathbf{X}_{i,\Omega}$ and  $\mathbf{X}_{i,E}$.}
    \label{fig:partition}
\end{figure}

Let $\{\psi^{k}_{i}\}_{i = 1}^{N_{\psi}}$ be a collection of $N_{\psi}$ compactly supported radial basis functions such that
$\psi^{k}_{i}(\x) = \psi^{k}(\x-\mathbf{p}_{i})$ for some choice of centres $\{\mathbf{p}_{i}\}_{i = 1}^{N_{\psi}}$. The superscript $k$ indicates the smallest subset $C^{k}_{0}$ of $C_{0}$ that $\psi^{k}$ is a member of. Define a partition $\Omega_{i}$ as the support of $\psi^{k}_{i}$, i.e.
$\Omega_{i} = \text{supp}(\psi^{k}_{i})$, which is a disc with radius $R$. We will return to the choice of $\psi^{k}$ in Section \ref{sss:PUXprop}. Note that all partitions have the same radius. The number of partitions $N_{\psi}$, the location of the partition centres $\{\mathbf{p}_{i}\}_{i = 1}^{N_{\psi}}$ and radius $R$ are chosen such that the partitions cover $\Gamma$ and that the partitions overlap with approximately a radius. The following notation will be useful. Each partition $\Omega_{i}$ has a set of points on the uniform grid within $R$ of
$\mathbf{p}_{i}$, which we denote $\mathbf{X}_{i}$, rather than $\mathbf{X}_{\Omega_{i}}$. It can be split into two disjoint subsets: $\mathbf{X}_{i,\Omega} = \{\mathbf{x}\in\Omega_{i}\cap \Omega\}$ and
$\mathbf{X}_{i,E} = \{\mathbf{x}\in\Omega_{i}\cap E\}$. Let $N_{i}$ denote the number of elements in $\X_{i}$. Analogously, let $N_{i,\Omega}$ and $N_{i,E}$ denote the number of elements in $\mathbf{X}_{i,\Omega}$ and $\mathbf{X}_{i,E}$, respectively. See Figure \ref{fig:partition} for a graphical example. Given a function $f:\Omega \rightarrow \mathbb{R}$, the function values at the locations $\mathbf{X}_{i,\Omega}$ are used to create a local extension $f^{e}_i$. We will return to the construction of the local extensions in Section \ref{sss:localinterp}, but for now assume their existence.

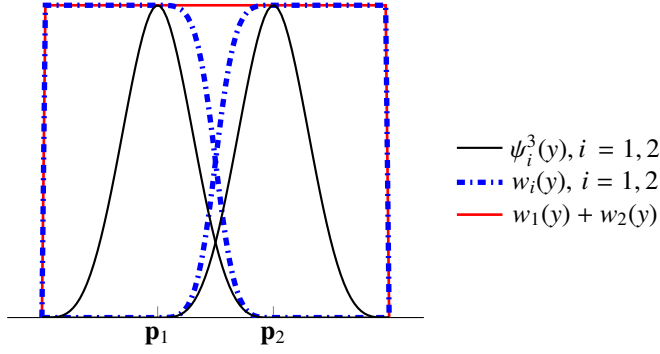
\begin{figure}
\begin{center}
\begin{tikzpicture}
\begin{axis}[
xticklabels={$\mathbf{p}_{1}$,$\mathbf{p}_{2}$},
reverse legend,
xtick={-1,0},
legend style={draw=none},
hide y axis,
axis x line*=bottom,
ymin=0,
no markers,
every axis plot/.append style={thick},
legend style={at={(1.6,0.4)},anchor=east},
scale=0.8,
]
\addplot[domain=-2:1,red,solid,line width=1.0pt, samples=101]{((max(1-abs(\x),0)^4)*(1/4)*(4+ 16*abs(\x)+12*abs(\x)^2+3*abs(\x)^3))/((max(1-abs(\x+1),0)^4)*(1/4)*(4+ 16*abs(\x+1)+12*abs(\x+1)^2+3*abs(\x+1)^3) + (max(1-abs(\x),0)^4)*(1/4)*(4+ 16*abs(\x)+12*abs(\x)^2+3*abs(\x)^3))+     ((max(1-abs(\x+1),0)^4)*(1/4)*(4+ 16*abs(\x+1)+12*abs(\x+1)^2+3*abs(\x+1)^3))/((max(1-abs(\x+1),0)^4)*(1/4)*(4+ 16*abs(\x+1)+12*abs(\x+1)^2+3*abs(\x+1)^3) + (max(1-abs(\x),0)^4)*(1/4)*(4+ 16*abs(\x)+12*abs(\x)^2+3*abs(\x)^3))};

\addplot[domain=-2:1,blue,dashdotted,line width=2.0pt, samples=101]{((max(1-abs(\x),0)^4)*(1/4)*(4+ 16*abs(\x)+12*abs(\x)^2+3*abs(\x)^3))/((max(1-abs(\x+1),0)^4)*(1/4)*(4+ 16*abs(\x+1)+12*abs(\x+1)^2+3*abs(\x+1)^3) + (max(1-abs(\x),0)^4)*(1/4)*(4+ 16*abs(\x)+12*abs(\x)^2+3*abs(\x)^3))};

\addplot[domain=-2:1,blue,dashdotted,line width=2.0pt,forget plot, samples=101]{((max(1-abs(\x+1),0)^4)*(1/4)*(4+ 16*abs(\x+1)+12*abs(\x+1)^2+3*abs(\x+1)^3))/((max(1-abs(\x+1),0)^4)*(1/4)*(4+ 16*abs(\x+1)+12*abs(\x+1)^2+3*abs(\x+1)^3) + (max(1-abs(\x),0)^4)*(1/4)*(4+ 16*abs(\x)+12*abs(\x)^2+3*abs(\x)^3))};

\addplot[forget plot, domain=-2:1, samples=101]{(max(1-abs(\x+1),0)^4)*(1/4)*(4+ 16*abs(\x+1)+12*abs(\x+1)^2+3*abs(\x+1)^3)};

\addplot[domain=-2:1, samples=101]{(max(1-abs(\x),0)^4)*(1/4)*(4+ 16*abs(\x)+12*abs(\x)^2+3*abs(\x)^3)};

      \legend{$w_{1}(y)+w_{2}(y)$\\$w_{i}(y),\,i=1,2$\\$\psi^{3}_{i}(y),i\,=1,2$\\}
      \end{axis}
    \end{tikzpicture}
\end{center}
\caption{Plot of weight functions \eqref{eq:weightfunction} and their sum.}
      \label{fig:weight}
\end{figure}

\begin{figure}[ht]
    \begin{center}

\begin{tikzpicture}
\node (bad) {\includegraphics[trim=2.4cm 1.7cm 1.4cm 1cm, clip, width=0.3\textwidth]{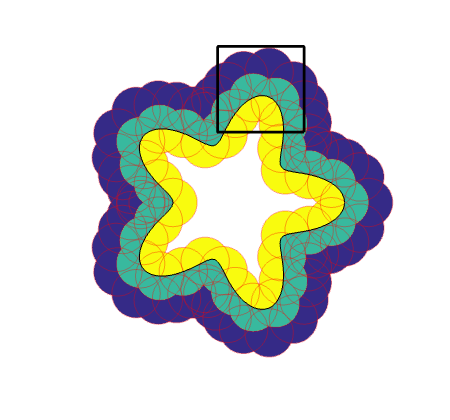}};
\node (good) [right=of bad] {\includegraphics[trim=1.8cm 1.2cm 1.2cm 0.8cm, clip, width =0.3\textwidth]{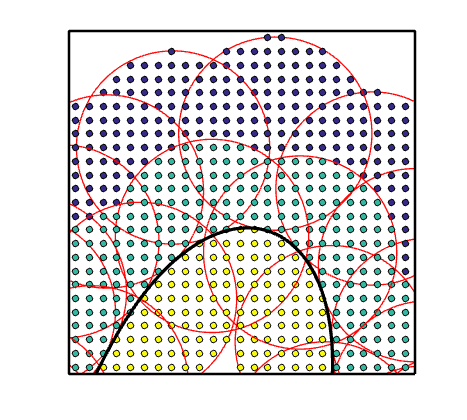}};
\draw [ultra thick,black,->] (1.2,2.5) to[bend left] (3.6,2.5);
\end{tikzpicture}

    \end{center}
\caption{Schematic image for function extension from a star shaped domain $\Omega$ given by the black border. Observe that in this figure the partitions are not centred at uniform grid points. The red overlapping circles are the partitions. The yellow section corresponds to the uniform data  points used for creating the local extension $f^{e}_{i}$, the blue section to points where $f^{e} = 0$ and the green sector is a blend of the two.}
\label{fig:schematic_extension_supression}

\end{figure}

For every partition $\Omega_{i}$ and its associated radial basis function $\psi^{k}_{i}$ define the corresponding weight function $w_{i}$ as
\begin{equation}
\label{eq:weightfunction}
w_{i}(\mathbf{x}) = \frac{\psi^{k}_{i}(\mathbf{x})}{\sum\limits_{j=1}^{N_{\psi}}\psi^{k}_{j}(\mathbf{x})},
\end{equation}
which belongs to the space $C^{k}_{0}$. By construction the set of weights $\{w_{i}\}_{i = 1}^{N_{\psi}}$ forms a partition of unity. That is
\begin{equation}
\label{eq:POU}
\sum\limits_{i=1}^{N_{\psi}}w_{i}(\mathbf{x})= 1,\quad \forall\mathbf{x}\in \bigcup\limits_{i = 1}^{N_{\psi}}\bar{\Omega}_{i},
\end{equation}
which is referred to in the literature as Shepard's method \cite{shepard}. See Figure \ref{fig:weight} for a visualisation. This construction is used to combine the local extension $\{f^{e}_{i}\}_{i = 1}^{N_{\psi}}$ into a global one,
\begin{equation}
f^{e}(\mathbf{x}) = \sum\limits_{i=1}^{N_{\psi}}w_{i}(\mathbf{x})f^{e}_{i}(\mathbf{x}).
\label{eq:fe_notcompact}
\end{equation}
However, \eqref{eq:fe_notcompact} is not used, as we want an extension that it is continuous or of higher regularity as it is extended by zero outside its support. Refer to the set of partitions $\{\Omega_{i}\}_{i = 1}^{N_{\psi}}$ as \emph{extension} partitions and now introduce also the \emph{zero} partitions
$\{\Omega^{0}_{i}\}_{i = 1}^{N_{\psi}^{0}}$. They are included in the partition of unity definition \eqref{eq:POU} and distributed such that they overlap the extension partitions, but do not intersect $\bar{\Omega}$. The associated local extension $f^{e}_i$ is set to be identically equal to zero for $i = 1,\ldots, N^{0}_{\psi}$. Hence, as the zero partitions are blended with the local extensions in the first layer of partitions $\{\Omega_{i}\}_{i = 1}^{N_{\psi}}$, the global extension will be forced to zero over the overlapping region. Therefore zero partitions should be placed such that
 $f^{e}$ has a controlled decay to zero and that the size of the overlap with extension partitions are about the same, see Figure \ref{fig:schematic_extension_supression}.  Thus the global extension will in these parts have the same regularity as $w^{p}$, as given by the regularity of the compactly supported radial basis function  $\psi^{k}$. The extension $f^{e}$ of $f$ is given by
\begin{equation}
  \label{eq:functionextension}
      f^{e}(\mathbf{x}) = \begin{cases}f(\mathbf{x}),\quad \mathbf{x}\in\mathbf{X}_{\Omega},\\\sum\limits_{i = 1}^{N_{\psi}+N_{\psi}^{0}}w_{i}(\mathbf{x})f^{e}_{i}(\mathbf{x}),\quad \mathbf{x} \in \bigcup\limits_{i=1}^{N_{\psi}}\mathbf{X}_{i, E},\\
      0, \quad \text{otherwise}.
      \end{cases}
\end{equation}

As $\psi^{k}$ we use one of the compactly supported Wu-functions, which are tabulated after their regularity $k$, see Table \ref{tab:Wu} or \cite{Fasshauer:2007:MAM:1506263}. There are other options, but the Wu-functions have compact support and are simple to implement. Note that they have lower regularity at the origin, e.g. the Wu-function listed as $C^{4}$ is only $C^{2}$ at that point. Moreover, the $k+1$th derivative of $\psi^{k}$ is of bounded variation. The partitions centres $\{\mathbf{p}_{i}\}_{i = 1}^{N_{\psi}}$ are set to be nodes on the regular grid that are the closets to be boundary, yet still in $\X_{i,\Omega}$.
 Thus evaluation of weight functions at the origin is omitted and higher regularity is maintained. With this, we have described how local extensions are combined into a global one. It remains to construct the local extensions $\{f^{e}_{i}\}_{i = 1}^{N_{\psi}}$.

\begin{table}
\centering
{\def\arraystretch{1.3}
\begin{tabular}{cc}
\toprule
\textbf{Regularity} & $\psi^{k}(r)$ \\
\midrule
$\psi^{1}\in C^{1}$ & $(1- r)_{+}^{2}(2+ r)$ \\
$\psi^{2}\in C^{2}$ & $(1- r)_{+}^{3}(8+9 r+3r^{2})$ \\
$\psi^{3}\in C^{3}$ & $(1-r)_{+}^{4}(4+16r+12r^{2}+3r^{3})$ \\
$\psi^{4}\in C^{4}$ & $(1-r)_{+}^{5}(8+40r+48r^{2}+25r^{3}+5r^{4})$ \\
$\psi^{5}\in C^{5}$ & $(1-r)_{+}^{6}(6+36r+82r^{2}+72r^{3}+30r^{4}+5r^{5})$ \\
\bottomrule
\end{tabular}
}
\caption{Wu-functions $\psi^{k}\in C^{k}_{0}$, with compact support in $r\in (0,1)$ \cite{Fasshauer:2007:MAM:1506263}. Here $(\cdot)_{+} = \max{(0,\cdot)}$. The listed regularity excludes evaluation at the origin.}
\label{tab:Wu}
\end{table}

\subsubsection{Local extensions}
\label{sss:localinterp}
We now return to the construction of the local extensions $f^{e}_{i}$ for each extension partition $i = 1,\ldots,N_{\psi}$. The local extension $f^{e}_{i}$ is created as a weighted sum of radial basis functions, that interpolates the values of $f$ at $\x\in \mathbf{X}_{i,\Omega}$ and is evaluated at $\x\in \mathbf{X}_{i,E}$. The radial basis functions are denoted
$\phi_{j}(\x) = \phi(\|\mathbf{z}_{j}-\x\|)$, to distinguish them from the radial basis functions $\psi^{k}$. The elements of the set $\mathbf{Z} = \myset{\mathbf{z}}{j}{N_{\phi}}\subset \text{supp}(\psi^{k}_{i})$ are the centres for the RBFs, whose distribution for now is left unspecified.
The standard form of an RBF interpolant at a point $\mathbf{x}$ is
\begin{equation}
\label{eq:sfOmegai}
f^{e}_{i}(\mathbf{x}) = \sum\limits_{j=1}^{N_{\phi}}\lambda_{j}\phi(\|\mathbf{z}_{j}-\mathbf{x}\|)
\end{equation}
where $\lambda_{j}$ are unknown coefficients to be determined. We use
 \begin{equation}
 \phi(\|\mathbf{z}_{j}-\mathbf{x}\|)  = e^{-(\varepsilon \|\mathbf{z}_{j}-\mathbf{x}\|)^{2}},
\label{eq:GaussianRBF}
\end{equation}
 where  $\varepsilon$ is a shape parameter setting the width of the Gaussian. The smallest interpolation error is obtained when $\varepsilon$ is small, yet nonzero, but no general value can be given \cite{Larsson2005}.

 With some abuse of notation let $\mathbf{X_{i}}$ and $\mathbf{Z}$ refer to vectors with the members of respective set as elements. Then, following the outline of \citep{Larsson1077950}, let $\Phi(\mathbf{X}_{i},\mathbf{Z})$ denote an $N_{i}\times N_{\phi}$-matrix  with elements
 $\Phi(\X_{i},\mathbf{Z})_{m,n}= \phi(\|\x_{m}-\mathbf{z}_{n}\|)$, for $m = 1,\ldots,N_{i}$ and $n = 1,\ldots,N_{\phi}$. Furthermore, let $\Lambda = (\lambda_{1}\,\lambda_{2}\,\ldots\, \lambda_{N_{\phi}})^{T}$.
  Consider a scenario when $f$ is known for all nodes in $\X_{i}$, then the associated interpolation problem to \eqref{eq:sfOmegai} can be written as
 \begin{equation}
 \label{eq:sfOmegaiMatrix}
 \Phi(\mathbf{X}_{i},\mathbf{Z})\Lambda = f_{\X_{i}},
 \end{equation}
with $f_{\X_{i}} = f(\X_{i})$. If $N_{i}\geq N_{\phi}$ then  $\Lambda$ can be solved for in a least-squares sense. However, this is an unstable problem for several reasons. First, the conditioning of the problem is heavily dependent on the shape parameter $\varepsilon$.  For small $\varepsilon$ the interpolation weights $\Lambda$ oscillate between positive and negative numbers of large magnitude \cite{Larsson2005}. Furthermore, it is not uncommon for the condition number for the interpolation matrix to be of order $10^{18}$ or more. These characteristics are common for interpolation with radial basis functions. Additionally, the data is represented on a uniform grid; collocating at these locations is the worst possible setting for interpolation, as with polynomials. These shortcomings can be circumvented by avoiding collocation and considering a least squares problem instead. Note that all problems mentioned above are purely numerical artifacts. The function space spanned by Gaussians is indeed a good approximation space.

 Decouple the centres $\mathbf{Z}$ of the radial basis functions from $\mathbf{X}_{i}$ and assume they are distributed in a near optimal way with respect to minimising the interpolation error. We wish to omit explicit use of the interpolation coefficients $\Lambda$ in \eqref{eq:sfOmegaiMatrix}. It can be achieved by formally solving for $\Lambda$ by collocating at the centres $\mathbf{Z}$:
\begin{equation}
\Phi(\mathbf{Z},\mathbf{Z})\Lambda = f_{\mathbf{Z}}\Leftrightarrow \Lambda = \Phi(\mathbf{Z},\mathbf{Z})^{-1}f_{\mathbf{Z}}.
\end{equation}
Here $f_{\mathbf{Z}}$ are the values of $f$ at the locations $\mathbf{Z}$, which are unknown. Due to the choice \eqref{eq:GaussianRBF} the matrix $\Phi(\mathbf{Z},\mathbf{Z})$ is symmetric and positive definite, thus the inverse $\Phi(\mathbf{Z},\mathbf{Z})^{-1}$ is well-defined. We can now reformulate \eqref{eq:sfOmegaiMatrix} as

\begin{equation}
 \label{eq:sfOmegaiMatrixReform}
 \Phi(\mathbf{X}_{i},\mathbf{Z})\Phi(\mathbf{Z},\mathbf{Z})^{-1}f_{\mathbf{Z}} = f_{\mathbf{X}_{i}}.
 \end{equation}
Henceforth we use the shorthand notation $A(\mathbf{X}_{i},\mathbf{Z}) = \Phi(\mathbf{X}_{i},\mathbf{Z})\Phi(\mathbf{Z},\mathbf{Z})^{-1}$. For the purpose of function extension, sort the data points in $\Omega_{i}$ such that

\begin{equation}
\mathbf{X}_{i} = \begin{bmatrix} \mathbf{X}_{i,\Omega}\\\mathbf{X}_{i,E}\end{bmatrix}.
\end{equation}
where the components are of length  $N_{i,\Omega}$ and $N_{i,E}$, respectively. Consequently, $A$ can also be rearranged and split into two block matrices
\begin{equation}
\label{eq:Amatrix}
A = \begin{bmatrix}
A_{i,\Omega}\\
A_{i,E}
\end{bmatrix},
\end{equation}
with $A_{i,\Omega} = \Phi(\mathbf{X}_{i,\Omega},\mathbf{Z})\Phi(\mathbf{Z},\mathbf{Z})^{-1}$ of size $N_{i,\Omega}\times N_{\phi}$ and
$A_{i,E}= \Phi(\mathbf{X}_{i,E},\mathbf{Z})\Phi(\mathbf{Z},\mathbf{Z})^{-1}$ of size $N_{i,E}\times N_{\phi}$. Since $f$ is  known at
$\mathbf{X}_{i,\Omega}$
it can replace the corresponding entries in $f_{\mathbf{X}_{i}}$ \eqref{eq:sfOmegaiMatrixReform} with $f_{i,\Omega} = f(\mathbf{X}_{i,\Omega})$. For each partition we obtain the system,
 \begin{equation}
   \label{eq:LS}
 \begin{bmatrix}
 A_{i,\Omega}\\
 A_{i,E}
 \end{bmatrix}
 f_{\mathbf{Z}} = \begin{bmatrix}
 f_{i,\Omega}\\
 f_{i,E}
 \end{bmatrix},
 \end{equation}
with $f_{i,E} = f(\mathbf{X}_{i,E})$ unknown. For each partition $i$ the values $f_{i,\Omega}$ are mapped to the nodes $\mathbf{Z}$ to obtain $f_{\mathbf{Z}}$. Thereafter we obtain $f_{i,E}$, which is the local extension. That is:
\begin{enumerate}
\item Solve the least-squares problem $A_{i,\Omega}f_{\mathbf{Z}}=f_{i,\Omega}$ for $f_{\mathbf{Z}}$.
\item Obtain the local extension $f^{e}_{i}(\mathbf{X}_{i,E}) = f_{i,E} = A_{i,E}f_{\mathbf{Z}}$.
\end{enumerate}

This approach allows us to use a non-uniform distribution of RBF centres which significantly improves the stability, but still lets the data be represented on the uniform grid. We also avoid explicit use of the interpolation weights $\Lambda$. It remains to address the notorious ill-conditioning of $\Phi$, associated with the shape parameter $\varepsilon$ set small. This is achieved by applying the algorithm \emph{RBF-QR}. It is intended for a formulation as \eqref{eq:LS}, since it computes $A$, rather than $\Phi^{-1}$, which acts as a mapping of data from non-uniformly  to uniformly distributed locations. Said algorithm performs a change of basis for $A$, and in process the condition number is reduced, see \cite{MR2801193}. By the use of RBF-QR the restrictions of choosing $\varepsilon$ is lifted.

\subsubsection{Properties of PUX}
\label{sss:PUXprop}
Four parameters need to be set for the PUX algorithm: the shape parameter $\varepsilon$ for the width of Gaussians \eqref{eq:GaussianRBF} used as interpolation basis, the partition radius $R$, the length $L$ for the computational domain $B = [-L,L]^2$ and $N_u$, where $N_{u}^{2}$ is the number of uniformly distributed nodes over $B$. The remaining parameters can be set based on these values. Here we give the most important relations. For a complete discussion see \cite{FRYKLUNDPUX}.

Due to RBF-QR the shape parameter can be set small without risk of suffering from ill-conditioning. A good value is $\varepsilon = 2$, but the error in solving the modified Helmholtz equation is relatively insensitive.

Let $P$ be
\emph{number of uniform grid points per partition radius}, denoted as
\begin{equation}
 \label{eq:Pdef}
 P = \frac{N_{u}}{2L}R.
\end{equation}
This measure is used to choose $\psi^{k}$ from Table \ref{tab:Wu}, and the number $N_{\phi}$ of basis functions \eqref{eq:GaussianRBF} per partition. To see how $P$ relates to $\psi^{k}$, consider the convergence of the error in solving the modified Helmholtz equation \eqref{eq:ModHelmEq}--\eqref{eq:ModHelmEqBC}, assuming that only resolving $u^{P}$
 limits the accuracy. If $f^e$ is smooth then the error has asymptotically spectral convergence. However, the extension inherits the regularity of the weight function $w$. Recall that by construction $w\in C^{k}_{0}$ for a fixed $\psi^{k}$  \eqref{eq:weightfunction}. Consequently the error has an asymptotic convergence of $4+k$, if the $k$th derivative of $f^{e}$ is of bounded variation. A Wu-function of high regularity is harder to resolve than one of lower regularity. This implies that given a resolution $P$ the error in resolving the Wu-function may hamper the convergence. As in \cite{FRYKLUNDPUX} we use the heuristic relation
\begin{equation}
  \label{eq:Ck}
  k = \min\left(\left\lfloor\sqrt{P}-0.9\right\rfloor,5\right)
\end{equation}
for choosing $\psi^{k}$. In Section \ref{s:results} we confirm that \eqref{eq:Ck} is a satisfactory estimate for an optimal $\psi^{k}$ given $P$.

Creating a local extension involves solving the least-squares problem $A_{i,\Omega}f_{\mathbf{Z}}=f_{i,\Omega}$ for $f_{\mathbf{Z}}$ for some $i$. It should be sufficiently overdetermined in order to be a well-posed problem. Given a $P$ the number of unknowns $N_{\phi}$ should be set accordingly to obtain a certain ratio of knowns and unknowns. Still, $P$ can be of such magnitude that $N_{\phi}$ is larger than required to obtain good results and the least-squares problem is more stable and cheaper to solve if the unknowns are few. Thus if the available data is abundant it can be downsampled to reduce $P$, and therefore $N_{\phi}$. Let $c$ be the sampling parameter, defined as
\begin{equation}
  \label{eq:coarsen}
c = \min\left(\left\lfloor\sqrt{\frac{P}{8}}\right\rfloor,1\right).
\end{equation}
 If $c = 1$ then all points are used, $c = 2$ means that every other point is removed, etc. Then, as in \cite{FRYKLUNDPUX}, we use
\begin{equation}
  \label{eq:setM}
  N_{\phi} = \left\lfloor\min\left(0.8\pi (P/c)^{2}/4,3(P/c)\right)\right\rfloor
\end{equation}
to set the number of radial basis functions per partition.
Note that choosing $\psi^{k}$ is question about resolution; Wu-functions of higher regularity require larger $P$ to be well resolved, while setting $N_{\phi}$ is related to solving a least-squares problem. These are two separate problems and two different values for $P$ may be used. So given a $P$ we set $\psi^{k}$ according to \eqref{eq:Ck} and then compute $c$ with \eqref{eq:coarsen}. Now $N_{\phi}$ is set by \eqref{eq:setM} for $P/c$. Thus the local least-squares problems are solved on a potentially coarser grid, but the local extensions are on the original grid.

The distribution of RBF-centres $\mathbf{Z}$ can be chosen freely, and we use the quasi uniform \emph{Vogel} node distribution defined as
\begin{equation}
  \label{eq:Vogel}
  \mathbf{z}_j = \sqrt{\frac{j}{N_{\phi}}}\left(\cos{\left(j\pi\left(3-\sqrt{5}\right)\right)},\sin{\left(j\pi\left(3-\sqrt{5}\right)\right)}\right),     \quad j = 1,\ldots,N_{\phi},
\end{equation}
in a unit disc. See Figure \ref{fig:partition} for a visualisation. The distribution \eqref{eq:Vogel} is near optimal and RBF-QR performs well up to about $400$ nodes. The locality of the weight functions guarantees that the least squares systems are of moderate size, which can be solved in parallel.

Constructing $A$ \eqref{eq:Amatrix} with RBF-QR is a computationally expensive operation, so employing it for every partition is undesirable. However, the matrix is the same for all partitions since $\mathbf{p}_{i}$ is centred at a grid point from the uniform distribution. Thus the pairwise distances for the elements in $\X_{i}$ are independent of $i$. Therefore a single matrix $A$ can be precomputed with RBF-QR and reused for all extension partitions. The only difference between them in terms of $A$ is the decomposition of $\mathbf{X}_i$ into $\mathbf{X}_{i,\Omega}$ and  $\mathbf{X}_{i,E}$,
 as it depends on how the boundary $\Gamma$ intersects the partition. Note that the zero partitions may individually have a radius different from $R$ in order to conform to the geometry of $\Omega$ and to overlap the extension partitions properly.
\subsection{The homogeneous problem}
\label{sss:discretisation_adaptive_modhelm_homo}
For simplicity, assume the number of countours $N_{\Gamma}$ to be one and write $\Gamma_{n,k} = \Gamma_{k}$, $\bm{y}_{n,k} = \bm{y}_{k}$ and $s_{n,k} = s_{k}$. We apply an $N_{Q}$-point, panel-based Nystr\"{o}m discretisation scheme based on
the composite Gauss-Legendre quadrature rule, with nodes $t_m^{G}$ and weights $W_m^{G}$, with
 $m = 1,\ldots,N_Q$. Let
 $\y_{k,m} = \bm{\gamma}_{k}(t_m^{G})$, $s_{k,m} = s_{k}(t_m^{G})$ and $\mu_{k,m} = \mu_{k}(t_m^{ G})$. An approximation of the solution $\mu$ to \eqref{eq:semidischomo}  is the solution of
\begin{equation}
  \label{eq:DensityDisc}
\mu_{i,j} + \sum\limits_{\substack{k=1}}^{N_{P}}\sum\limits_{\substack{m=1}}^{N_{Q}}M\left(\x_{i},\y_{k,m}\right)\mu_{k,m}s_{k,m}W^{G}_{m} = - \frac{2\tilde{g}_{i,j}}{\alpha^{2}},\quad i = 1,\ldots,N_{P},\quad i = j,\ldots,N_{Q}.
\end{equation} and correspondingly for \eqref{eq:dblPotpan} we have

\begin{equation}
  \label{eq:dblPotMDisc}
  u^{H}(\x) = \sum\limits_{\substack{k=1}}^{N_{P}}\sum\limits_{\substack{m=1}}^{N_{Q}}M\left(\x,\y_{k,m}\right)\mu_{k,m}s_{k,m}W^{G}_{k,m}, \quad \x \in \Omega.
\end{equation}
An important observation is that the kernel $M$ \eqref{eq:kerM} is not smooth and can contain singularities, depending on how $\x$ approaches $\y$. Here the Gauss-Legendre quadrature rule is insufficient, as the resulting loss of accuracy can be critical enough to render the result useless. We elaborate on this topic in Section \eqref{sss:discretisation_adaptive_modhelm_homo_specquad}.

In matrix notation \eqref{eq:DensityDisc} can be written as $(\bm{I} + \bm{M})\bm{\mu}=\mathbf{\tilde{g}}$, where $\bm{I}$ is the identity matrix and $\bm{M}$ a compact operator. The density $\bm{\mu}$ can be efficiently obtained with GMRES, in terms of numbers of iterations. The condition number for $\bm{I} + \bm{M}$ is typically small or moderate and uniformly bounded.
A fast multipole method (FMM) can be used for efficient computation of the involved potentials in \eqref{eq:DensityDisc} and \eqref{eq:dblPotMDisc} \cite{doi:10.1137/0909044}. We use the point to point FMM for the two-dimensional Yukawa kernel presented \cite{Kropinski2011modHelm}. It is based  on the volume equivalent in \cite{CHENG2006616}. For the corresponding three-dimensional version see \cite{GREENGARD2002642}.

Finally a note on the restriction of the boundaries being smooth. For non-smooth boundaries the integrand of  \eqref{eq:semidischomo}
 is not compact and the Fredholm alternative fails. While there are theoretical results on the solvability with Lipschitz continuous boundaries \cite{VERCHOTA1984572}, they require the implementation of sophisticated quadrature techniques, such as \cite{2012arXiv1207.6737H}, which we have not implemented. These methods also allow cusps, i.e. non-Lipschitz boundaries, and mixed boundary conditions.

\subsubsection{Special purpose quadrature}
\label{sss:discretisation_adaptive_modhelm_homo_specquad}
When solving for $\mu$ in \eqref{eq:DensityDisc} or evaluating the layer potential \eqref{eq:dblPotMDisc} several orders of accuracy may be lost, since the kernel $M$ \eqref{eq:kerM} is not smooth. Moreover, $M$ can be singular, depending on if $\x$ approaches some $\y\in\Gamma$ along $\Gamma$ or from $\Omega$.
One of the most efficient methods to circumvent this loss of accuracy is explicit kernel-split quadrature with product integration by \citeauthor{HELSING20098892}, see \cite{HELSING20098892}. However, for the modified Helmholtz equation with \emph{large} $\alpha$,  i.e. for high temporal resolution, it can fail completely. Below we sketch the problem, its relation to $\alpha$ and how to circumvent it.

We start by explaining product integration, which requires the involved integrals to be expressed in complex notation. To keep these paragraphs brief and simple, the reformulations are omitted. Consider a single panel $\Gamma_{k}\in \mathbb{C}$ with endpoints at $-1$ and $1$, but the panel does not have to follow the real axis. Let $\varphi: \Gamma_{k} \rightarrow \mathbb{R}$ be a smooth function and $s:\Gamma_{k} \times \mathbb{C} \rightarrow \mathbb{R}$ a non-smooth kernel that may be singular or nearly singular. The goal is to compute
\begin{equation}
\label{eq:specquadtarget}
    \int_{\Gamma_{k}} \varphi(\tau) s(\tau_{0},\tau) \,d\tau
\end{equation}
accurately for some fixed $\tau_0\in\mathbb{C}$ arbitrarily close to or on $\Gamma_k$. To do this, approximate $\varphi$ with a polynomial of degree $N_{Q}-1$ in $\tau\in \Gamma_{k}$, such that
\begin{equation}
    \label{eq:polyapprox}
    \varphi(\tau_{0}) \approx \sum_{n=1}^{N_{Q}}c_{n}  \tau^{n-1},
\end{equation}
with unknown coefficients $\{c_n\}$. Inserting this into \eqref{eq:specquadtarget} gives
\begin{equation}
\label{eq:specQuad}
    \int_{\Gamma_{k}} \varphi(\tau) s(\tau_{0},\tau) \,d\tau\approx \sum_{n=1}^{N_{Q}}c_{n} \int_{\Gamma_k} \tau^{n-1} s(\tau_{0},\tau) \,d\tau.
\end{equation}
The integrals on the right hand side can can be computed analytically through recursive formulas. The unknown coefficients $\{c_n\}$ are obtained by solving a Vandermonde system. If $\varphi$ can be accurately represented as a $N_{Q}-1$ degree polynomial over $\Gamma_{k}$, then product integration allows evaluation of integrals such as \eqref{eq:specQuad} without loss of accuracy as $\tau_{0}$ and $\tau$ approach each other.

Kernel-split means that a kernel is decomposed into smooth and singular terms. Leaving complex notation, by \cite[\S10]{NIST:DLMF} the first-order modified Bessel function of the second kind $K_{1}$, appearing in \eqref{eq:kerM}, can be decomposed as
\begin{equation}
  \label{eq:K1decomp}
K_{1}\left(x\right) = \frac{1}{x} + I_{1}\left(x\right)\log\left(x\right) + K^{S}_{1}\left(x\right),\quad x\in \mathbb{R}^{+}.
\end{equation}
This form is attractive since the singular terms are separated and can be studied individually. Here $I_{1}$ is the modified Bessel function of the first kind of order one and $K^{S}_{1}$ is a power series in $x$. For the kernel $M$, see \eqref{eq:kerM}, the situation is slightly more involved, as the singularity structure depends on how $\x$ approaches $\y\in \Gamma$.  To distinguish between the two cases, for any $\y \in \Gamma$ denote $M(\x,\y)$ as $M_{\Gamma}(\x,\y)$
for $\x\in \Gamma$ and $M_{\Omega}(\x,\y)$ for $\x\in  \Omega$. We first study $M_{\Gamma}$; the decomposition \eqref{eq:K1decomp} motivates the formulation
  \begin{equation}
  \label{eq:Mpartial}
  M_{\Gamma}(\mathbf{x},\mathbf{y}) = M_{\Gamma,0}(\mathbf{x},\mathbf{y})+\log(\|\y-\x\|)M_{\Gamma,L}(\mathbf{x},\mathbf{y}),\quad\mathbf{x},\mathbf{y}\in\Gamma,
  \end{equation}
  with $M_{\Gamma,L}$ identified as
  \begin{equation}
  \label{eq:MpartialL}
  M_{\Gamma,L}(\mathbf{x},\mathbf{y}) = \frac{\alpha}{\pi}I_{1}\left(\alpha\|\bf{y}-\bf{x}\|\right)\frac{\mathbf{y}-\mathbf{x}}{\|\mathbf{y}-\mathbf{x}\|}\cdot \nu_{\mathbf{y}},\quad\mathbf{x},\mathbf{y}\in\Gamma.
  \end{equation}
The term $M_{\Gamma,0}$ is smooth and by \eqref{eq:Mlimit} we have
  \begin{equation}
  \label{eq:MGamma0}
  M_{\Gamma,0}(\mathbf{y},\mathbf{y}) = -\frac{1}{2\pi}\kappa(\mathbf{y}),\quad\mathbf{y}\in\Gamma,
  \end{equation}
    since the term $\log(\|\x-\y\|)M_{\Gamma,L}(\x,\y)$ goes to zero in the limit $\x\rightarrow \y$. But in this limit the derivative of $\log(\|\y-\x\|)M_{\Gamma,L}(\mathbf{x},\mathbf{y})$ has a log-type singularity. Thus standard quadrature rules that relies on smoothness fail to be accurate. To maintain accuracy product integration is needed, even though the limit is well-defined. In terms of \eqref{eq:specQuad} $\phi$ and $s$ correspond to $\mu M_{\Gamma,L}$ and $\log$. This approach is used to compute the involved integrals in \eqref{eq:semidischomo}.

In the case $\bf{x}\in \Omega$, corresponding to computing \eqref{eq:dblPotpan}, the kernel $M(\x,\y)$ is singular in the limit $\x\rightarrow\y$ and product integration is required. We have
\begin{equation}
\label{eq:Momega}
M_{\Omega}(\mathbf{x},\mathbf{y}) = M_{\Omega,0}(\mathbf{x},\mathbf{y})+\log(\|\y-\x\|)M_{\Omega,L}(\mathbf{x},\mathbf{y})+\frac{(\y-\x)\cdot\nu_{\bf{y}}}{\|\y-\x\|^{2}}M_{\Omega,C}(\mathbf{x},\mathbf{y}),\quad\mathbf{x}\in \Omega,\mathbf{y}\in\Gamma,
\end{equation}
where $M_{\Omega,0}$ is a smooth function, $M_{\Omega,L} = -\alpha^{2}/2M_{\Gamma,L}$ and $M_{\Omega,C} = -\alpha^{2}/2$. Again, we identify $\varphi$ from \eqref{eq:specQuad} as $\mu$ multiplied with $M_{\Omega,L}$ or $M_{\Omega,C}$ and the singular function $s$ corresponds to either $\log(\|\y-\x\|)$
 or $(\y-\x)\cdot\nu_{\bf{y}}/\|\y-\x\|^{2}$. In complex notation, the latter is reduced to a Cauchy-type singularity.

Both $M_{\Gamma,L}$ and $M_{\Omega,L}$ contain the factor $I_{1}(\alpha \|\x-\y\|)$, which grows like $e^{\alpha \|\x-\y\|}/\sqrt{\alpha \|\x-\y\|}$. The scaling with $\alpha$ can make $I_{1}$ grow too fast over a single panel to be accurately approximated by e.g. a $15$th degree polynomial or even a $31$th degree polynomial. The product integration relies on $\phi$ being well approximated by such a polynomial \eqref{eq:polyapprox}, otherwise the result may be very inaccurate. An adaptive time stepper will adjust the time step to satisfy the given tolerance, potentially decreasing it until the algorithm stalls.

This problem is not unique to the modified Helmholtz equation, but appears for biharmonic and Stokes equations as well. One solution is an algorithm presented in a separate paper, see \cite{2019arXiv190607713K}. By local refinement of panels through adaptive recursive bisection a kernel-split quadrature with product integration can be used successfully for a wide range of $\alpha$. It is ensured that the new panels are of adequate size to accurately approximate $\phi$ with polynomial interpolation. The method is effective in terms of computations, as the increased cost scales as $\log(\alpha)$. Moreover,
$K_{1}(\alpha \|\x-\y\|)\sim \sqrt{\pi/(\alpha \|\x-\y\|)}\,e^{-\alpha \|\x-\y\|}$ for large arguments, i.e. $K_{1}$ is very localised for large $\alpha$ and only a small portion of the boundary $\Gamma$ needs to be upsampled.

\section{Numerical results}
\label{s:results}

In this section we present numerical results, starting with a study of the modified Helmholtz equation to confirm that the parameters for PUX can be set as in \cite{FRYKLUNDPUX} for the Poisson equation. It forms the basis for the second numerical experiment, where the modified Helmholtz equation is solved on a more complex domain. The heat equation is solved on the same domain, for a range of set tolerances with an adaptive time stepper for different grid resolutions. Finally, the Allen-Cahn equation , a reaction-diffusion problem, is solved with randomised initial data.

To compute the errors we consider an evaluation grid. It consists of $N_{\text{eval}}^2$ uniformly distributed nodes over the computational domain $B$. We evaluate the numerical solution and an analytical or computed reference solution on the nodes that fall inside $\Omega$. The cardinality of this set of nodes as $N_{\text{eval},\Omega}$. Two different errors are computed: the relative $\ell_{2}$-error and the relative discrete $\ell_{\infty}$-error, defined as $\|\mathbf{u}_{\text{solution}}-\mathbf{u}_{\text{numerical}}\|_{\ell_{p}}/\|\mathbf{u}_{\text{solution}}\|_{\ell_{p}}$
where
\begin{equation}
\|\mathbf{u}\|_{\ell_{2}} = \frac{1}{N_{\text{eval},\Omega}}\sqrt{\sum\limits_{i=1}^{N_{\text{eval},\Omega}}|u_{i}|^{2}}
\label{eq:ell2}
\end{equation}
and
\begin{equation}
\|\mathbf{u}\|_{\ell_{\infty}} = \max |u_{i}|,\quad i = 1,\ldots,N_{\text{eval},\Omega},
\label{eq:ellinfty}
\end{equation}
for a vector $\bm{u}$ of length $N_{\text{eval},\Omega}$. When referring to the errors we mean both of them.

The following parameters are user specified in the numerical experiments: the length $L$ for the computation domain $B = [-L,L]^{2}$, the resolution $N_{u}$, partition radius $R$ and the number $N_{P,n}$ of Gauss-Legendre panels for each component curve $\Gamma_{n}$.
We use set the shape parameter $\varepsilon = 2$ for all numerical experiments and set the number of Gauss-Legendre nodes $N_Q =
16$.
\subsection{Example $1$: Study of weight functions}
\label{ss:eaxmple1}
We now solve the modified Helmholtz equation \eqref{eq:ModHelmEq}--\eqref{eq:ModHelmEqBC} for
\begin{equation}
  \label{eq:results1sol}
  u(x,y) = \sin(2\pi x)\sin(2\pi y)\exp(-(x^2+y^2)),
\end{equation}
to confirm that the parameters $N_{\phi}$ and $c$ and the function $\psi^{k}$ can be set by \eqref{eq:Pdef},  \eqref{eq:setM} and \eqref{eq:coarsen}, as in \cite{FRYKLUNDPUX} for the Poisson equation. To reduce the complexity of the problem, assume the corresponding right hand side to be known in all of $\mathbb{R}^{2}$, not just $\Omega$. To isolate the influence of the choice of weight function $\psi$, see Table \ref{tab:Wu}, the actual values of $f$ are used as values for the local extensions $f_{i,E}$, instead of the extrapolated ones $A_{i,E}f_{\mathbf{X}_{i}}$. Compact support is still enforced via PUX, but blending with the zero partitions reduces the regularity of $f^{e}$ to $k$.

The computational domain is the unit circle centred at $(17/701,5/439)$, contained in the box $B  =[-L,L]^{2}$, with $L  = 1.5$. The resolution $N_u$ attains values between $40$ and $500$ and for the evaluation grid use $N_{\text{e}} = 1000$. The partition radius and the number of panels are set such that only the resolution of the uniform grid $\mathbf{X}$ limits the accuracy. In this case the partition radius is $R = 0.4$ and the number of panels $N_P = 32$. This means that the rate of convergence is only dependent on the regularity of the extension  and we can study the influence of choice of Wu-function. Furthermore, we set $\alpha^{2} = 10$.

In Figure \ref{fig:wul2max} the errors for solving the modified Helmholtz equation are plotted as functions of the number of grid points for different Wu-functions. The behaviour of the errors is as for the Poisson equation in \cite{FRYKLUNDPUX}: $\psi^{k}$ with few continuous derivatives requires less points to be represented then $\psi^{k}$ with a larger $k$. Consequently, high regularity can increase the error, since $\psi^{k}$ is not sufficiently resolved. Compare the errors for using $\psi^{1}$ and $\psi^{5}$ in Figure \ref{fig:wul2max} for $N_{u}\sim 40$.
 As the grid is refined the decay is spectral until the errors is limited by an algebraic tail. The algebraic tail has a slope of $4+k$, as expected. The $\ell_{\infty}$-error is about one to two digits less accurate than the $\ell_{2}$-error, which is consistent for all numerical experiments in this paper. The reason is that there is almost always some target points close to the boundary for which the special quadrature does not give optimal results, e.g. at the intersection of two panels.

\begin{figure}
  \centering
%
\definecolor{mycolor1}{rgb}{0.00000,0.44700,0.74100}%
\definecolor{mycolor2}{rgb}{0.85000,0.32500,0.09800}%
\definecolor{mycolor3}{rgb}{0.92900,0.69400,0.12500}%
\definecolor{mycolor4}{rgb}{0.49400,0.18400,0.55600}%
\definecolor{mycolor5}{rgb}{0.46600,0.67400,0.18800}%
\definecolor{mycolor6}{rgb}{0.00000,0.74902,0.74902}%
\begin{tikzpicture}
\begin{loglogaxis}[%
xticklabel style={/pgf/number format/fixed},
xmin=40,
xmax=500,
xminorticks=true,
xlabel style={font=\Large\color{white!15!black}},
xlabel={$N_{u}$},
xtick={100,200,300,400,500},
xticklabels={$100$,$200$,$$,$400$,$$},
xticklabel style = {font=\Large},
ymin=1e-15,
ymax=1e-3,
yminorticks=true,
ylabel style={font=\Large\color{white!15!black}},
ylabel={Relative error},
ytick={1e-4,1e-5,1e-6,1e-7,1e-8,1e-9,1e-10,1e-11,1e-12,1e-13,1e-14},
yticklabels={$10^{-4}$,$$,$10^{-6}$,$$,$10^{-8}$,$$,$10^{-10}$,$$,$10^{-12}$,$$,$10^{-14}$},
yticklabel style = {font=\Large},
axis background/.style={fill=white},
xmajorgrids,
xminorgrids,
ymajorgrids,
legend style={font=\Large,at={(0.712,0.831)}, anchor=south west, legend cell align=left, align=left, draw=white!15!black},
reverse legend,
clip mode=individual,
]
\addplot [color=black, dashed, line width=1.0pt]
  table[row sep=crcr]{%
52	5.27294025687297e-06\\
60	1.26058978366927e-06\\
68	3.60576279945235e-07\\
78	9.14408246899256e-08\\
90	2.18605578083952e-08\\
116	1.72786940311198e-09\\
294	1.70530256582424e-13\\
};
\addlegendentry{$\mathcal{O}(N_{u}^{-10})$}

\addplot [color=mycolor2, line width=1.0pt, mark=*, mark options={solid, fill=mycolor2}, mark size=2.0pt]
  table[row sep=crcr]{%
  40	8.11026721751921e-05\\
  46	2.09599299600757e-05\\
  52	7.49995866963641e-06\\
  60	3.53201094763984e-06\\
  68	1.9957877669138e-06\\
  78	5.46101459090552e-07\\
  90	1.18388830678251e-07\\
  102	3.20092462111397e-08\\
  116	1.88963209389891e-08\\
  132	2.53900764685601e-09\\
  152	6.77583349798731e-10\\
  174	2.46205147548264e-10\\
  198	3.39332353064431e-11\\
  226	6.62057859080878e-12\\
  258	1.63721285839264e-12\\
  294	1.57678996091115e-13\\
  336	3.02220593867358e-14\\
  384	1.95671948828294e-14\\
  438	1.9119940279242e-14\\
  500	1.90550413795787e-14\\
};
\addlegendentry{$\ell_{2}$}

\addplot [color=mycolor1, line width=1.0pt, mark=square*, mark options={solid, fill=mycolor1}, mark size=2.0pt]
  table[row sep=crcr]{%
  40	0.000275474542317884\\
  46	8.43351813860872e-05\\
  52	4.1397775226551e-05\\
  60	1.83350843527973e-05\\
  68	1.00323409941675e-05\\
  78	3.7528112527237e-06\\
  90	5.04093088863687e-07\\
  102	2.21164592700261e-07\\
  116	1.29609218147426e-07\\
  132	1.66395851085275e-08\\
  152	4.75761185832281e-09\\
  174	1.70274061796918e-09\\
  198	2.60540336968883e-10\\
  226	4.71711737407606e-11\\
  258	1.55311182815372e-11\\
  294	1.87415266464433e-12\\
  336	1.89313051995576e-13\\
  384	1.43846793471683e-13\\
  438	1.38872184501159e-13\\
  500	1.42512820116458e-13\\
};
\addlegendentry{$\ell_{\infty}$}

\end{loglogaxis}
\end{tikzpicture}%
\caption{Error in numerical solution for the modified Helmholtz equation with \eqref{eq:results1sol} and $\alpha^{2} = 10$, but with local extensions given by analytic expression. The errors are plotted as functions of $N_{u}$ in loglog-scale with $\psi^{k}$ chosen according to \eqref{eq:Ck} for each different value of $N_{u}$.}
\label{fig:wuauto}
\end{figure}
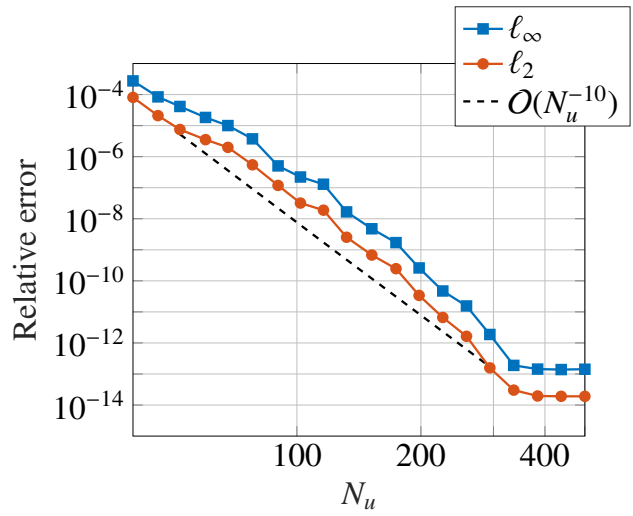

We now solve the modified Helmholtz equation in the same numerical setting, but let $\psi^{k}$ be set automatically by \eqref{eq:Ck}. The result is presented in Figure  \ref{fig:wuauto} and the lines follows the corresponding lowest errors in Figure \ref{fig:wul2max}. Thus  \eqref{eq:Ck} indeed chooses $\psi^{k}$ correctly for a given $N_{u}$ and we can set the PUX parameters for the modified Helmholtz equation as for the Poisson equation. This holds for $\alpha^2$ from $10$ to $10^{5}$ as well, as is shown in the following numerical experiment. Moreover, the error decreases as that of a tenth order method. For the subsequent numerical experiments $\psi^{k}$, $N_{\phi}$ and $c$
 are set by \eqref{eq:Ck}, \eqref{eq:setM} and \eqref{eq:coarsen}.
\FloatBarrier

\begin{figure}
  \centering
%
%
\definecolor{mycolor1}{rgb}{0.00000,0.44700,0.74100}%
\definecolor{mycolor2}{rgb}{0.85000,0.32500,0.09800}%
\definecolor{mycolor3}{rgb}{0.92900,0.69400,0.12500}%
\definecolor{mycolor4}{rgb}{0.49400,0.18400,0.55600}%
\definecolor{mycolor5}{rgb}{0.46600,0.67400,0.18800}%
\definecolor{mycolor6}{rgb}{0.00000,0.74902,0.74902}%
\begin{tikzpicture}

\begin{loglogaxis}[%
xticklabel style={/pgf/number format/fixed},
xmin=40,
xmax=500,
xminorticks=true,
xlabel style={font=\Large\color{white!15!black}},
xlabel={$N_{u}$},
xtick={100,200,300,400,500,600,700,800,900,1000},
xticklabels={$100$,$200$,$$,$400$,$$,$600$,$$,$800$,$$,$1000$},
xticklabel style = {font=\Large},
ymin=1e-15,
ymax=0.1,
yminorticks=true,
ylabel style={font=\Large\color{white!15!black}},
ylabel={Relative $\ell_{2}$--error},
ytick={1e-1,1e-2,1e-3,1e-4,1e-5,1e-6,1e-7,1e-8,1e-9,1e-10,1e-11,1e-12,1e-13,1e-14},
yticklabels={$$,$10^{-2}$,$$,$10^{-4}$,$$,$10^{-6}$,$$,$10^{-8}$,$$,$10^{-10}$,$$,$10^{-12}$,$$,$10^{-14}$},
yticklabel style = {font=\Large},
axis background/.style={fill=white},
xmajorgrids,
xminorgrids,
ymajorgrids,
legend style={font=\normalsize,at={(0.674,0.594)}, anchor=south west, legend cell align=left, align=left, draw=white!15!black},
clip mode=individual,
]
\addplot [color=mycolor1, line width=1.0pt, mark=x, mark options={solid, mycolor1}, mark size=3.0pt]
  table[row sep=crcr]{%
  40	8.11026721751921e-05\\
  46	2.09599299600757e-05\\
  52	7.49995866963641e-06\\
  60	3.53201094763984e-06\\
  68	2.10966439722754e-06\\
  78	1.1798107100431e-06\\
  90	4.64040874463504e-07\\
  102	2.85993817580419e-07\\
  116	1.52722409666461e-07\\
  132	8.73139876193984e-08\\
  152	3.77581184529023e-08\\
  174	1.91857558075753e-08\\
  198	1.11478406432104e-08\\
  226	5.51250626996815e-09\\
  258	3.10869787580898e-09\\
  294	1.56423224330004e-09\\
  336	8.58615018247101e-10\\
  384	4.22557440604256e-10\\
  438	2.56620116371661e-10\\
  500	1.48979148869209e-10\\
};
\addlegendentry{$\psi^{1}$}

\addplot [color=mycolor2, line width=1.0pt, mark=triangle, mark options={solid, rotate=180, mycolor2}, mark size=2.0pt]
  table[row sep=crcr]{%
  40	0.000238153178066147\\
  46	7.55745448979556e-05\\
  52	3.61219372045641e-05\\
  60	8.19252330215551e-06\\
  68	1.9957877669138e-06\\
  78	5.46101459090552e-07\\
  90	1.18388830678251e-07\\
  102	3.20092462111397e-08\\
  116	1.64142059242487e-08\\
  132	6.65346467288964e-09\\
  152	3.82100938181814e-09\\
  174	1.65607644561615e-09\\
  198	7.15518345205048e-10\\
  226	3.89988497445154e-10\\
  258	1.59686059667012e-10\\
  294	7.73391905391978e-11\\
  336	3.28943322607672e-11\\
  384	1.54029571854314e-11\\
  438	6.86747722353057e-12\\
  500	3.39817964961677e-12\\
};
\addlegendentry{$\psi^{2}$}

\addplot [color=mycolor3, line width=1.0pt, mark=square*, mark options={solid, mycolor3}, mark size=2.0pt]
  table[row sep=crcr]{%
  40	0.000357830894682471\\
  46	0.000156009353239782\\
  52	7.63437079709362e-05\\
  60	2.17735220199271e-05\\
  68	6.24601919483422e-06\\
  78	2.03224779226462e-06\\
  90	4.65819000468834e-07\\
  102	1.01444242560521e-07\\
  116	1.88963209389891e-08\\
  132	2.53900764685601e-09\\
  152	6.77583349798731e-10\\
  174	2.46205147548264e-10\\
  198	1.08823519877219e-10\\
  226	3.73776503631255e-11\\
  258	1.78023125782223e-11\\
  294	6.5348926280681e-12\\
  336	2.74955023331007e-12\\
  384	9.93466151536616e-13\\
  438	4.4945327549947e-13\\
  500	1.9912079661066e-13\\
};
\addlegendentry{$\psi^{3}$}

\addplot [color=mycolor4, line width=1.0pt, mark=+, mark options={solid, mycolor4}, mark size=3.0pt]
  table[row sep=crcr]{%
  40	0.000466693163709611\\
  46	0.000264205762033149\\
  52	0.000124890294871057\\
  60	4.50345739457359e-05\\
  68	1.60453436787075e-05\\
  78	6.28533350388946e-06\\
  90	1.79949492271873e-06\\
  102	5.29281203165131e-07\\
  116	1.24714528945196e-07\\
  132	1.47878379744415e-08\\
  152	3.60745678322932e-09\\
  174	3.75465821690614e-10\\
  198	3.39332353064431e-11\\
  226	6.62057859080878e-12\\
  258	1.63721285839264e-12\\
  294	5.84432766266613e-13\\
  336	1.93707069905459e-13\\
  384	7.01276336866323e-14\\
  438	2.90942111263859e-14\\
  500	2.0484639559957e-14\\
};
\addlegendentry{$\psi^{4}$}

\addplot [color=mycolor5, line width=1.0pt, mark=*, mark options={solid, mycolor5}, mark size=2.0pt]
  table[row sep=crcr]{%
  40	0.000540925044420423\\
  46	0.000350692890327795\\
  52	0.000163359274152613\\
  60	6.76928460904015e-05\\
  68	2.67068515381524e-05\\
  78	1.15565201967908e-05\\
  90	3.83675567370214e-06\\
  102	1.33361978564413e-06\\
  116	3.62355260160731e-07\\
  132	5.50644076367508e-08\\
  152	1.51639164682651e-08\\
  174	2.37942951208565e-09\\
  198	2.53419935632716e-10\\
  226	4.07249965086432e-11\\
  258	2.79109457760682e-12\\
  294	1.57678996091115e-13\\
  336	3.02220593867358e-14\\
  384	1.95671948828294e-14\\
  438	1.9119940279242e-14\\
  500	1.90550413795787e-14\\
};
\addlegendentry{$\psi^{5}$}

\addplot [color=mycolor1, dashed, line width=1.0pt]
  table[row sep=crcr]{%
  40	4.89440657839133e-05\\
  46	2.43338508386708e-05\\
  52	1.31820599321596e-05\\
  60	6.4453090744248e-06\\
  68	3.44711233482761e-06\\
  78	1.73590912686875e-06\\
  90	8.48764980994212e-07\\
  102	4.53940718989645e-07\\
  116	2.38621765089383e-07\\
  132	1.25063432438032e-07\\
  152	6.17706342011079e-08\\
  174	3.14234423162973e-08\\
  198	1.64692585926627e-08\\
  226	8.50076096397747e-09\\
  258	4.38431218159681e-09\\
  294	2.28172524751887e-09\\
  336	1.17031726791533e-09\\
  384	6.00266184138581e-10\\
  438	3.10906246373171e-10\\
  500	1.60379914760727e-10\\
};
\addlegendentry{$\mathcal{O}(N_{u}^{-5})$}

\addplot [color=mycolor2, dotted, line width=1.0pt]
  table[row sep=crcr]{%
  40	2.44140625e-05\\
  46	1.05548729470595e-05\\
  52	5.05801296467293e-06\\
  60	2.14334705075446e-06\\
  68	1.01145490252146e-06\\
  78	4.44050520904071e-07\\
  90	1.88167642315892e-07\\
  102	8.87971382186192e-08\\
  116	4.10442254667416e-08\\
  132	1.89041177063835e-08\\
  152	8.10846174955253e-09\\
  174	3.60333392300612e-09\\
  198	1.65962075885945e-09\\
  226	7.50497699114281e-10\\
  258	3.39064064844926e-10\\
  294	1.5485171426828e-10\\
  336	6.94967244413335e-11\\
  384	3.1189804587549e-11\\
  438	1.41630026684205e-11\\
  500	6.4e-12\\
};
\addlegendentry{$\mathcal{O}(N_{u}^{-6})$}

\addplot [color=mycolor3, dashdotted, line width=1.0pt]
  table[row sep=crcr]{%
  40	6.103515625e-06\\
  46	2.29453759718684e-06\\
  52	9.72694800898641e-07\\
  60	3.57224508459076e-07\\
  68	1.48743368017862e-07\\
  78	5.69295539620604e-08\\
  90	2.09075158128769e-08\\
  102	8.70560178613914e-09\\
  116	3.53829529885703e-09\\
  132	1.43213012927148e-09\\
  152	5.33451430891614e-10\\
  174	2.07088156494605e-10\\
  198	8.3819230245427e-11\\
  226	3.32078627926673e-11\\
  258	1.31420180172452e-11\\
  294	5.267065111166e-12\\
  336	2.06835489408731e-12\\
  384	8.12234494467421e-13\\
  438	3.23356225306404e-13\\
  500	1.28e-13\\
};
\addlegendentry{$\mathcal{O}(N_{u}^{-7})$}

\addplot [color=mycolor4, dashdotdotted, line width=1.0pt]
  table[row sep=crcr]{%
  40	6.07463333974452e-06\\
  46	1.98580841422755e-06\\
  52	7.44686105765106e-07\\
  60	2.37022730525011e-07\\
  68	8.70820608532183e-08\\
  78	2.90564918573186e-08\\
  90	9.24825773729658e-09\\
  102	3.39780636769149e-09\\
  116	1.21432821552648e-09\\
  132	4.31925207294455e-10\\
  152	1.39717657749983e-10\\
  174	4.73811954236821e-11\\
  198	1.68530487833227e-11\\
  226	5.84968508695466e-12\\
  258	2.0278804682979e-12\\
  294	7.13216458682767e-13\\
  336	2.4506753410285e-13\\
  384	8.42073896922852e-14\\
  438	2.93905095747926e-14\\
  500	1.01915435661695e-14\\
};
\addlegendentry{$\mathcal{O}(N_{u}^{-8})$}

\end{loglogaxis}
\end{tikzpicture}%
%
%
\definecolor{mycolor1}{rgb}{0.00000,0.44700,0.74100}%
\definecolor{mycolor2}{rgb}{0.85000,0.32500,0.09800}%
\definecolor{mycolor3}{rgb}{0.92900,0.69400,0.12500}%
\definecolor{mycolor4}{rgb}{0.49400,0.18400,0.55600}%
\definecolor{mycolor5}{rgb}{0.46600,0.67400,0.18800}%
\definecolor{mycolor6}{rgb}{0.00000,0.74902,0.74902}%
\begin{tikzpicture}

\begin{loglogaxis}[%
xticklabel style={/pgf/number format/fixed},
xmin=40,
xmax=500,
xminorticks=true,
xlabel style={font=\Large\color{white!15!black}},
xlabel={$N_{u}$},
xtick={100,200,300,400,500,600,700,800,900,1000},
xticklabels={$100$,$200$,$$,$400$,$$,$600$,$$,$800$,$$,$1000$},
xticklabel style = {font=\Large},
ymin=1e-15,
ymax=0.1,
yminorticks=true,
ylabel style={font=\Large\color{white!15!black}},
ylabel={Relative $\ell_{\infty}$--error},
ytick={1e-1,1e-2,1e-3,1e-4,1e-5,1e-6,1e-7,1e-8,1e-9,1e-10,1e-11,1e-12,1e-13,1e-14},
yticklabels={$$,$10^{-2}$,$$,$10^{-4}$,$$,$10^{-6}$,$$,$10^{-8}$,$$,$10^{-10}$,$$,$10^{-12}$,$$,$10^{-14}$},
yticklabel style = {font=\Large},
axis background/.style={fill=white},
xmajorgrids,
xminorgrids,
ymajorgrids,
legend style={font=\normalsize,at={(0.674,0.594)}, anchor=south west, legend cell align=left, align=left, draw=white!15!black},
clip mode=individual,
]
\addplot [color=mycolor1, line width=1.0pt, mark=x, mark options={solid, mycolor1}, mark size=3.0pt]
  table[row sep=crcr]{%
  40	0.000275474542317884\\
  46	8.43351813860872e-05\\
  52	4.1397775226551e-05\\
  60	1.83350843527973e-05\\
  68	1.25007092664863e-05\\
  78	6.50277200159115e-06\\
  90	2.42728175491164e-06\\
  102	2.04378243854308e-06\\
  116	1.15894536833788e-06\\
  132	7.38997370231323e-07\\
  152	3.58286780577342e-07\\
  174	1.98388950819128e-07\\
  198	1.09188524137578e-07\\
  226	6.36271306359103e-08\\
  258	3.4109860770847e-08\\
  294	2.11538255339499e-08\\
  336	1.04547787243129e-08\\
  384	4.62730378276219e-09\\
  438	3.20563123957431e-09\\
  500	2.41550759366591e-09\\
};
\addlegendentry{$\psi^{1}$}

\addplot [color=mycolor2, line width=1.0pt, mark=triangle, mark options={solid, rotate=180, mycolor2}, mark size=2.0pt]
  table[row sep=crcr]{%
  40	0.000814079577639214\\
  46	0.000286259742720819\\
  52	0.000171212352003751\\
  60	3.50963091031691e-05\\
  68	1.00323409941675e-05\\
  78	3.7528112527237e-06\\
  90	5.04093088863687e-07\\
  102	2.21164592700261e-07\\
  116	1.25890474612587e-07\\
  132	4.19047417020257e-08\\
  152	3.33579878104196e-08\\
  174	1.69176898726647e-08\\
  198	7.48206601359197e-09\\
  226	3.6475828373691e-09\\
  258	1.7030707763746e-09\\
  294	9.00354644321772e-10\\
  336	4.75221782589651e-10\\
  384	1.76601508729198e-10\\
  438	8.89452302418417e-11\\
  500	4.72943578427821e-11\\
};
\addlegendentry{$\psi^{2}$}

\addplot [color=mycolor3, line width=1.0pt, mark=square*, mark options={solid, mycolor3}, mark size=2.0pt]
  table[row sep=crcr]{%
  40	0.00126779594753728\\
  46	0.000539991422894946\\
  52	0.000398370098854147\\
  60	8.16171351855435e-05\\
  68	4.06656836868635e-05\\
  78	1.2996142243682e-05\\
  90	2.687287969227e-06\\
  102	6.73572156455648e-07\\
  116	1.29609218147426e-07\\
  132	1.66395851085275e-08\\
  152	4.75761185832281e-09\\
  174	1.70274061796918e-09\\
  198	1.06748094030881e-09\\
  226	4.06192552212415e-10\\
  258	1.57264564465206e-10\\
  294	6.57054134338022e-11\\
  336	2.77679893623482e-11\\
  384	9.47521274215791e-12\\
  438	4.09474237830713e-12\\
  500	2.57768118007834e-12\\
};
\addlegendentry{$\psi^{3}$}

\addplot [color=mycolor4, line width=1.0pt, mark=+, mark options={solid, mycolor4}, mark size=3.0pt]
  table[row sep=crcr]{%
  40	0.00171425507631862\\
  46	0.000947675259694775\\
  52	0.000663491223219709\\
  60	0.000166691311243456\\
  68	0.000101535384215434\\
  78	3.63041021930765e-05\\
  90	1.01263111530364e-05\\
  102	3.03898024974368e-06\\
  116	7.69726764517963e-07\\
  132	1.04193418583078e-07\\
  152	2.92027106209896e-08\\
  174	3.02005191647001e-09\\
  198	2.60540336968883e-10\\
  226	4.71711737407606e-11\\
  258	1.55311182815372e-11\\
  294	4.77951536732225e-12\\
  336	2.28832012644095e-12\\
  384	6.47866392853915e-13\\
  438	2.54788910847833e-13\\
  500	1.43652255690712e-13\\
};
\addlegendentry{$\psi^{4}$}

\addplot [color=mycolor5, line width=1.0pt, mark=*, mark options={solid, mycolor5}, mark size=2.0pt]
  table[row sep=crcr]{%
  40	0.00203339852528292\\
  46	0.00128529617233745\\
  52	0.000878525192310621\\
  60	0.000248030200413092\\
  68	0.000164275741583011\\
  78	6.52207220846238e-05\\
  90	2.07660379605379e-05\\
  102	7.27758367053953e-06\\
  116	2.27513443287829e-06\\
  132	3.96446212776324e-07\\
  152	1.22459880671672e-07\\
  174	1.82906807465636e-08\\
  198	2.32877245157145e-09\\
  226	4.42438859353895e-10\\
  258	3.21968086928043e-11\\
  294	1.87415266464433e-12\\
  336	1.89313051995576e-13\\
  384	1.43846793471683e-13\\
  438	1.38872184501159e-13\\
  500	1.42512820116458e-13\\
};
\addlegendentry{$\psi^{5}$}

\addplot [color=mycolor1, dashed, line width=1.0pt]
  table[row sep=crcr]{%
  40	0.000308816177750818\\
  46	0.000153536219061448\\
  52	8.31731752957417e-05\\
  60	4.06671509795316e-05\\
  68	2.17498084490775e-05\\
  78	1.09528461294804e-05\\
  90	5.3553449849589e-06\\
  102	2.86417230605136e-06\\
  116	1.50560155235974e-06\\
  132	7.89096912226787e-07\\
  152	3.89746353223884e-07\\
  174	1.98268517183175e-07\\
  198	1.0391399667184e-07\\
  226	5.36361756389211e-08\\
  258	2.76631397147248e-08\\
  294	1.43967130300792e-08\\
  336	7.38420275563174e-09\\
  384	3.78742357525338e-09\\
  438	1.96168579593921e-09\\
  500	1.01192885125388e-09\\
};
\addlegendentry{$\mathcal{O}(N_{u}^{-5})$}

\addplot [color=mycolor2, dotted, line width=1.0pt]
  table[row sep=crcr]{%
  40	0.000154042320429735\\
  46	6.65967460600248e-05\\
  52	3.19139042853642e-05\\
  60	1.35236056344349e-05\\
  68	6.38184899356413e-06\\
  78	2.80176937484679e-06\\
  90	1.18725755912734e-06\\
  102	5.60272065278606e-07\\
  116	2.58971555067416e-07\\
  132	1.19276919077607e-07\\
  152	5.11609349331688e-08\\
  174	2.27355000333534e-08\\
  198	1.04714990685416e-08\\
  226	4.73532035271642e-09\\
  258	2.13934961963215e-09\\
  294	9.77048264229198e-10\\
  336	4.38494687035755e-10\\
  384	1.9679436277416e-10\\
  438	8.93625055353248e-11\\
  500	4.03812700467324e-11\\
};
\addlegendentry{$\mathcal{O}(N_{u}^{-6})$}

\addplot [color=mycolor3, dashdotted, line width=1.0pt]
  table[row sep=crcr]{%
  40	6.103515625e-05\\
  46	2.29453759718684e-05\\
  52	9.72694800898641e-06\\
  60	3.57224508459076e-06\\
  68	1.48743368017862e-06\\
  78	5.69295539620604e-07\\
  90	2.09075158128769e-07\\
  102	8.70560178613914e-08\\
  116	3.53829529885703e-08\\
  132	1.43213012927148e-08\\
  152	5.33451430891614e-09\\
  174	2.07088156494605e-09\\
  198	8.3819230245427e-10\\
  226	3.32078627926673e-10\\
  258	1.31420180172452e-10\\
  294	5.267065111166e-11\\
  336	2.06835489408731e-11\\
  384	8.12234494467421e-12\\
  438	3.23356225306404e-12\\
  500	1.28e-12\\
};
\addlegendentry{$\mathcal{O}(N_{u}^{-7})$}

\addplot [color=mycolor4, dashdotdotted, line width=1.0pt]
  table[row sep=crcr]{%
  40	6.07463333974452e-05\\
  46	1.98580841422755e-05\\
  52	7.44686105765106e-06\\
  60	2.37022730525011e-06\\
  68	8.70820608532183e-07\\
  78	2.90564918573186e-07\\
  90	9.24825773729658e-08\\
  102	3.39780636769149e-08\\
  116	1.21432821552648e-08\\
  132	4.31925207294455e-09\\
  152	1.39717657749983e-09\\
  174	4.73811954236821e-10\\
  198	1.68530487833227e-10\\
  226	5.84968508695466e-11\\
  258	2.0278804682979e-11\\
  294	7.13216458682767e-12\\
  336	2.4506753410285e-12\\
  384	8.42073896922852e-13\\
  438	2.93905095747926e-13\\
  500	1.01915435661695e-13\\
};
\addlegendentry{$\mathcal{O}(N_{u}^{-8})$}

\end{loglogaxis}
\end{tikzpicture}%
%
\caption{Error in numerical solution for the modified Helmholtz equation with \eqref{eq:results1sol} and $\alpha^{2} = 10$, but with local extensions given by analytic expression. Errors are plotted as a function of $N_{u}$ in loglog-scale for $\psi^{k}$, where $k = 1,\,2,\,3,\,4,\,5$. See Table \ref{tab:Wu}. Left: relative $\ell_{2}$-error. Right: relative $\ell_{\infty}$-error.}
\label{fig:wul2max}
\end{figure}
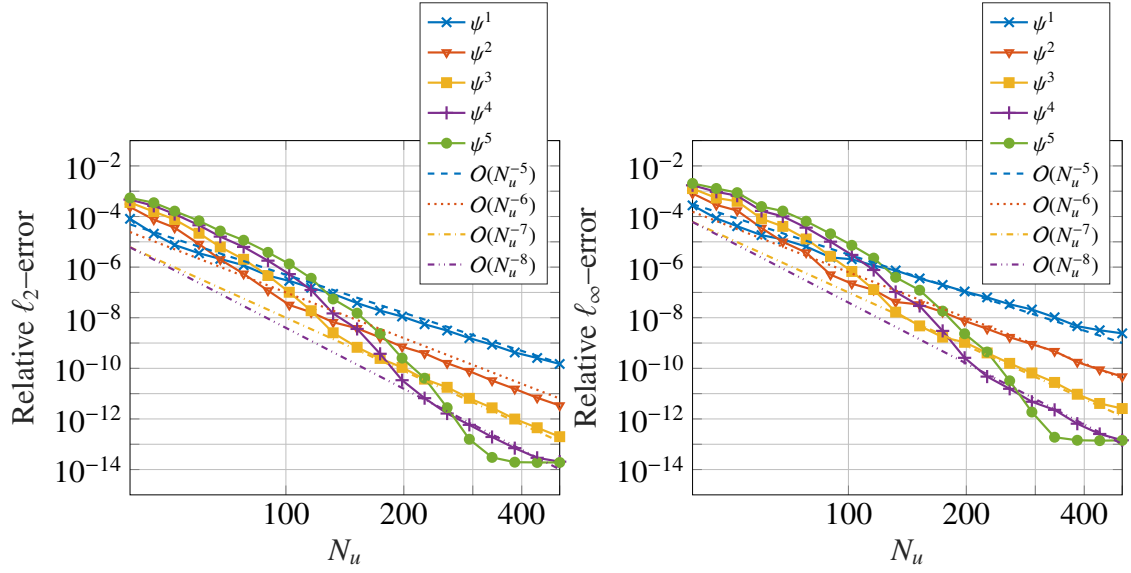

%

\FloatBarrier
\subsection{Example $2$: the modified Helmholtz equation on a multiply connected domain}
\label{ss:example2}
\begin{figure}
  \centering

  \includegraphics[width=0.49\textwidth]{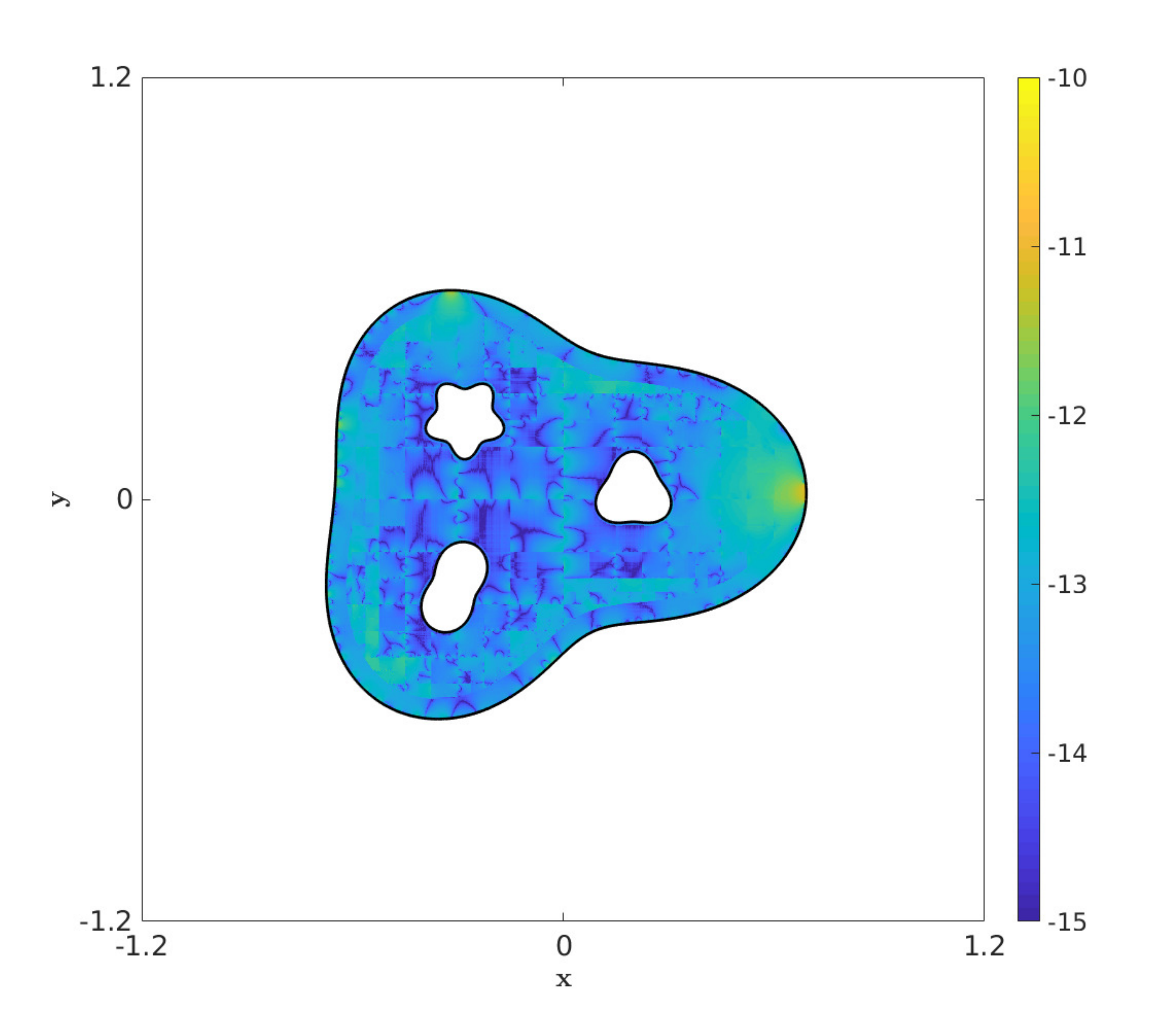}
  \includegraphics[width=0.49\textwidth]{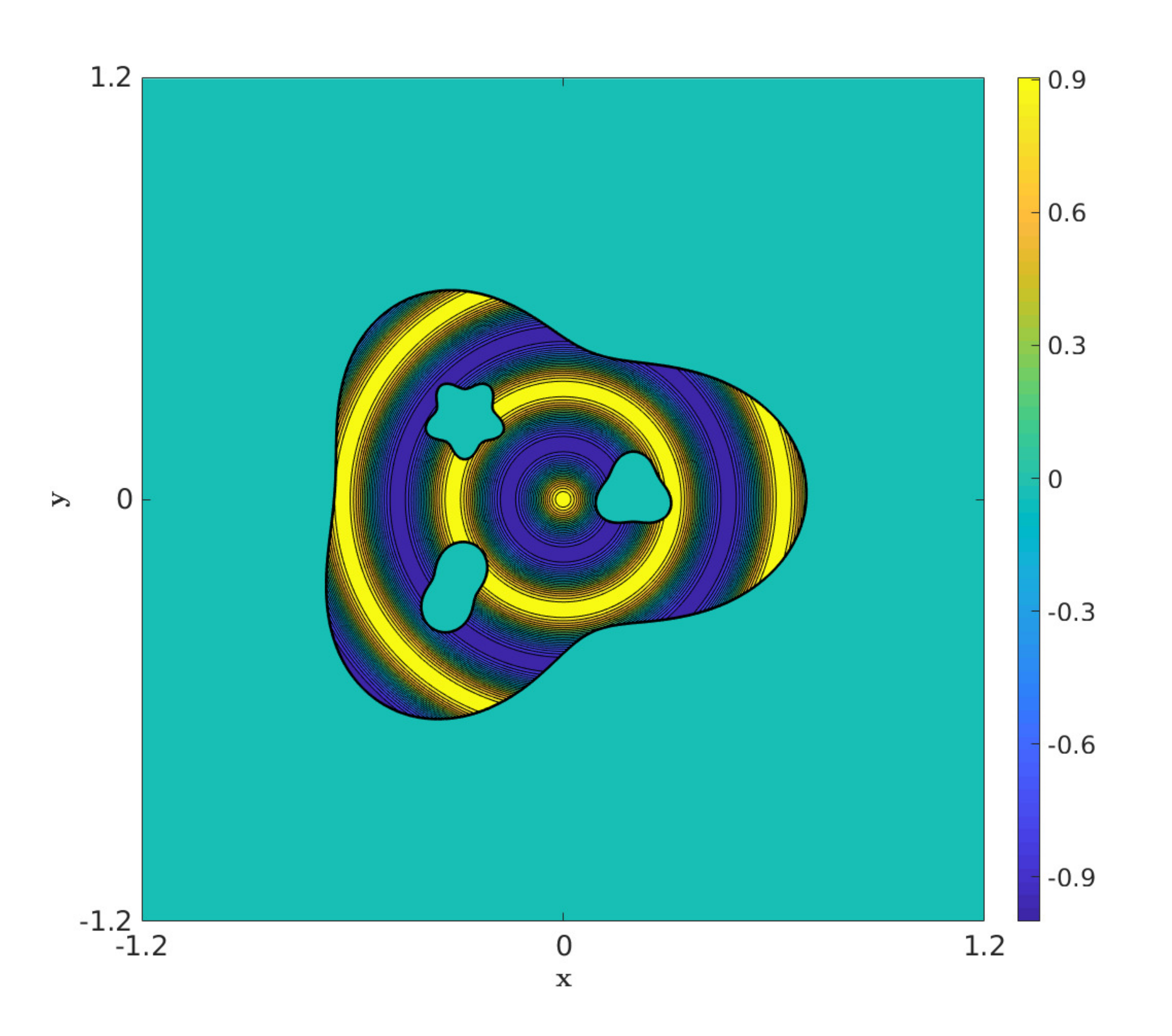}
\caption{For $N_u = 1000$ and $\alpha^2 = 10$: the left image shows the pointwise relative $\ell_2$-error for solving the modified Helmholtz equation for \eqref{eq:ex2}. The right image shows the solution \eqref{eq:ex2}.}
\label{fig:e2_errorplot}
\end{figure}

\begin{figure}
  \centering
  \includegraphics[width=0.49\textwidth]{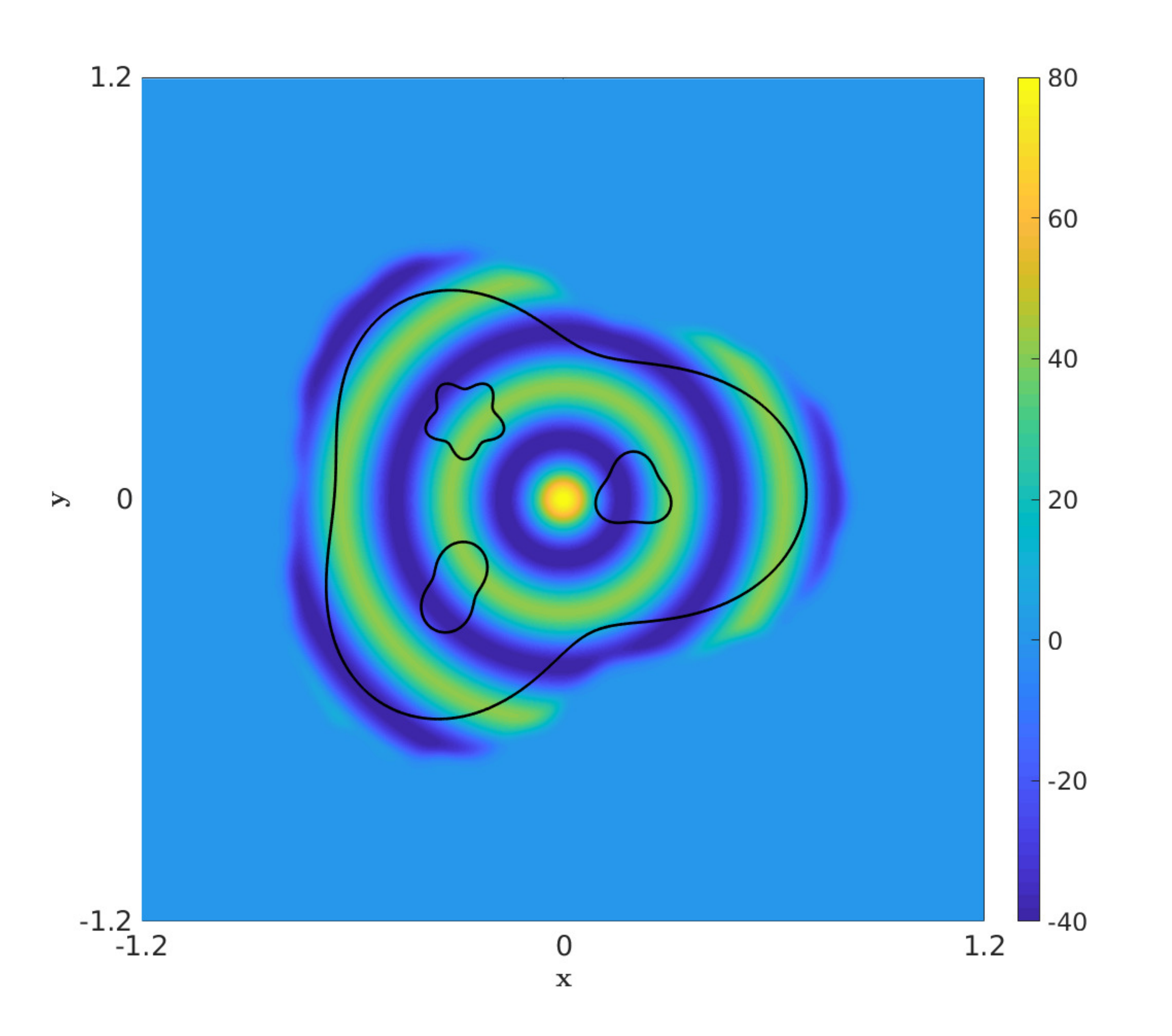}
    \includegraphics[width=0.49\textwidth]{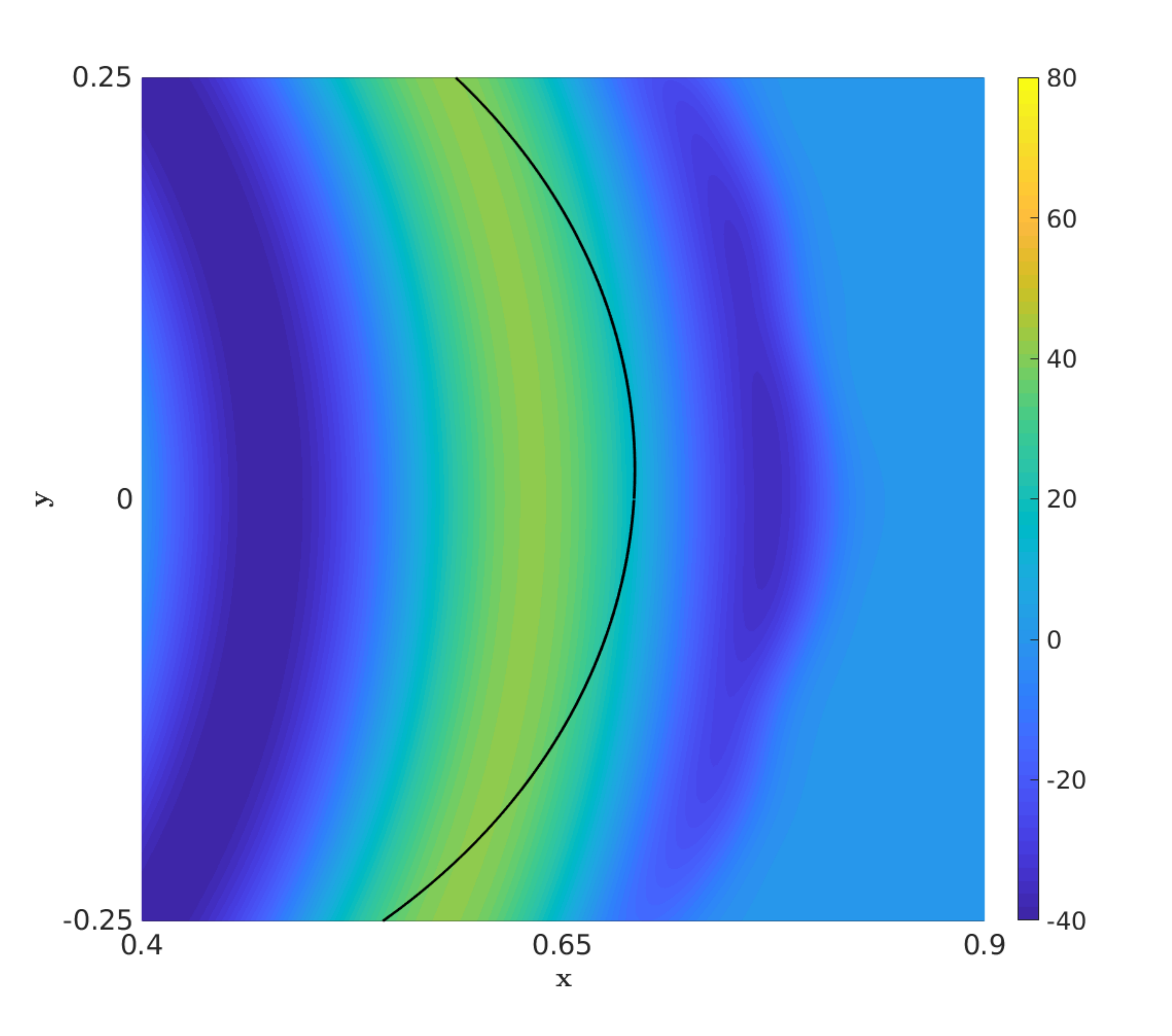}
\caption{Left: The right hand side of the modified Helmholtz equation for \eqref{eq:ex2}, extended with PUX. Right: magnification of left image. For both $N_u = 1000$ and $\alpha^2 = 10$}
\label{fig:ex2_solfe}
\end{figure}

 We now study the modified Helmholtz equation with a more complex setup for $\alpha^2 = 10^{n}$, $n = 1,2,3,4,5$. We take the solution to be
\begin{equation}
  \label{eq:ex2}
  u(x,y) = \cos\left(20\sqrt{x^2 + y^2}\right),
\end{equation}
on
 the multiply connected domain shown in Figure \ref{fig:e2_errorplot} and evaluate the right hand side in \eqref{eq:ModHelmEq} accordingly. The corresponding extension by PUX is shown in \ref{fig:ex2_solfe}, where $N_u = 1000$ and $k = 5$. The outer boundary is discretised into $80$ panels, and the boundaries of the cavities are discretised with $20$ panels each. Again all parameters are set such that only $N_{u}$ sets the bound for the error. We set $R = 0.23$ and $L=1.2$. The parameters $\psi^{k}$, $N_{\phi}$
  and $c$ are set by \eqref{eq:Ck}, \eqref{eq:setM} and \eqref{eq:coarsen}. The parameter $\alpha^2$ ranges from $10$ to $10^{5}$. The evaluation grid has a resolution of $N_{\text{eval}} = 1000$.

\begin{figure}
  \centering
%
%
\definecolor{mycolor1}{rgb}{0.00000,0.44700,0.74100}%
\definecolor{mycolor2}{rgb}{0.85000,0.32500,0.09800}%
\definecolor{mycolor3}{rgb}{0.92900,0.69400,0.12500}%
\definecolor{mycolor4}{rgb}{0.49400,0.18400,0.55600}%
\definecolor{mycolor5}{rgb}{0.46600,0.67400,0.18800}%
\definecolor{mycolor6}{rgb}{0.00000,0.74902,0.74902}%
\begin{tikzpicture}

\begin{loglogaxis}[%
xticklabel style={/pgf/number format/fixed},
xmin=100,
xmax=1000,
xminorticks=true,
xlabel style={font=\Large\color{white!15!black}},
xlabel={$N_{u}$},
xtick={100,200,300,400,500,600,700,800,900,1000},
xticklabels={$100$,$200$,$$,$400$,$$,$600$,$$,$800$,$$,$1000$},
xticklabel style = {font=\normalsize},
ymin=1e-15,
ymax=0.1,
yminorticks=true,
ylabel style={font=\Large\color{white!15!black}},
ylabel={Relative $\ell_{2}$--error},
ytick={1e-1,1e-2,1e-3,1e-4,1e-5,1e-6,1e-7,1e-8,1e-9,1e-10,1e-11,1e-12,1e-13,1e-14},
yticklabels={$$,$10^{-2}$,$$,$10^{-4}$,$$,$10^{-6}$,$$,$10^{-8}$,$$,$10^{-10}$,$$,$10^{-12}$,$$,$10^{-14}$},
yticklabel style = {font=\normalsize},
axis background/.style={fill=white},
xmajorgrids,
xminorgrids,
ymajorgrids,
legend style={font=\normalsize,at={(0.674,0.594)}, anchor=south west, legend cell align=left, align=left, draw=white!15!black},
clip mode=individual,
]
\addplot [color=mycolor1, line width=1.0pt, mark=x, mark options={solid, mycolor1}, mark size=3.0pt]
  table[row sep=crcr]{%
  100	0.000225782214261774\\
  114	0.000102929167683722\\
  128	4.38901946809698e-05\\
  144	2.41104332551026e-05\\
  162	6.98929352804723e-06\\
  184	3.27857701627264e-06\\
  208	7.04159530704679e-07\\
  234	3.49579399077685e-07\\
  264	7.97147302214591e-08\\
  298	1.09606900339661e-08\\
  336	3.70225707397689e-09\\
  380	8.57603531810569e-10\\
  428	9.56523766100256e-11\\
  484	1.01270688651425e-11\\
  546	6.65748158527511e-13\\
  616	1.46359817231979e-13\\
  696	1.3821047680121e-13\\
  786	1.37830919079752e-13\\
  886	1.37178111257187e-13\\
  1000	1.44051044450096e-13\\
};
\addlegendentry{$\alpha^{2} = 10$}

\addplot [color=mycolor2, line width=1.0pt, mark=triangle, mark options={solid, rotate=180, mycolor2}, mark size=2.0pt]
  table[row sep=crcr]{%
  100	0.000245818141577165\\
  114	0.000107029610267463\\
  128	4.96590095909307e-05\\
  144	2.65971967540886e-05\\
  162	7.33676738519671e-06\\
  184	3.63548670139635e-06\\
  208	7.87705924013685e-07\\
  234	3.70664743384383e-07\\
  264	8.57794566906167e-08\\
  298	1.23193437494359e-08\\
  336	4.00617884140342e-09\\
  380	9.56604001164798e-10\\
  428	1.09159951378141e-10\\
  484	1.11343245726196e-11\\
  546	7.48844531036495e-13\\
  616	6.35106782865687e-14\\
  696	2.73482057891514e-14\\
  786	3.26140761230481e-14\\
  886	2.51351514020699e-14\\
  1000	4.55328686038816e-14\\
};
\addlegendentry{$\alpha^{2} = 10^{2}$}

\addplot [color=mycolor3, line width=1.0pt, mark=square*, mark options={solid, mycolor3}, mark size=2.0pt]
  table[row sep=crcr]{%
  100	0.000520293834905425\\
  114	0.000225871509421597\\
  128	0.000111237296183673\\
  144	6.08013756826845e-05\\
  162	1.61090428653186e-05\\
  184	8.54388936940567e-06\\
  208	1.83517795466705e-06\\
  234	8.28556821312198e-07\\
  264	1.92051360066946e-07\\
  298	2.95592379832439e-08\\
  336	9.35134192037398e-09\\
  380	2.25490576985357e-09\\
  428	2.70406826721731e-10\\
  484	2.64356868665732e-11\\
  546	1.863463912935e-12\\
  616	1.5489667803323e-13\\
  696	4.88601333441789e-14\\
  786	6.46022592943389e-14\\
  886	3.9173559503965e-14\\
  1000	9.25453959474552e-14\\
};
\addlegendentry{$\alpha^{2} = 10^{3}$}

\addplot [color=mycolor4, line width=1.0pt, mark=+, mark options={solid, mycolor4}, mark size=3.0pt]
  table[row sep=crcr]{%
  100	0.0020765875369424\\
  114	0.000979298272259503\\
  128	0.000509648340353767\\
  144	0.000306978504591998\\
  162	7.9952836956823e-05\\
  184	4.60588465695873e-05\\
  208	1.00481801078406e-05\\
  234	4.57852245535203e-06\\
  264	1.08889624579304e-06\\
  298	1.75442040081831e-07\\
  336	5.69442918065915e-08\\
  380	1.34668728472113e-08\\
  428	1.72772719556165e-09\\
  484	1.63771271891743e-10\\
  546	1.20746881435446e-11\\
  616	1.03583214538052e-12\\
  696	3.25473865237281e-13\\
  786	3.77381015492488e-13\\
  886	2.24971885617828e-13\\
  1000	5.50769135351345e-13\\
};
\addlegendentry{$\alpha^{2} = 10^{4}$}

\addplot [color=mycolor5, line width=1.0pt, mark=*, mark options={solid, mycolor5}, mark size=2.0pt]
  table[row sep=crcr]{%
  100	0.0045264031523765\\
  114	0.00245709477531243\\
  128	0.00139252006177891\\
  144	0.000940781969743638\\
  162	0.000262431318818641\\
  184	0.00016867463296332\\
  208	4.00216793906553e-05\\
  234	1.9983479650203e-05\\
  264	5.21376554482039e-06\\
  298	9.00290896683963e-07\\
  336	3.23998858338306e-07\\
  380	7.94677657800863e-08\\
  428	1.11355284190986e-08\\
  484	1.08859146986765e-09\\
  546	8.5434752696188e-11\\
  616	7.68909971591013e-12\\
  696	2.46994838574283e-12\\
  786	2.42375926351394e-12\\
  886	1.57243032165402e-12\\
  1000	3.72423648192546e-12\\
};
\addlegendentry{$\alpha^{2} = 10^{5}$}

\end{loglogaxis}
\end{tikzpicture}%
%
%
\definecolor{mycolor1}{rgb}{0.00000,0.44700,0.74100}%
\definecolor{mycolor2}{rgb}{0.85000,0.32500,0.09800}%
\definecolor{mycolor3}{rgb}{0.92900,0.69400,0.12500}%
\definecolor{mycolor4}{rgb}{0.49400,0.18400,0.55600}%
\definecolor{mycolor5}{rgb}{0.46600,0.67400,0.18800}%
\definecolor{mycolor6}{rgb}{0.00000,0.74902,0.74902}%
\begin{tikzpicture}

\begin{loglogaxis}[%
xticklabel style={/pgf/number format/fixed},
xmin=100,
xmax=1000,
xminorticks=true,
xlabel style={font=\Large\color{white!15!black}},
xlabel={$N_{u}$},
xtick={100,200,300,400,500,600,700,800,900,1000},
xticklabels={$100$,$200$,$$,$400$,$$,$600$,$$,$800$,$$,$1000$},
xticklabel style = {font=\normalsize},
ymin=1e-15,
ymax=0.1,
yminorticks=true,
ylabel style={font=\Large\color{white!15!black}},
ylabel={Relative $\ell_{\infty}$--error},
ytick={1e-1,1e-2,1e-3,1e-4,1e-5,1e-6,1e-7,1e-8,1e-9,1e-10,1e-11,1e-12,1e-13,1e-14},
yticklabels={$$,$10^{-2}$,$$,$10^{-4}$,$$,$10^{-6}$,$$,$10^{-8}$,$$,$10^{-10}$,$$,$10^{-12}$,$$,$10^{-14}$},
yticklabel style = {font=\normalsize},
axis background/.style={fill=white},
xmajorgrids,
xminorgrids,
ymajorgrids,
legend style={font=\normalsize,at={(0.674,0.594)}, anchor=south west, legend cell align=left, align=left, draw=white!15!black},
clip mode=individual,
]
\addplot [color=mycolor1, line width=1.0pt, mark=x, mark options={solid, mycolor1}, mark size=3.0pt]
  table[row sep=crcr]{%
  100	0.00209744341641153\\
  114	0.00066153791565088\\
  128	0.000469110863012094\\
  144	0.000297253958465151\\
  162	7.53421059095405e-05\\
  184	4.31007997673835e-05\\
  208	1.03442560079885e-05\\
  234	6.13566769180923e-06\\
  264	1.02736898699445e-06\\
  298	1.73493619914929e-07\\
  336	5.82307392235663e-08\\
  380	1.99697625102866e-08\\
  428	1.98519878313883e-09\\
  484	2.12041051383958e-10\\
  546	1.39558364864462e-11\\
  616	1.45117251548761e-12\\
  696	7.95460919356152e-13\\
  786	8.33055846527523e-13\\
  886	8.02122257503903e-13\\
  1000	1.55764290354916e-12\\
};
\addlegendentry{$\alpha^{2} = 10$}

\addplot [color=mycolor2, line width=1.0pt, mark=triangle, mark options={solid, rotate=180, mycolor2}, mark size=2.0pt]
  table[row sep=crcr]{%
  100	0.00242812583628654\\
  114	0.000781215578750473\\
  128	0.000556802793354333\\
  144	0.000354830218778151\\
  162	9.00953160571967e-05\\
  184	5.15663286538871e-05\\
  208	1.23752959821927e-05\\
  234	7.37903755398983e-06\\
  264	1.23622057279787e-06\\
  298	2.0923966220777e-07\\
  336	7.00578059245647e-08\\
  380	2.40442911136906e-08\\
  428	2.39922437472739e-09\\
  484	2.5622515220648e-10\\
  546	1.68837166469878e-11\\
  616	1.02190478301627e-12\\
  696	3.82804898890771e-13\\
  786	7.42406136566875e-13\\
  886	3.45445894112133e-13\\
  1000	1.56796797767818e-12\\
};
\addlegendentry{$\alpha^{2} = 10^{2}$}

\addplot [color=mycolor3, line width=1.0pt, mark=square*, mark options={solid, mycolor3}, mark size=2.0pt]
  table[row sep=crcr]{%
  100	0.00542243774927383\\
  114	0.00188622636303912\\
  128	0.00136517023444418\\
  144	0.000903924608977409\\
  162	0.000231037145343126\\
  184	0.000137506414259254\\
  208	3.19609799194223e-05\\
  234	1.94401543400233e-05\\
  264	3.27478615347335e-06\\
  298	5.61719100233862e-07\\
  336	1.86446224709211e-07\\
  380	6.37078891863069e-08\\
  428	6.52192461148034e-09\\
  484	6.93976653920489e-10\\
  546	4.5920378610733e-11\\
  616	2.82679435414947e-12\\
  696	1.09701137063217e-12\\
  786	2.0110579868061e-12\\
  886	9.32531829533942e-13\\
  1000	4.19231316328724e-12\\
};
\addlegendentry{$\alpha^{2} = 10^{3}$}

\addplot [color=mycolor4, line width=1.0pt, mark=+, mark options={solid, mycolor4}, mark size=3.0pt]
  table[row sep=crcr]{%
  100	0.0218030689334785\\
  114	0.0094563218129394\\
  128	0.00660620403426773\\
  144	0.00487802719339613\\
  162	0.00129880073832062\\
  184	0.000778024466043838\\
  208	0.000182918070789858\\
  234	0.000116035991007351\\
  264	2.10102915099162e-05\\
  298	3.71550643790702e-06\\
  336	1.21550868703585e-06\\
  380	3.93942560406053e-07\\
  428	4.61041529953319e-08\\
  484	4.77802053477937e-09\\
  546	3.22321391799122e-10\\
  616	1.99815719525995e-11\\
  696	7.52897744149583e-12\\
  786	1.3708312263106e-11\\
  886	6.47970566092255e-12\\
  1000	3.08034708851337e-11\\
};
\addlegendentry{$\alpha^{2} = 10^{4}$}

\addplot [color=mycolor5, line width=1.0pt, mark=*, mark options={solid, mycolor5}, mark size=2.0pt]
  table[row sep=crcr]{%
  100	0.050767701110088\\
  114	0.0311386185448961\\
  128	0.0191526515310423\\
  144	0.0173480733006062\\
  162	0.00484268488599526\\
  184	0.00284788757243072\\
  208	0.000691330103277185\\
  234	0.00053459820324014\\
  264	0.00010697288364892\\
  298	2.07406217843662e-05\\
  336	7.1883071218196e-06\\
  380	2.40187577538197e-06\\
  428	3.35806996487624e-07\\
  484	3.25880161744522e-08\\
  546	2.53547538520149e-09\\
  616	1.6739409858247e-10\\
  696	5.29997157272556e-11\\
  786	9.48975342751666e-11\\
  886	5.52378143225976e-11\\
  1000	2.50879095276002e-10\\
};
\addlegendentry{$\alpha^{2} = 10^{5}$}

\end{loglogaxis}
\end{tikzpicture}%
\caption{The errors for example $2$ for various resolutions of the uniform grid, over a range of values for $\alpha$.}
\label{fig:ex2error}
\end{figure}
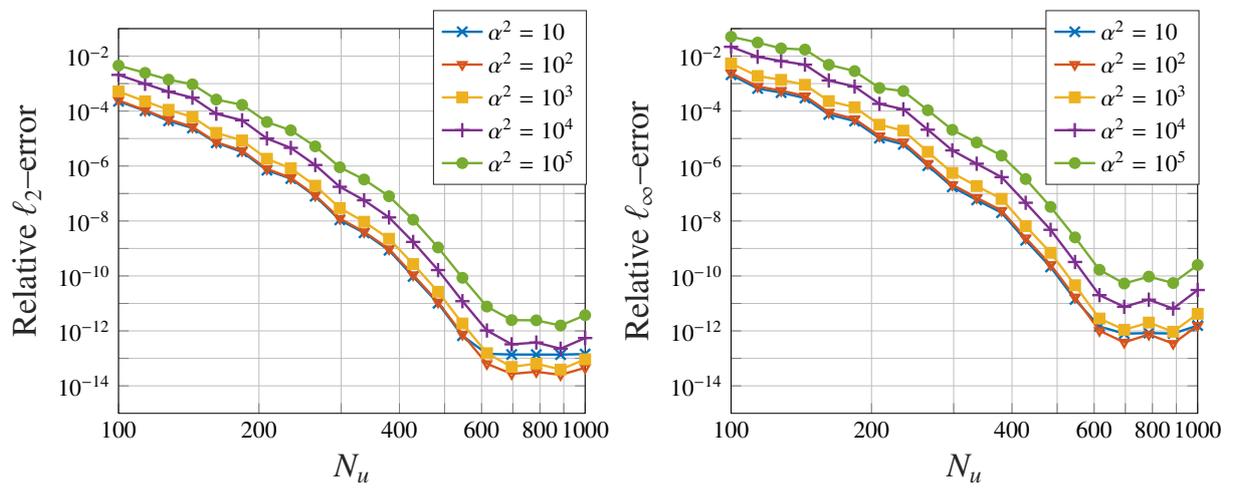

The results in Figure \ref{fig:ex2error}  suggest that \eqref{eq:Ck} is a good estimate for setting $\psi^{k}$ for more complex problems as well. We obtain 10th order convergence with grid refinement. Note that slightly better results can be achieved; the same parameters are used for the entire range of $\alpha$ and are therefore potentially not optimal. As in the previous example the relative $\ell_{\infty}$-error is about two orders of magnitude larger than the relative $\ell_2$-error. In Figure \ref{fig:e2_errorplot} the largest error is by the rightmost point in $\Omega$, at the intersection of two panels. The special purpose quadrature is know to struggle with maintaining full accuracy in such situations.

The modified Helmholtz equation becomes significantly harder to solve for increasing $\alpha^{2}$. This is due to the rapid decay of the kernel \eqref{eq:kerM}, which requires a very fine resolution of the boundary to be resolved. We also suffer from cancellation errors due to the scaling of terms with $\alpha$ or $\alpha^{-1}$. Still, this is not alarming, as an relative $\ell_{\infty}$-error of about $10^{-10}$ can still be obtained for $\alpha^{2} = 10^{5}$. In terms of the heat equation this corresponds to a time step of about $10^{-5}$.
\FloatBarrier

%

\FloatBarrier

\subsection{Example $3$: Adaptive time stepper}
We now test the solver for the heat equation \eqref{eq:HeatEq}--\eqref{eq:HeatEqBC} by setting a tolerance for the time stepping error and investigate if it can be maintained for different resolutions $N_{u}$. For this purpose we use the IMEXRK34 scheme with an adaptive time stepper, see \ref{sss:discretisation_temporal_adaptivity} and \ref{sss:discretisation_temporal_imexrk4}. The smaller the time step, the harder the modified Helmholtz equation is to solve, as concluded above. Thus a high order time stepping scheme, such as IMEXRK34 of fourth order, is a suitable choice since larger time steps can be used. However, other time marching methods can be used as well.

The domain and all parameters are set as for the previous experiment. The heat equation \eqref{eq:HeatEq}--\eqref{eq:HeatEqBC} is solved with right hand side $F$, initial condition and Dirichlet boundary data prescribed by the analytical solution
\begin{equation}
  \label{eq:ex3}
  U(t,x,y) = \exp(-t) \sin((x \cos(\pi/4) + y \sin(\pi/4)))
  + \cos\left(20\sqrt{x^2 + y^2}\right),
\end{equation}
where the time ranges from $0$ to $1$. For the evaluation grid we set $N_{\text{eval}} = N_{u}$ and measure the error at terminal time $t=1$.

\begin{figure}
  \centering
%
%
\definecolor{mycolor1}{rgb}{0.00000,0.44700,0.74100}%
\definecolor{mycolor2}{rgb}{0.85000,0.32500,0.09800}%
\definecolor{mycolor3}{rgb}{0.92900,0.69400,0.12500}%
\definecolor{mycolor4}{rgb}{0.49400,0.18400,0.55600}%
\definecolor{mycolor5}{rgb}{0.46600,0.67400,0.18800}%
\definecolor{mycolor6}{rgb}{0.00000,0.74902,0.74902}%
\begin{tikzpicture}

\begin{axis}[%
width=4.521in,
height=3.566in,
at={(0.758in,0.481in)},
scale only axis,
xmin=80,
xmax=480,
xminorticks=true,
xlabel style={font=\Large\color{white!15!black}},
xlabel={$N_{u}$},
ymode=log,
ymin=1e-11,
ymax=0.01,
yminorticks=true,
ylabel style={font=\Large\color{white!15!black}},
ylabel={Relative $\ell_{2}$--error},
axis background/.style={fill=white},
legend style={legend cell align=left, align=left, draw=white!15!black}
]
%

\addplot [color=mycolor1, line width=1.0pt, mark=x, mark options={solid, mycolor1}, mark size=3.0pt]
  table[row sep=crcr]{%
98	5.21759020112352e-05\\
120	1.61544945743669e-05\\
146	4.13382576428522e-06\\
178	6.97823357826669e-07\\
216	1.75915203660487e-07\\
264	1.77155322985958e-07\\
322	1.48937280909132e-07\\
394	1.48256748343286e-07\\
480	1.48728185750123e-07\\
};
\addlegendentry{Tolerance $10^{-6}$}

\addplot [color=mycolor2, line width=1.0pt, mark=triangle, mark options={solid, rotate=180, mycolor2}, mark size=2.0pt]
  table[row sep=crcr]{%
178	7.21063537005957e-07\\
216	8.11499562115625e-08\\
264	3.16562793493161e-08\\
322	2.54822098774081e-08\\
394	2.5143979339697e-08\\
480	2.54037749861316e-08\\
};
\addlegendentry{Tolerance $10^{-7}$}

\addplot [color=mycolor3, line width=1.0pt, mark=square*, mark options={solid, mycolor3}, mark size=2.0pt]
  table[row sep=crcr]{%
216	7.56981156455107e-08\\
264	9.67550831492913e-09\\
322	2.35517306540817e-09\\
394	2.41387069445615e-09\\
480	2.32143917722796e-09\\
};
\addlegendentry{Tolerance $10^{-8}$}

\addplot [color=mycolor4, line width=1.0pt, mark=+, mark options={solid, mycolor4}, mark size=3.0pt]
  table[row sep=crcr]{%
264	1.10068943028501e-08\\
322	4.05054711193076e-10\\
394	2.06029257570879e-10\\
480	1.92443406386903e-10\\
};
\addlegendentry{Tolerance $10^{-9}$}

\addplot [color=mycolor5, line width=1.0pt, mark=*, mark options={solid, mycolor5}, mark size=2.0pt]
  table[row sep=crcr]{%
394	2.35084700265429e-11\\
480	2.11015519105861e-11\\
};
\addlegendentry{Tolerance $10^{-10}$}

\addplot [color=red]
  table[row sep=crcr]{%
80	1e-06\\
98	1e-06\\
120	1e-06\\
146	1e-06\\
178	1e-06\\
216	1e-06\\
264	1e-06\\
322	1e-06\\
394	1e-06\\
480	1e-06\\
};

\addplot [color=red]
  table[row sep=crcr]{%
80	1e-07\\
98	1e-07\\
120	1e-07\\
146	1e-07\\
178	1e-07\\
216	1e-07\\
264	1e-07\\
322	1e-07\\
394	1e-07\\
480	1e-07\\
};

\addplot [color=red]
  table[row sep=crcr]{%
80	1e-08\\
98	1e-08\\
120	1e-08\\
146	1e-08\\
178	1e-08\\
216	1e-08\\
264	1e-08\\
322	1e-08\\
394	1e-08\\
480	1e-08\\
};

\addplot [color=red]
  table[row sep=crcr]{%
80	1e-09\\
98	1e-09\\
120	1e-09\\
146	1e-09\\
178	1e-09\\
216	1e-09\\
264	1e-09\\
322	1e-09\\
394	1e-09\\
480	1e-09\\
};

\addplot [color=red]
  table[row sep=crcr]{%
80	1e-10\\
98	1e-10\\
120	1e-10\\
146	1e-10\\
178	1e-10\\
216	1e-10\\
264	1e-10\\
322	1e-10\\
394	1e-10\\
480	1e-10\\
};

\end{axis}
\end{tikzpicture}%
%
%
\definecolor{mycolor1}{rgb}{0.00000,0.44700,0.74100}%
\definecolor{mycolor2}{rgb}{0.85000,0.32500,0.09800}%
\definecolor{mycolor3}{rgb}{0.92900,0.69400,0.12500}%
\definecolor{mycolor4}{rgb}{0.49400,0.18400,0.55600}%
\definecolor{mycolor5}{rgb}{0.46600,0.67400,0.18800}%
\definecolor{mycolor6}{rgb}{0.00000,0.74902,0.74902}%
\begin{tikzpicture}

\begin{axis}[%
width=4.521in,
height=3.566in,
at={(0.758in,0.481in)},
scale only axis,
xmin=80,
xmax=480,
xminorticks=true,
xlabel style={font=\Large\color{white!15!black}},
xlabel={$N_{u}$},
ymode=log,
ymin=1e-11,
ymax=0.01,
yminorticks=true,
ylabel style={font=\Large\color{white!15!black}},
ylabel={Relative $\ell_{\infty}$--error},
axis background/.style={fill=white},
legend style={legend cell align=left, align=left, draw=white!15!black}
]
%

\addplot [color=mycolor1, line width=1.0pt, mark=x, mark options={solid, mycolor1}, mark size=3.0pt]
  table[row sep=crcr]{%
98	0.000364085717805901\\
120	0.000235252764277243\\
146	5.77488423659335e-05\\
178	8.84229597174475e-06\\
216	2.41310548083693e-06\\
264	4.50025586616405e-07\\
322	3.85698544576379e-07\\
394	3.84173407574805e-07\\
480	3.85251599683235e-07\\
};
\addlegendentry{Tolerance $10^{-6}$}

\addplot [color=mycolor2, line width=1.0pt, mark=triangle, mark options={solid, rotate=180, mycolor2}, mark size=2.0pt]
  table[row sep=crcr]{%
178	9.12590674831361e-06\\
216	1.80310897703536e-06\\
264	1.36991607114933e-07\\
322	6.75571593136641e-08\\
394	6.67531117517696e-08\\
480	6.74316819041888e-08\\
};
\addlegendentry{Tolerance $10^{-7}$}

\addplot [color=mycolor3, line width=1.0pt, mark=square*, mark options={solid, mycolor3}, mark size=2.0pt]
  table[row sep=crcr]{%
216	1.62506730187398e-06\\
264	2.74348546395707e-07\\
322	6.15364981788288e-09\\
394	6.52179835837028e-09\\
480	6.19406810223975e-09\\
};
\addlegendentry{Tolerance $10^{-8}$}

\addplot [color=mycolor4, line width=1.0pt, mark=+, mark options={solid, mycolor4}, mark size=3.0pt]
  table[row sep=crcr]{%
264	3.21005310656022e-07\\
322	5.38054630086699e-09\\
394	5.8093678476999e-10\\
480	5.31940872055256e-10\\
};
\addlegendentry{Tolerance $10^{-9}$}

\addplot [color=mycolor5, line width=1.0pt, mark=*, mark options={solid, mycolor5}, mark size=2.0pt]
  table[row sep=crcr]{%
394	3.43131220619345e-10\\
480	7.48203262649225e-11\\
};
\addlegendentry{Tolerance $10^{-10}$}
\addplot [color=red]
  table[row sep=crcr]{%
80	1e-06\\
98	1e-06\\
120	1e-06\\
146	1e-06\\
178	1e-06\\
216	1e-06\\
264	1e-06\\
322	1e-06\\
394	1e-06\\
480	1e-06\\
};

\addplot [color=red]
  table[row sep=crcr]{%
80	1e-07\\
98	1e-07\\
120	1e-07\\
146	1e-07\\
178	1e-07\\
216	1e-07\\
264	1e-07\\
322	1e-07\\
394	1e-07\\
480	1e-07\\
};

\addplot [color=red]
  table[row sep=crcr]{%
80	1e-08\\
98	1e-08\\
120	1e-08\\
146	1e-08\\
178	1e-08\\
216	1e-08\\
264	1e-08\\
322	1e-08\\
394	1e-08\\
480	1e-08\\
};

\addplot [color=red]
  table[row sep=crcr]{%
80	1e-09\\
98	1e-09\\
120	1e-09\\
146	1e-09\\
178	1e-09\\
216	1e-09\\
264	1e-09\\
322	1e-09\\
394	1e-09\\
480	1e-09\\
};

\addplot [color=red]
  table[row sep=crcr]{%
80	1e-10\\
98	1e-10\\
120	1e-10\\
146	1e-10\\
178	1e-10\\
216	1e-10\\
264	1e-10\\
322	1e-10\\
394	1e-10\\
480	1e-10\\
};

\end{axis}
\end{tikzpicture}%
\caption{The errors for different resolutions of the uniform grid, at terminal time $t=1$. The red lines are the set tolerances for the relative $\ell_2$-error, used by adaptive time stepper.}
\label{fig:imexrk34tol}
\end{figure}
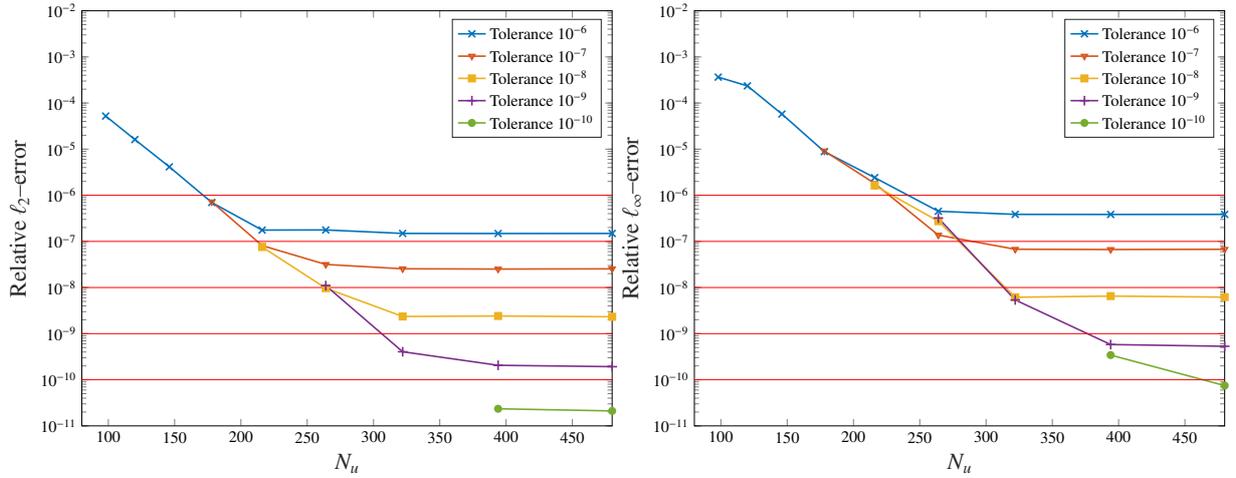

In Figure \ref{fig:imexrk34tol} the red lines correspond to set tolerances. It is clear that the adaptive time stepper works as intended, even for tolerances down to $10^{-10}$. The relative $\ell_{\infty}$-error is more sensitive to the resolution and exceeds the set tolerance earlier in terms of spatial resolution, roughly with one order in magnitude.
\FloatBarrier

%

\FloatBarrier

\subsection{Example $4$: The Allen-Cahn equation, a reaction diffusion problem}
\begin{figure}
  \centering
  \includegraphics[width=0.49\textwidth]{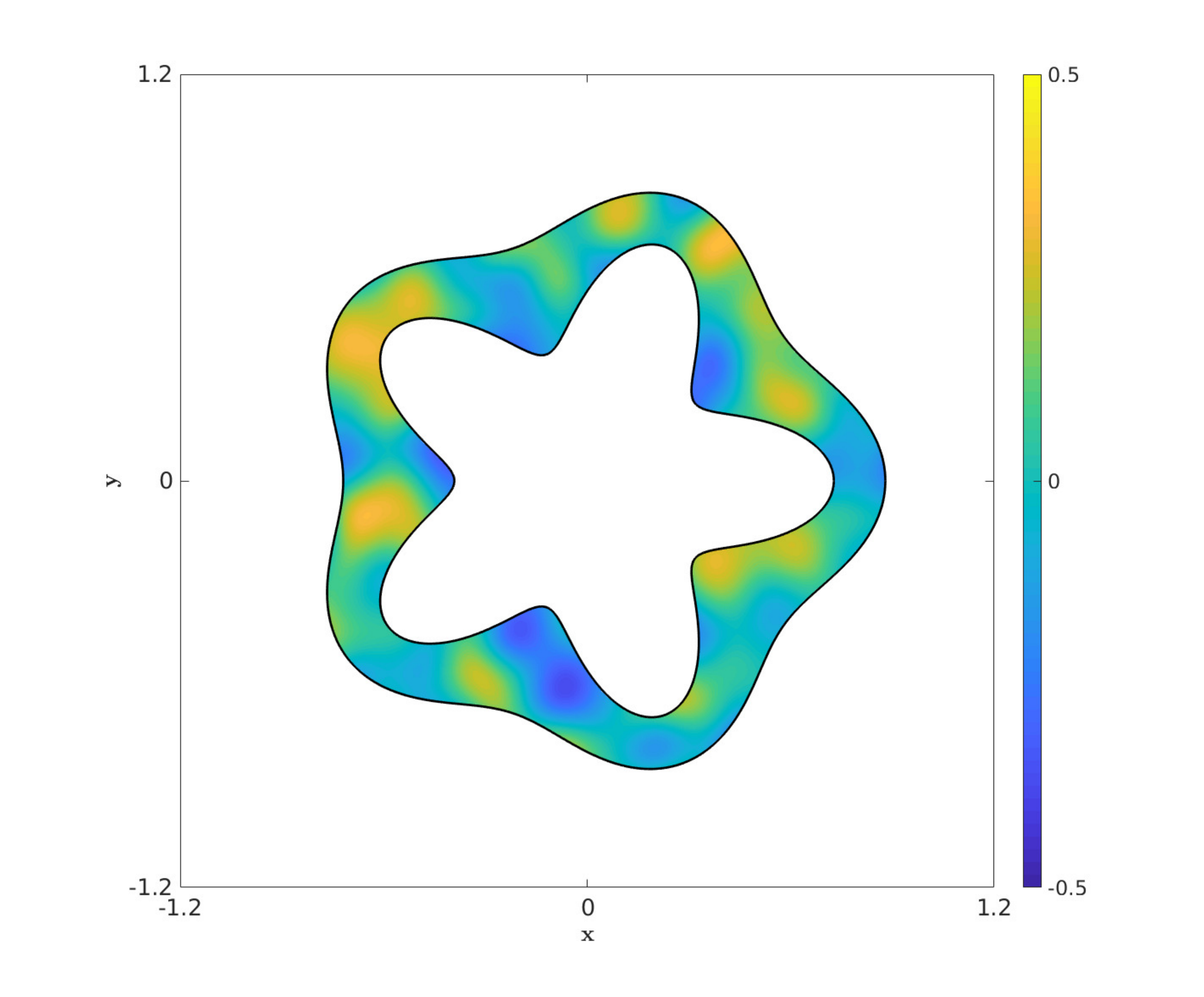}
    \includegraphics[width=0.49\textwidth]{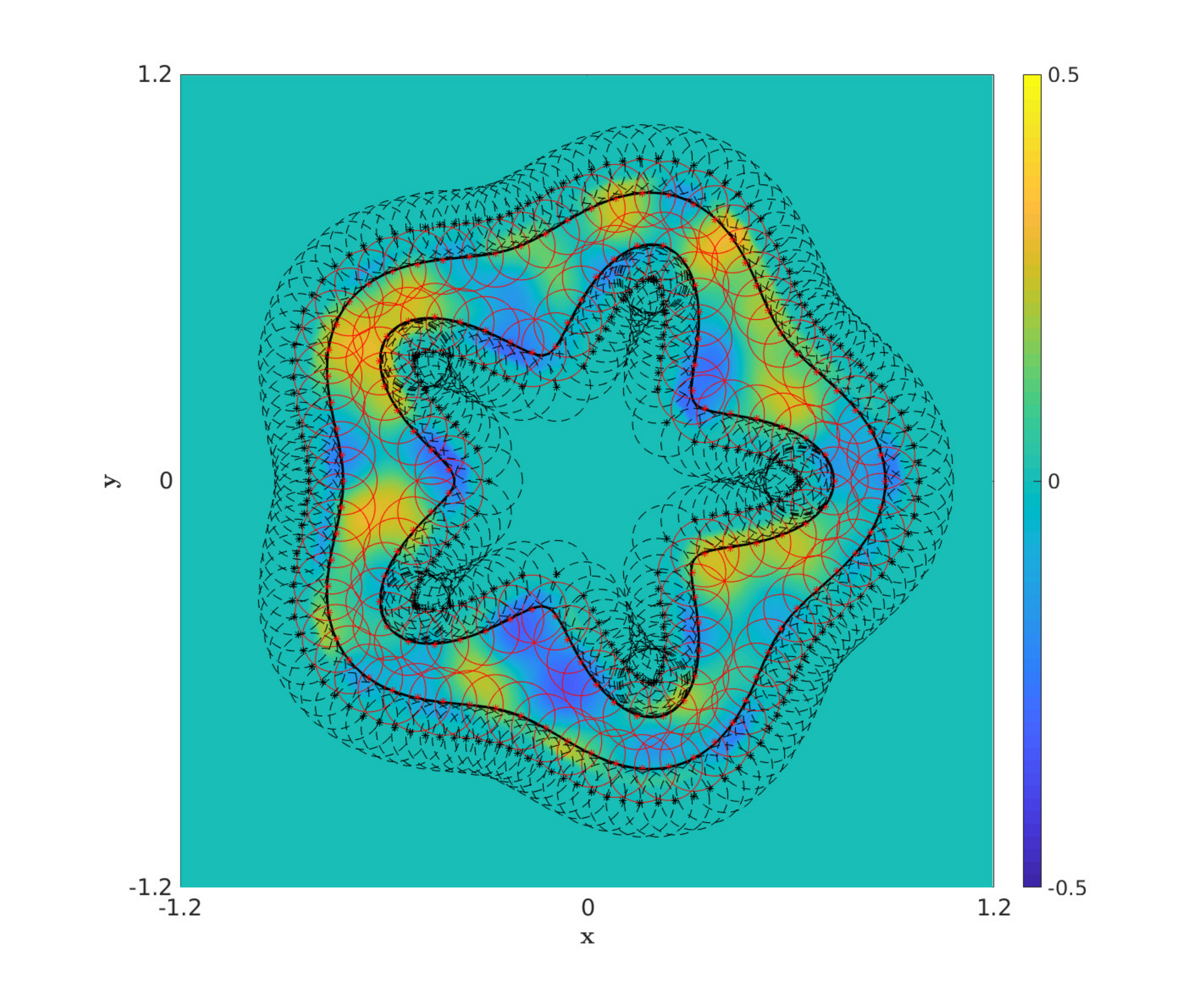}
\caption{Left: The initial data $U_0$ \eqref{eq:allencahnEqIC}. Right: The right hand side of \eqref{eq:allencahnEq} at $t_0$, extended with PUX. Black corresponds to zero partitions and red to interpolation partitions. Note that to increase visibility of the field a different scaling is used than for \ref{fig:t1}--\ref{fig:t6}}
\label{fig:ex4_u0fe}
\end{figure}

The Allen-Cahn equation is stated as
\begin{align}
  \frac{\partial U(t,\mathbf{x})}{\partial t} - C\Delta U(t,\mathbf{x}) &= U(t,\mathbf{x})(1-U(t,\mathbf{x})^2),\quad t_{0}<t,\quad \mathbf{x}\in \Omega\subset\mathbb{R}^{2},\label{eq:allencahnEq}\\
  U(t_{0},\mathbf{x}) &= U_{0}(\mathbf{x}),\quad \mathbf{x}\in  \Omega,\label{eq:allencahnEqIC}\\
  U(t,\mathbf{x}) &= e^{-t/2}U_{0}(\mathbf{x}),\quad \mathbf{x}\in \Gamma,\label{eq:allencahnEqBC}
\end{align}
with $C = 10^{-3}$. The right hand side of \eqref{eq:allencahnEq} is nonlinear and has three stationary points: $U = -1,\, 0,\,1$. For randomised initial data the solution creates over time patterns with zones attaining these values. The initial data is not entirely randomised, since we need smoothness to discuss convergence and accuracy. Instead, we create smooth data by uniformly distributing $50$ Gaussians \eqref{eq:GaussianRBF} with $\varepsilon = 10$ over the computational domain with $L = 1.2$. Each Gaussian is assigned a coefficient drawn randomly from a uniform distribution over $-0.5$ to $0.5$. The partition size $R$ is set to $0.1$; the domain, the extended right hand and the distribution of partitions are shown in Figure \ref{fig:ex4_u0fe}. Each boundary component is discretised with $80$ panels.

\begin{figure}[ht]
  \begin{subfigure}[b]{0.35\linewidth}
    \centering
    \includegraphics[width=\linewidth]{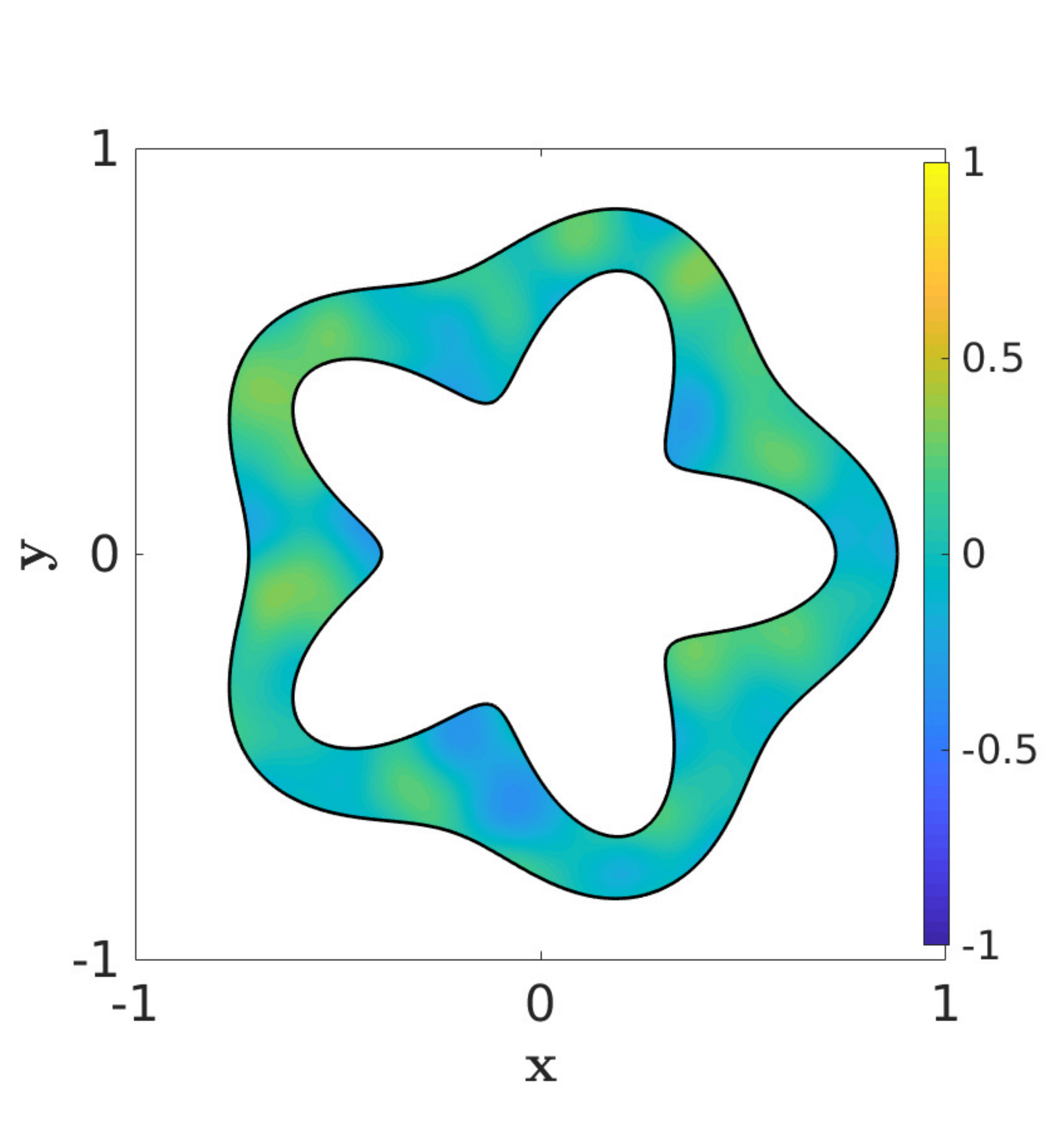}
    \caption{$t\sim 0.005$}
    \label{fig:t1}
  \end{subfigure}
  \begin{subfigure}[b]{0.35\linewidth}
    \centering
    \includegraphics[width=\linewidth]{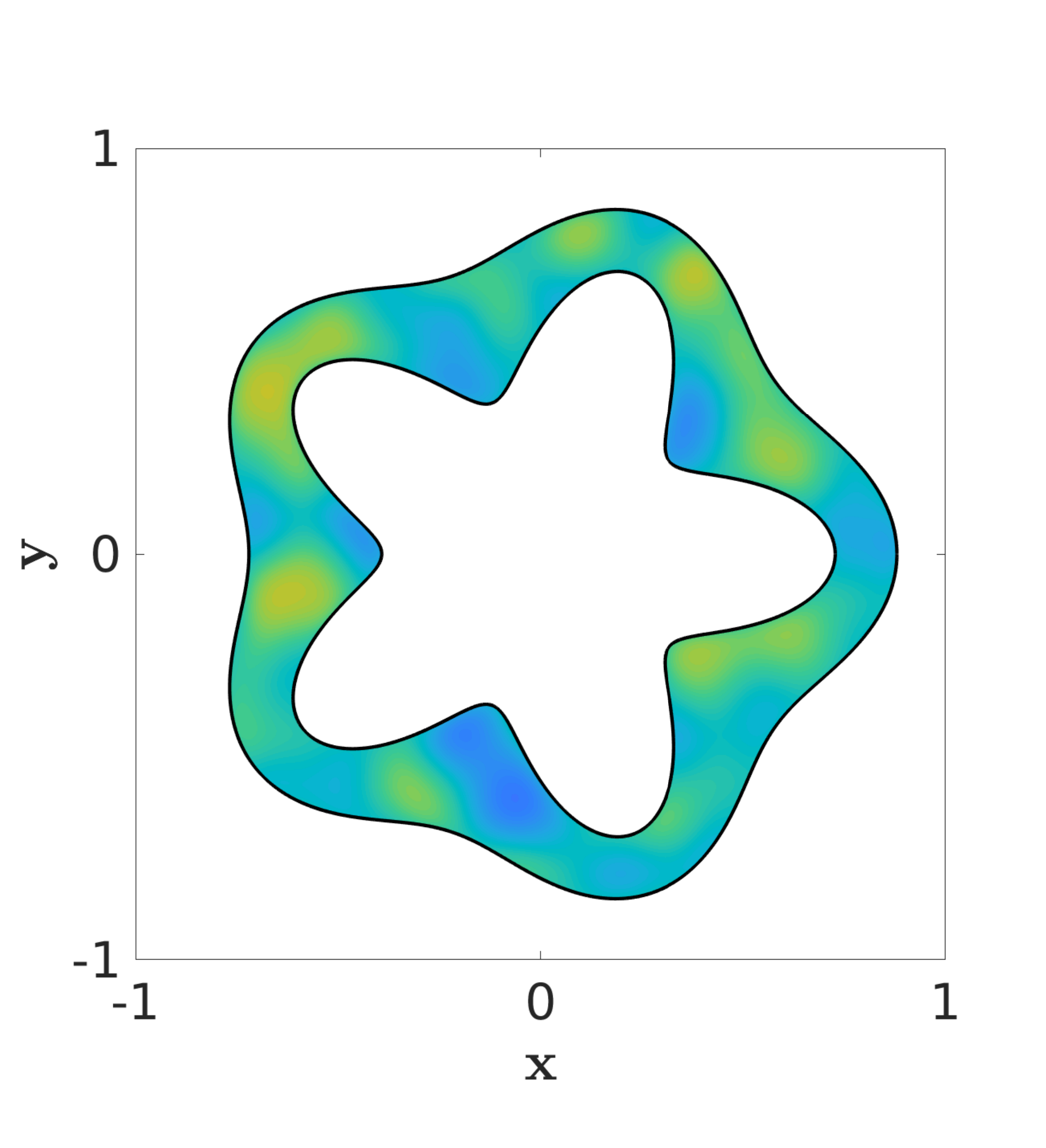}
    \caption{$t\sim 0.5$}
    \label{fig:t2}
  \end{subfigure}
  \begin{subfigure}[b]{0.35\linewidth}
    \centering
    \includegraphics[width=\linewidth]{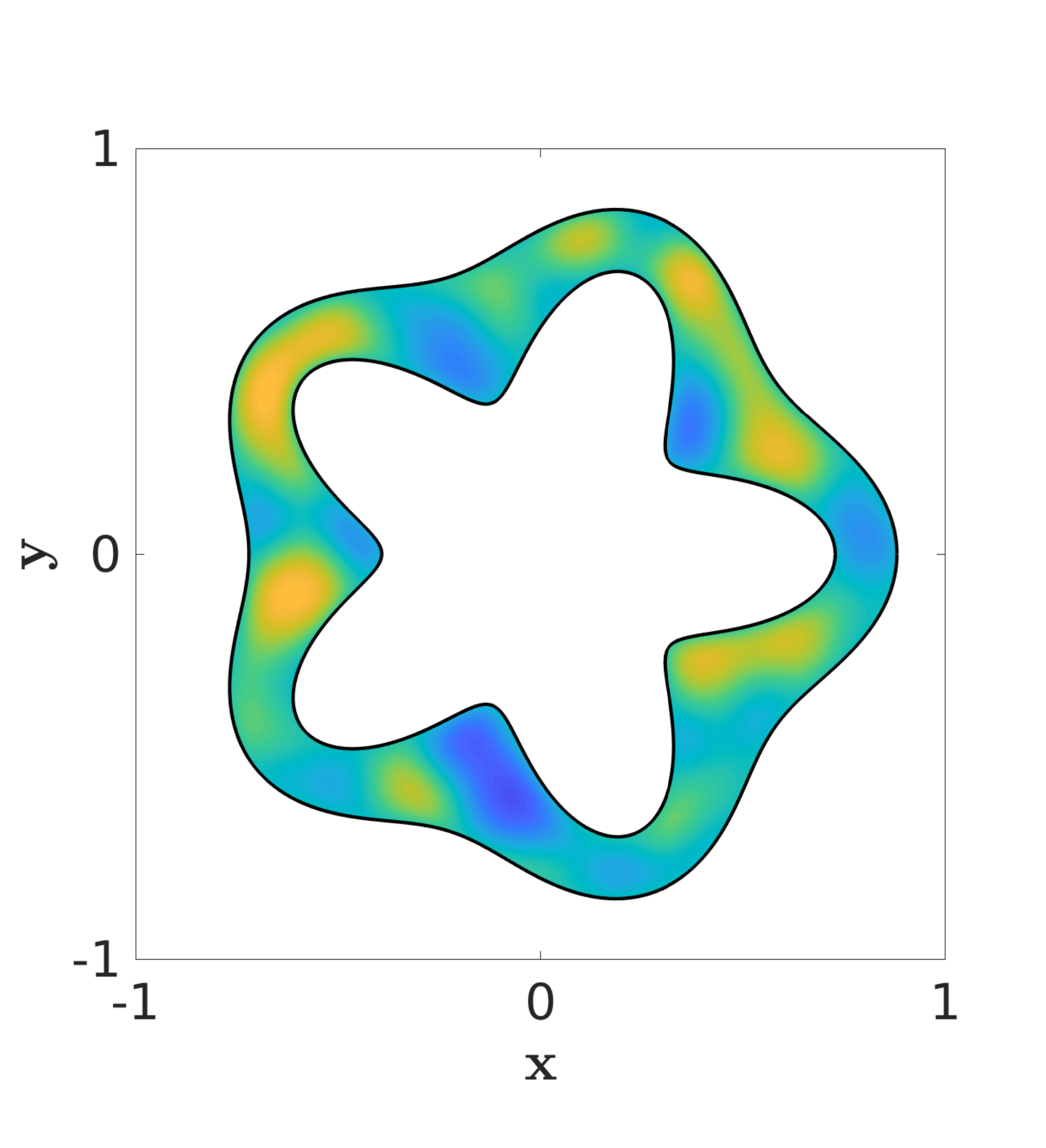}
    \caption{$t\sim 1.3$}
    \label{fig:t3}
  \end{subfigure}
  \linebreak
  \begin{subfigure}[b]{0.35\linewidth}
    \centering
    \includegraphics[width=\linewidth]{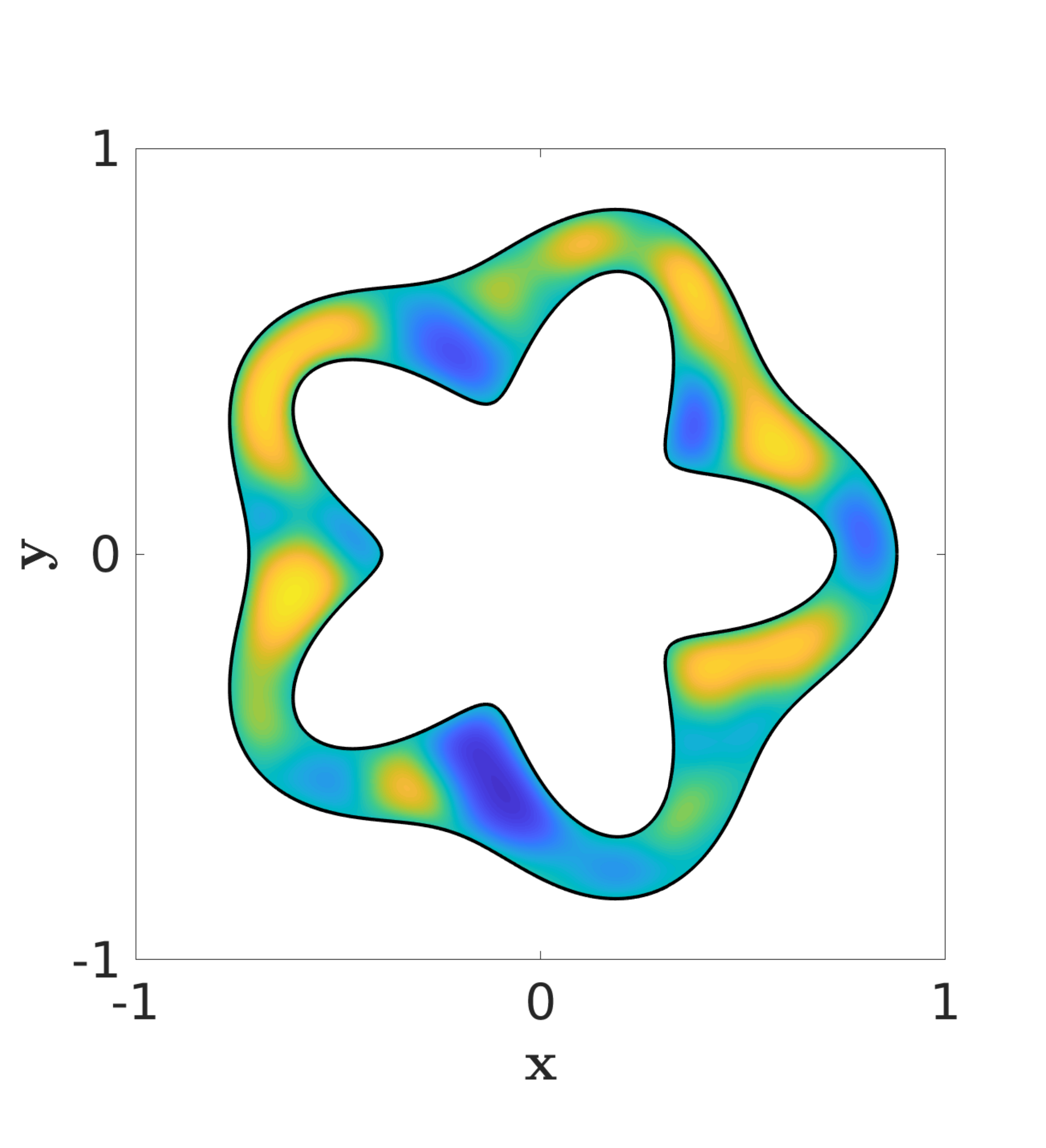}
    \caption{$t\sim 2.5$}
    \label{fig:t4}
  \end{subfigure}
  \begin{subfigure}[b]{0.35\linewidth}
    \centering
    \includegraphics[width=\linewidth]{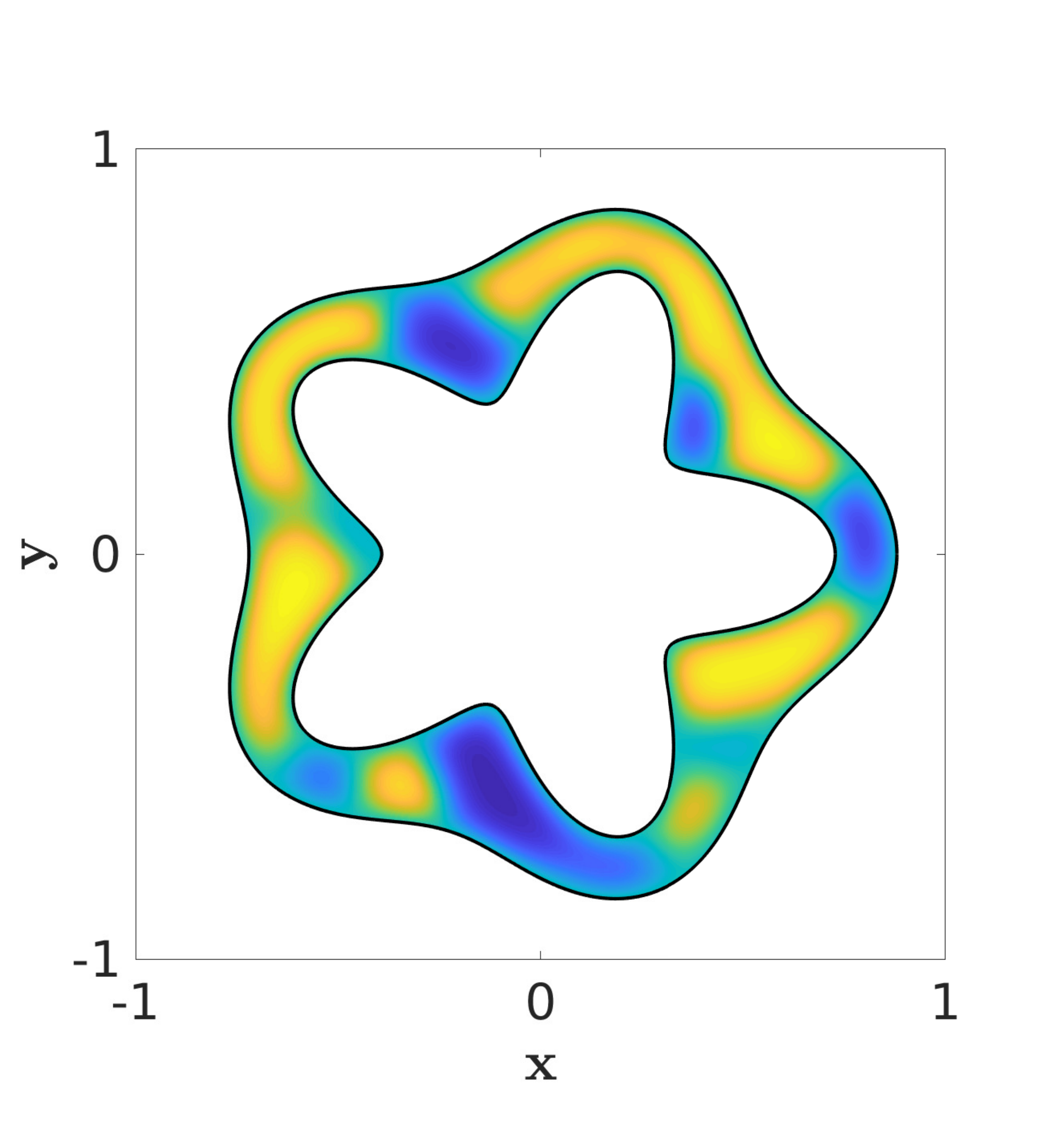}
    \caption{$t\sim 4.2$}
    \label{fig:t5}
  \end{subfigure}
  \begin{subfigure}[b]{0.35\linewidth}
    \centering
    \includegraphics[width=\linewidth]{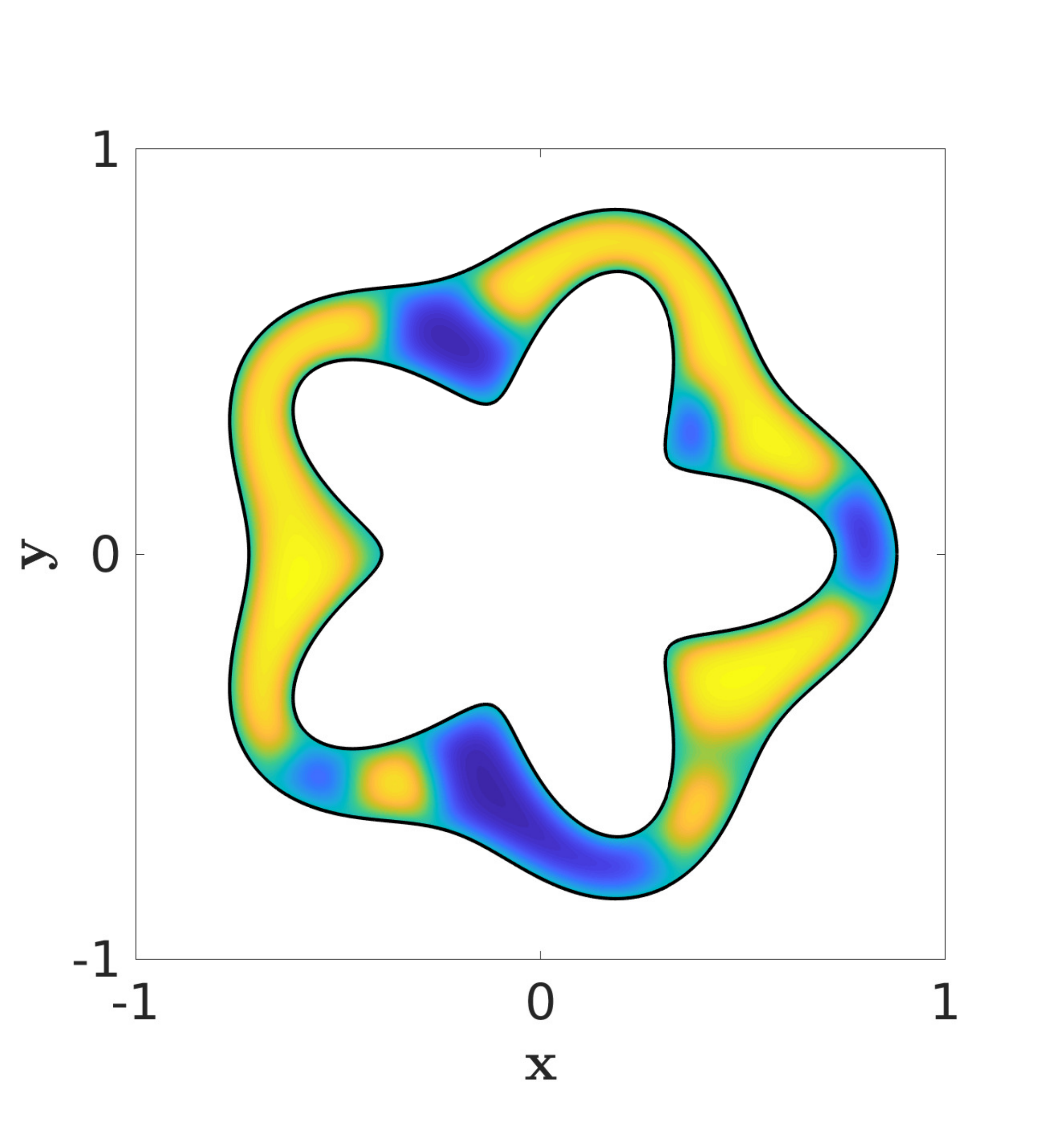}
    \caption{$t=6$}
    \label{fig:t6}
  \end{subfigure}

  \caption{Numerical solution to \eqref{eq:allencahnEq} with $N_u = 800$ and tolerance $10^{-6}$ at $t\sim 0.005,\,0.5,\,1.3,\,2.5,\,4.2$ and terminal time $t = 6$.}
  \label{fig:ex4_snapshots}
\end{figure}

\begin{figure}
  \centering
      \includegraphics[width=0.40\textwidth]{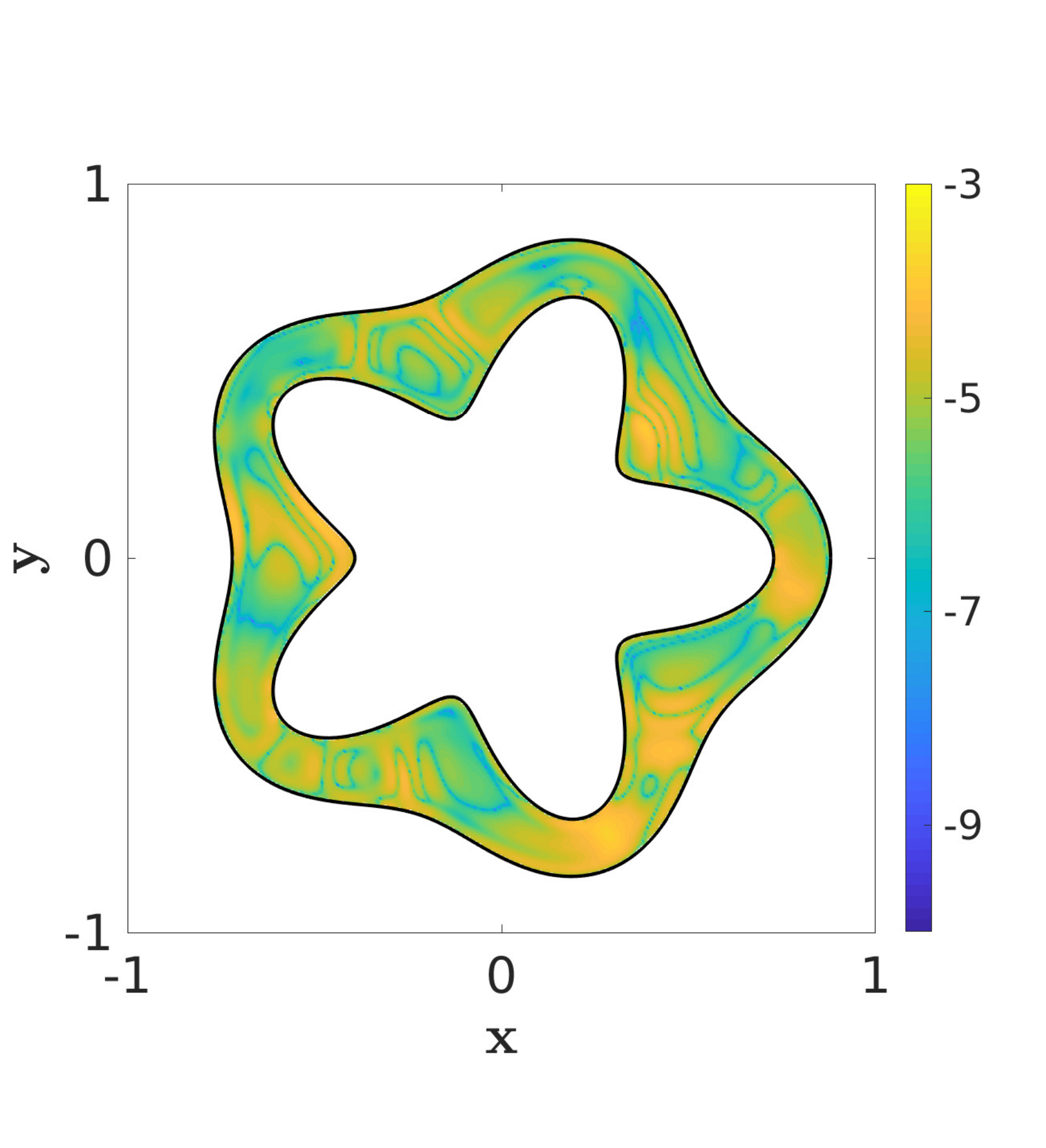}
      \includegraphics[width=0.40\textwidth]{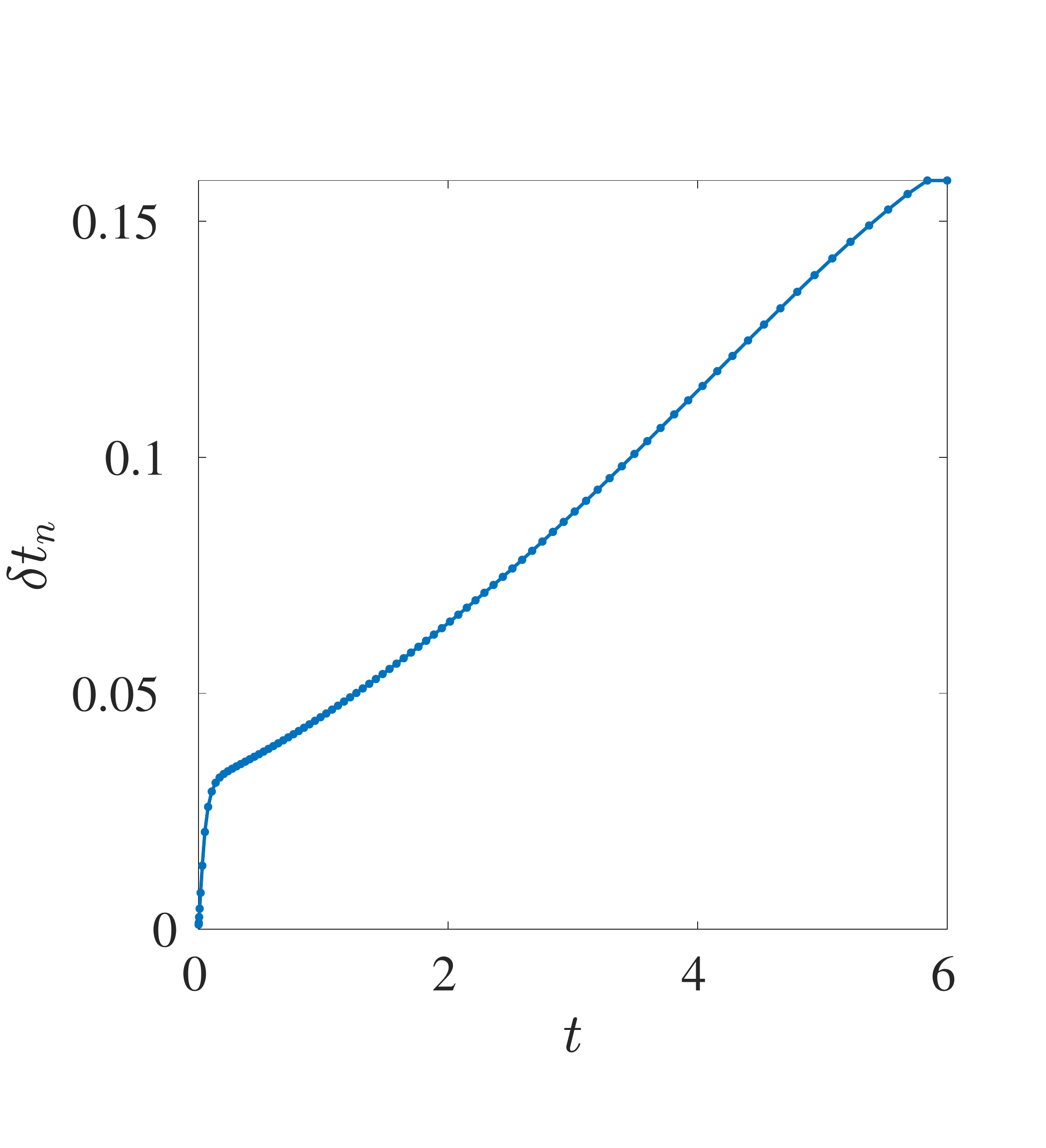}
\caption{Left: Pointwise relative error for $N_u = 400$, tolerance $10^{-5}$. Right: Evolution of time step $\delta t$ over time.}
\label{fig:ex4_erranddt}
\end{figure}

We create a reference solution by solving the Allen-Cahn equation with tolerance $10^{-6}$ and $N_u = 800$, from time $0$ to $6$. The errors are measured on grids with $N_{\text{eval}}=200,\,400$ at terminal time $t=6$. Snapshots of this solution are shown in Figures \ref{fig:t1} to \ref{fig:t6}. Indeed the solution forms a pattern of patches with the values $-1$, $0$ and $1$. The results are shown in Table \ref{tab:ex4error}. For $N_u = 400$ the relative $\ell_2$-error stays under the set tolerance. However, unlike example $3$ the relative $\ell_{2}$-error is always a factor ten larger than the set tolerance. For $N_u = 200$
 only the tolerance $10^{-3}$ can be obtained. Clearly this resolution is insufficient to resolve the spatial problem more accurately than that. The error at the terminal time $t=6$ for $N_u = 400$ with tolerance $10^{-5}$ is shown in Figure
\ref{fig:ex4_erranddt}. In this figure the evolution of the time step is also shown; as the solution advances in time the time step becomes larger. Initially it grows faster, compared to later, as the initial time step was intentionally set small.

\begin{table}
\centering
\small
\caption{Errors at terminal time $t = 6$ for solving \eqref{eq:allencahnEq}--\eqref{eq:allencahnEqBC} for resolutions $N_u = 200,\,400$ and set tolerances for the adaptive time stepper. The reference solution is computed with a tolerance of $10^{-6}$ with $N_u = 800$. The resolution $N_u=200$ is insufficient to reach errors below $10^{-3}$, while the relative $\ell_2$-error for $N_u = 400$ satisfies tolerances $10^{-3}$, $10^{-4}$ and $10^{-5}$.}
\begin{threeparttable}
{\def\arraystretch{1.3}
\begin{tabular}{ccccc}
  \toprule
& \multicolumn{2}{c}{\textbf{Relative max error}}&\multicolumn{2}{c}{\textbf{Relative $\ell_2$ error}} \\
\midrule
\textbf{Tolerance} & $N_u = 200$ & $N_u = 400$ & $N_u = 200$ & $N_u = 400$ \\
\midrule
$10^{-3}$ & $3.7314\times 10^{-3}$ & $3.6123\times 10^{-3}$ & $7.2035\times 10^{-4}$ & $7.0625\times 10^{-4}$ \\
$10^{-4}$ & $2.6619\times 10^{-2}$ & $3.6519\times 10^{-4}$ & $6.6734\times 10^{-3}$ & $6.8990\times 10^{-5}$ \\
$10^{-5}$ & $2.1136\times 10^{-3}$ & $3.9993\times 10^{-5}$ & $3.1405\times 10^{-4}$ & $9.3987\times 10^{-6}$ \\
\bottomrule
\end{tabular}
}
\end{threeparttable}
\label{tab:ex4error}
\end{table}

%

\FloatBarrier

\section{Conclusions}
\label{s:conclusions}
We present a framework built around a panel-based Nystr\"{o}m boundary integral method for solving the forced isotropic heat equation in two dimensions, on multiply connected complex domains.  We have addressed several of the issues listed in \cite{Kropinski2011Heat}, thereby increasing the class of solvable problems as well as the accuracy in the solutions.

We show how any IMEX method can be applied as time stepping scheme, and employ an adaptive fourth order Runge-Kutta scheme in our examples, to accurately solve the heat equation as well as the Allen-Cahn equation, a reaction-diffusion problem with a nonlinear forcing term.
Regardless of the specific details of the chosen method, a time step in solving the heat equation is reduced to solving one or a sequence, for a multi-stage method, of modified Helmholtz equations.

As in \cite{Kropinski2011modHelm} we formulate the modified Helmholtz equation as a boundary integral problem. Utilising the linearity of the differential operator, the solution is split into a particular- and  homogeneous  problem. Solving the former to high accuracy relies on extending the given right hand side from the domain it is given on to the entire plane. It is achieved with a partition of unity extension (PUX), that only requires known data at uniform point locations inside the domain. The extension that is computed on a uniform grid in a rectangular domain has compact support and a specified global regularity, making spectral methods very efficient and simple to use. We confirm that the various parameters for PUX, in the context of the modified Helmholtz equation, indeed can be set as for the Poisson equation in \cite{FRYKLUNDPUX}.  This yields an automated selection for the global regularity to balance different errors, leading to a method which converges with an order $10$ in the grid size.

A panel-based Nystr\"{o}m boundary integral method is used to solve the homogeneous problem with modified Dirichlet data, such that the total solution is the sum of the particular- and homogeneous solution. The boundary values of the particular solution are computed using a non-uniform FFT. For evaluation of singular and nearly singlar integrals, we have introduced a methodology based on product integration and an explicit kernel split that has given highly accurate results for the Helmholtz \cite{Helsing2015} and Stokes equations \cite{OJALA2015145}.  For large $\alpha$ (small time steps), the method in its original form, would however fail completely if an unfeasibly high upsampling of the boundary was not applied.
We however realized that this upsampling is only needed very locally, and
developed an adaptive approach \cite{2019arXiv190607713K} to achive a computationally efficient method with high accuarcy.

In total, these developments yields a method for very accurately solving the heat equation on comlex domains. The highest attainable accuracy in the solution of the modified Helmholtz equation does show a weak depenence on $\alpha$, but even for the largest values, solutions can typically be attained with at least ten correct digits, meaning that strict time stepping tolerances for the heat equation can be satisfied.

In terms of future developments, it would be useful for some problems to replace the uniform grids and FFT-based method for the particular solution with a volume potential evaluation based on an adaptive FMM. This would however need an integration of the PUX method into the adaptive procedure. Another development is to consider the solution of the heat equation and the closely related advection-diffusion equation on time-dependent domains. The motivation for this is the need to solve such an equation for the concentration of surfactants in the oil-phase of a micro-system with water drops in oil. These surfactants, or surface active agents, have an exchange with surfactants on the drop surfaces, that alters the surface tension of the drop.  Numerical methods for simulating surfactant advection and diffusion on the boundary of drops have been understood and implemented successfully, see \cite{PALSSON2019218,KROPINSKI20114466}. An important  extension would be to allow also for surfactants in the oil-phase. One strength of these methods is the accurate treatment of interface conditions, something that is absolutely essential at these small scales where the interface dynamics is of key importance.

\section{Acknowledgements}
\label{s:acknowl}

We thankfully acknowledge the support of the Swedish Research Council
under Grant No. $2015$-$04998$ and funding from the G\"{o}ran Gustafsson
Foundation for Research in Natural Sciences and Medicine. We are humbly grateful for the support from the
Natural Science and Engineering Research Council of Canada.

\appendix
\section{Adaptive time-stepping with IMEX Runge-Kutta methods}

\subsection{Adaptive discretisation in time}
This appendix shows how applying implicit-explicit Runge-Kutta (IMEXRK) schemes from \citep{KENNEDY2003139} to the heat equation reduces it to a sequence of modified Helmhotlz equations to solve at each time step.
\label{ss:discretisation_temporal}
Formulate the heat equation \eqref{eq:HeatEq}--\eqref{eq:HeatEqBC} as
\begin{align}
  \label{eq:IMEXgeneral}
  \frac{\partial U(t,\x)}{\partial t} =  F^{I}(t,\x,U) + F^{E}(t,\x,U),\quad \x\in\Omega, \\
  \label{eq:IMEXFIFE}
  F^{I}(t,\x,U) = \Delta U(t,\x),\quad F^{E}(t,\x,U) = F(t,\x),
\end{align}
where the superscripts denote implicit and explicit, referring to the term being classified as stiff or nonstiff, respectively.

Let $t_{N}$ denote an instance in time that is the sum of previous discrete time steps $\{\delta t_{i}\}_{i = 1}^{N}$ that may be of different size:
\begin{equation}
 \label{eq:timedisc}
t_{N} = \sum\limits_{i=1}^{N}\delta t_{i} + t_{0},
\end{equation}
for some initial time $t_{0}$. Let $U_{N}$ be the approximation of $U(t_{N})$, then the approximated solution at time $t_{N+1}$ is
\begin{equation}
 \label{eq:imexsol}
   U_{N+1} = U_{N} + \delta t_{N+1}\sum_{\sigma \in \{I,E\}} \sum_{j = 1}^{N_{S}} b_{j}^{\sigma}k_{j}^{\sigma},
\end{equation}
where $N_{S}$ is the number of stages for $k^{\sigma}$, $\sigma\in \{I,E\}$, computed as

\begin{equation}
 \label{eq:stagek}
 k_{i}^{\sigma} = F^{\sigma}\left(t_{N}+\delta t_{N+1}c_{i}^{\sigma},\bar{U}_{i}\right),\quad i = 1,\ldots,N_{S}.
\end{equation}
The second argument of $F^{\sigma}$ in \eqref{eq:stagek} is defined as
\begin{equation}
 \label{eq:Ubar}
 \bar{U}_{i} = U_{N} + \delta t_{N+1}\sum_{\sigma \in \{\text{I, E}\}}\sum_{j = 1}^{i} a_{i,j}^{\sigma}k_{j}^{\sigma} = U_{N} + \delta t_{N+1}\sum_{\sigma \in \{\text{I, E}\}}\sum_{j = 1}^{i-1} a_{i,j}^{\sigma}k_{j}^{\sigma} + \delta t_{N+1}a_{i,i}^{I}k_{i}^{I}, \quad i>1,
\end{equation}
and $\bar{U}_{1} = U_{N}$. The coefficients $\{a_{i,j}^{\sigma}\}_{i,j = 1}^{N_{S}}$, $\{b_{j}^{\sigma}\}_{j = 1}^{N_{S}}$ and $\{c_{i}^{\sigma}\}_{i= 1}^{N_{S}}$ are tabulated in the two associated Butcher tableaus for $\sigma = I$
 and $\sigma =E$, see Table \ref{tab:butcher} for a general IMEXRK scheme. The principal difference between the coefficients for implicit and explicit methods is that $a_{i,j}^{E} = 0$ for $i\leq j$ while $a_{i,j}^{I} \neq 0$ for $i=j$, excluding $i = 1$.  The quantity $\bar{U}_{i}$ is unknown for every $i=2,\ldots,N_{S}$, since the corresponding implicit stage $k^{I}_{i}$ is unknown.

 The implicit stage at $i$ is $k^{I}_{i} = F^{I} = \Delta \bar{U}_{i}$ by definition \eqref{eq:IMEXFIFE}. To avoid approximating the differential operator replace $k^{I}_{i}$ in \eqref{eq:Ubar} with $\Delta \bar{U}_{i}$ and reformulate as
 \begin{equation}
  \label{eq:IMEXmodHelm}
  \frac{1}{\delta t_{N+1}a_{i,i}^{I}} \bar{U}_{i} - \Delta \bar{U}_{i} = \frac{1}{\delta t_{N+1}a_{i,i}^{I}}U_{N} + \sum_{\sigma \in \{\text{I, E}\}}\sum_{j = 1}^{i-1} \frac{a_{i,j}^{\sigma}}{a_{i,i}^{I}}k_{j}^{\sigma}.
 \end{equation}
 The idea is to solve for $\bar{U}_{i}$ and since the right hand side is known, $\Delta \bar{U}_i$ can be extracted from the expression above.

The equation \eqref{eq:IMEXmodHelm} has the form of the modified Helmholtz equation \eqref{eq:ModHelmEq}--\eqref{eq:ModHelmEqBC}: $f(\x)$ corresponds to the right hand side, $u(\mathbf{x}) = \bar{U}_{i}(\x)$ and
$\alpha^{2} = (\delta t_{N+1}a_{i,i}^{I})^{-1}$. We stress that $\alpha^{2}\sim (\delta t_{N+1})^{-1}$; the larger $\alpha^{2}$ is the  harder \eqref{eq:ModHelmEq}--\eqref{eq:ModHelmEqBC} is to solve accurately in terms of numerics, see Section \ref{sss:discretisation_adaptive_modhelm_homo_specquad}. The associated boundary condition $g$ is \eqref{eq:HeatEqBC} evaluated at $t_{N}+\delta t_{N+1}c_{i}^{I}$.

To obtain the next stage $k^{I}_{i}$ the equation \eqref{eq:IMEXmodHelm}  must be solved for $\bar{U}_{i}$ in $\Omega$. Once $\bar{U}_{i}$ is known, reformulate \eqref{eq:IMEXmodHelm} and compute
\begin{equation}
 \label{eq:computeimplicitstage}
 k_{i}^{I} = F^{I}(t_{N}+\delta t_{N+1}c_{i}^{\sigma},\bar{U}_{i}) = \Delta \bar{U}_{i} = \frac{1}{\delta t_{N+1}a_{i,i}^{I}} \bar{U}_{i} -  \frac{1}{\delta t_{N+1}a_{i,i}^{I}}U_{N} + \sum_{\sigma \in \{\text{I, E}\}}\sum_{j = 1}^{i-1} \frac{a_{i,j}^{\sigma}}{a_{i,i}^{I}}k_{j}^{\sigma}.
\end{equation}
With $k^{I}_{i}$ known the stage $k^{E}_{i}$, that is $F^{E}$, can be computed explicitly. Note that for \eqref{eq:HeatEq}--\eqref{eq:HeatEqBC} $F^{E} = F(t,x)$, so the explicit stage $k^{E}_{i}$ is independent of the implicit stages, thus it is computed directly. Note that this is not the case if e.g. an advection term $\nabla U$ is added, as it would be included in $F^{E}$. In order to keep the formulation general we think of $F^{E}$ as function of $U$.

To summarise: the approximate solution $U^{N+1}$ at time $t^{N+1}$ is given by \eqref{eq:imexsol}. The stages $k^{I}_{i}$, for $i = 1,\ldots,N_{S}$ are obtained by solving \eqref{eq:ModHelmEq}--\eqref{eq:ModHelmEqBC}, corresponding to \eqref{eq:IMEXmodHelm}, and explicit computation of \eqref{eq:computeimplicitstage}. Once $\bar{U}_{i}$  is known $k_{i}^{E} = F^{E}\left(t_{N}+\delta t_{N+1}c_{i}^{E},\bar{U}_{i}\right)$ is computed explicitly. See the flowchart in \ref{appendix:flowchart} for a graphical overview.
\begin{table}
\centering
{\def\arraystretch{1.3}
\begin{tabular}{c|cccc}
$0$ & $0$ &{} &{} & {}\\
$c^{\sigma}_{2}$ & $a^{\sigma}_{21}$ & $a^{\sigma}_{22}$ & {} & {}\\
$\vdots$ & $\vdots$ & {} &$\ddots$ & {}\\
$c^{\sigma}_{N_{S}}$ & $a^{\sigma}_{N_{S}1}$ &$\cdots$ & $\cdots$ & $a^{\sigma}_{N_{S}N_{S}}$\\
\hline
& $b^{\sigma}_{1}$ & $\cdots$ & $\cdots$ & $b^{\sigma}_{N_{S}}$
\end{tabular}
}
\caption{Coefficients for an IMEXRK scheme, where $\sigma \in \{ I,E\}$, denoting implicit or
explicit, applied to the stiff and nonstiff term, respectively. In general $a^{E}_{ij}=0$ for
 $i\leq j$ and $a^{I}_{ij}\neq 0$ for $i= j$, excluding $i = 1$.}
\label{tab:butcher}
\end{table}

\subsubsection{IMEXRK2}
\label{sss:discretisation_temporal_imexrk2}
\label{ssec:IMEXRK2}
This scheme is never used in this paper, but serves as a simple example of applying an IMEX Runge-Kutta scheme. The stencil for IMEXRK2, with coefficients tabulated in Table \ref{tab:butcherIMEXRK2},  involves taking a half time step $\delta t_{N+1}/2$ and solving for $\bar{U}_{2}$ satisfying
\begin{equation}
  \label{eq:deltauNplushalf}
\frac{2}{\delta t_{N+1}}\bar{U}_{2}(\mathbf{x}) -\Delta \bar{U}_{2}(\mathbf{x}) = \frac{2}{\delta t_{N+1}}U^{N}(\mathbf{x}) + F^{E}(t_{N},\mathbf{x},\bar{U}_{2}),\quad\mathbf{x}\in \Omega.
\end{equation}
By \eqref{eq:imexsol} the solution at the next time-step $t_{N+1} = \delta t_{N+1} + t_{N}$ for every $\x\in \Omega$ is
\begin{align}
  U_{N+1} &= U_{N} + \delta t_{N+1}\left(k^{I}_{2} + k^{E}_{2}\right) = U_{N} + \delta t_{N+1}\left(\Delta\bar{U}_{2} + F^{E}\left(t_{N} + \frac{\delta t_{N+1}}{2},\mathbf{x},\bar{U}_{2}\right)\right)\\
  &=  U_{N} + \delta t_{N+1}\left(\frac{2}{\delta t_{N+1}}\bar{U}_{2} -\frac{2}{\delta t^{N+1}}U^{N} - F^{E}(t_{N},\mathbf{x},\bar{U}_{1})+ F^{E}\left(t_{N} + \frac{\delta t_{N+1}}{2},\mathbf{x},\bar{U}_{2}\right)\right)\\
  &= 2U_{2}^{I} - U_{N} + \delta t_{N+1}\left(F^{E}\left(t_{N} + \frac{\delta t_{N+1}}{2},\mathbf{x},\bar{U}_{2}\right)  - F^{E}(t_{N},\mathbf{x},U_{N})\right).
\end{align}
An important aspect of IMEXRK2 is that we obtain a second order method by only solving \eqref{eq:ModHelmEq}-\eqref{eq:ModHelmEqBC} once, i.e. only one intermediate stage is required.

 An adaptive time-stepper can be constructed by coupling IMEXRK2 with a method of lower order. A simple IMEX scheme of first order is the Forward-Backward Euler scheme, with coefficients given in Table \ref{tab:butcherFBE}. Applied to  the heat equation \eqref{eq:HeatEq} we have
\begin{equation}
  \label{eq:imexBEHeat}
  \frac{U_{N+1}(\mathbf{x})}{\delta t_{N+1}} - \Delta U_{N+1}(\mathbf{x}) = G(t_{N},\mathbf{x}) + \frac{U_{N}(\mathbf{x})}{\delta t_{N+1}}.
\end{equation}

\subsubsection{IMEXRK34}
\label{sss:discretisation_temporal_imexrk4}
The IMEKRK34 scheme is a coupled third and forth order scheme, see \ref{tab:butcherIMEXRK4E} and \ref{tab:butcherIMEXRK4I} for the associated Butcher tableaus. It has two sets of six stages $\{k_{i}^{\sigma}\}_{i = }$ for
$\sigma=I,E$, but only five implicit stages need to be solved for every iterate in time \cite{KENNEDY2003139}.
This is due to $k^{I}_{6}$ at $t_{N}$ is equal to $k^{I}_{1}$ at $t_{N+1}$ for $N>1$, a property sometimes referred to as \emph{first same as last}, or FSAL. For $N = 0$ the first stage must be given by supplementary initial data $\Delta U_{0}$. Otherwise the procedure is exactly as described in \ref{ss:discretisation_temporal}: for a given $i$ solve \eqref{eq:IMEXmodHelm} for $\bar{U}_{i}$. Once known extract $k^{I}_{i} = \Delta \bar{U}_{i}$ from \eqref{eq:IMEXmodHelm} and compute $k_{i}^{E} = F^{E}\left(t_{N}+\delta t_{N+1}c_{i}^{E},\bar{U}_{i}\right)$ explicitly and start over for $i + 1$ until all six stages are known. An approximate solution  $U_{N+1}$ at $t_{N+1}$ is given by  \eqref{eq:imexsol}, which is a fourth order approximation. The third order approximation $\tilde{U}_{N+1}$ is given by \eqref{eq:imexsol} as well, but with the coefficients $\{\tilde{b}_{j}^{\sigma}\}_{j = 1}^{N_{S}}$ instead of $\{b_{j}^{\sigma}\}_{j = 1}^{N_{S}}$.
\subsubsection{Adaptivity}
\label{sss:discretisation_temporal_adaptivity}
Denote the solution given by Forward-Backward Euler or the third order method in IMEXRK34 as $\tilde{U}(\mathbf{x})$. At each discrete time instance $t_{N+1} = \delta t_{N+1} + t_{N}$ for some $\delta t_{N+1}$ we compute $U_{N+1}(\mathbf{x})$ and $\tilde{U}_{N+1}(\mathbf{x})$. The relative temporal error is approximated by
\begin{equation}
  \label{eq:estError}
r = \frac{\|U_{N+1}-\tilde{U}_{N+1}\|}{\|U_{N+1}\|},
\end{equation}
where $\|\cdot\|$ will henceforth represent the standard discrete $\ell_{2}$-norm \eqref{eq:ell2}, unless stated otherwise.  If $r$ is less than some tolerance $\text{TOL}$, then $U_{N+1}(\mathbf{x})$ is accepted as solution at time $t_{N+1}$, otherwise  $\delta t_{N+1}$ is updated via
\begin{equation}
\label{eq:adaptiveScheme}
\delta t_{N+1,\textit{NEW}} = \delta t_{N+1,\textit{OLD}} * (0.9*TOL/r)^{\frac{1}{p+1}},
\end{equation}
where $p = 2$ from the order of the IMEXRK2 scheme and $p = 4$ for IMEXRK34. The value $0.9$ is a safety factor. Even if the solution is accepted the step size is updated by the scheme \eqref{eq:adaptiveScheme}, thus growth is possible if appropriate.  See the flowchart in appendix \ref{appendix:flowchart} for a graphical overview.\\

\begin{algorithm}\captionsetup{labelfont={sc,bf}, labelsep=newline}
\caption{Adaptive time stepping}\label{euclid}
\begin{algorithmic}[1]
\Procedure{Step in time }{step size $\delta t_{N+1}$}
\BState \emph{top}:
\For {$i = 1:N_{S}$}
\State Solve for $k_{i}^{I}$, compute $k_{i}^{E}$
\EndFor
\State {Compute $U_{N+1}$ and $\tilde{U}_{N+1}$}
\State $r \gets \frac{\|U-\tilde{U}\|}{\|U\|}$
\State $\delta t_{N+1} \gets \delta t_{N+1} * (0.9*TOL/r)^{\frac{1}{p+1}}$
\If {$r < \text{TOL}$}
\State $\textit{Keep solution } U_{N+1}$
\State $t_{N}\gets t_{N} + \delta t_{N+1}$
\EndIf
\State \textbf{goto} \emph{top}.
\EndProcedure
\end{algorithmic}

\end{algorithm}

\label{s:appendixbutchertableau}
\FloatBarrier
\begin{table}
\centering
{\def\arraystretch{1.3}
\begin{tabular}{c|cc}
$0$ & $0$ & $0$\\
$1$ & $0$ & $1$ \\
\hline
 & $0$ & $1$
\end{tabular}
}\hspace{3cm}
{\def\arraystretch{1.3}
\begin{tabular}{c|cc}
$0$ & $0$ & $0$\\
$1$ & $1$ & $0$ \\
\hline
 & $1$ & $0$
\end{tabular}
}
\caption{Coeffcients for the IMEX scheme Forward-Backward Euler. The left and right tables correspond to $\sigma = I$ and to $\sigma = E$, respectively.}
\label{tab:butcherFBE}
\end{table}

\begin{table}
\centering
{\def\arraystretch{1.3}
\begin{tabular}{c|cc}
$0$ & $0$ & $0$\\
$\frac{1}{2}$ & $0$ & $\frac{1}{2}$ \\
\hline
 & $0$ & $1$
\end{tabular}
}\hspace{3cm}
{\def\arraystretch{1.3}
\begin{tabular}{c|cc}
$0$ & $0$ & $0$\\
$\frac{1}{2}$ & $\frac{1}{2}$ & $0$ \\
\hline
 & $0$ & $1$
\end{tabular}
}
\caption{Coeffcients for the IMEXRK2 scheme. The left and right tables correspond to $\sigma = I$ and to $\sigma = E$, respectively.}
\label{tab:butcherIMEXRK2}
\end{table}

\FloatBarrier
\begin{table}
\centering
{\def\arraystretch{1.3}
\begin{tabular}{c|cccccc}
$0$ & $0$ & $0$ & $0$ & $0$ & $0$ & $0$ \\
$\frac{1}{2}$ & $\frac{1}{2}$ & $0$ & $0$ & $0$ & $0$ & $0$\\
$\frac{83}{250}$ & $\frac{13861}{62500}$ & $\frac{6889}{62500}$ & $0$ & $0$ & $0$ & $0$ \\
$\frac{31}{50}$ & $\frac{-116923316275}{2393684061468}$ & $\frac{-2731218467317}{15368042101831}$ & $\frac{9408046702089}{11113171139209}$ & $0$ & $0$ & $0$ \\
$\frac{17}{20}$ &  $\frac{-451086348788}{2902428689909}$ & $\frac{-2682348792572}{7519795681897}$ & $\frac{12662868775082}{11960479115383}$ & $\frac{3355817975965}{11060851509271}$ & $0$ & $0$ \\
$1$ & $\frac{647845179188}{3216320057751}$ & $\frac{73281519250}{8382639484533}$ & $\frac{552539513391}{3454668386233}$ & $\frac{3354512671639}{8306763924573}$ & $\frac{4040}{17871}$ & $0$\\
\hline
 $b^{E}_{i}$ & $\frac{82889}{524892}$ & $0$ & $\frac{15625}{83664}$ & $\frac{69875}{102672}$ & $\frac{-2260}{8211}$ & $\frac{1}{4}$\\
 \hline
  $\tilde{b}^{E}_{i}$ & $\frac{4586570599}{29645900160}$ & $0$ & $\frac{178811875}{945068544}$ & $\frac{814220225}{1159782912}$ & $\frac{-3700637}{11593932}$ & $\frac{61727}{225920}$
\end{tabular}
}
\caption{The coefficients $\{a_{i,j}^{E}\}_{i,j = 1}^{6}$, $\{b_{j}^{E}\}_{j = 1}^{6}$, $\{\tilde{b}_{j}^{E}\}_{j = 1}^{6}$ and $\{c_{i}^{E}\}_{i= 1}^{6}$ for IMEXRK4.}
\label{tab:butcherIMEXRK4E}
\end{table}


%
%

\begin{table}
\centering
{\def\arraystretch{1.3}
\begin{tabular}{c|cccccc}
$0$ & $0$ & $0$ & $0$ & $0$ & $0$ & $0$ \\
$\frac{1}{2}$ & $\frac{1}{4}$ & $\frac{1}{4}$ & $0$ & $0$ & $0$ &  $0$\\
$\frac{83}{250}$ & $\frac{8611}{62500}$ & $\frac{-1743}{31250}$ & $\frac{1}{4}$ & $0$ & $0$ &  $0$\\
$\frac{31}{50}$ & $\frac{5012029}{34652500}$ & $\frac{-654441}{2922500}$ & $\frac{174375}{388108}$ & $\frac{1}{4}$ & $0$ & $0$\\
$\frac{17}{20}$ & $\frac{15267082809}{155376265600}$ & $\frac{-71443401}{120774400}$ & $\frac{730878875}{902184768}$ & $\frac{2285395}{8070912}$ &  $\frac{1}{4}$ & $0$\\
$1$ & $\frac{82889}{524892}$ & $0$ & $\frac{15625}{83664}$ &  $\frac{69875}{102672}$ & $\frac{-2260}{8211}$ & $\frac{1}{4}$ \\
\hline
 $b^{E}_{i}$ & $\frac{82889}{524892}$ & $0$ & $\frac{15625}{83664}$ & $\frac{69875}{102672}$ & $\frac{-2260}{8211}$ & $\frac{1}{4}$\\
 \hline
  $\tilde{b}^{E}_{i}$ & $\frac{4586570599}{29645900160}$ & $0$ & $\frac{178811875}{945068544}$ & $\frac{814220225}{1159782912}$ & $\frac{-3700637}{11593932}$ & $\frac{61727}{225920}$
\end{tabular}
}
\caption{The coefficients $\{a_{i,j}^{I}\}_{i,j = 1}^{6}$, $\{b_{j}^{I}\}_{j = 1}^{6}$, $\{\tilde{b}_{j}^{E}\}_{j = 1}^{6}$ and $\{c_{i}^{I}\}_{i= 1}^{6}$ for IMEXRK4.}
\label{tab:butcherIMEXRK4I}
\end{table}
\FloatBarrier

\FloatBarrier

\begin{figure}
  \label{chart:imex}
  \centering
\begin{tikzpicture}[node distance=2cm]

\node (start) [startstop] {Given $U_{N}$, the next solution $U_{N+1}$ is obtained by \eqref{eq:imexsol}. The implicit stages \eqref{eq:stagek} must be solved for. Explicit stages are computed directly. The same stages are used to compute low order approximation $\hat{U}_{N+1}$, used for adaptive time stepper.};
\node (input) [io, below of=start] {Input: time step $\delta t_{N+1}$, Butcher tableau \ref{tab:butcher}, solution $U_{N}$ at $t_{N}$, first implicit stage $k_{1}^{I}$ and Dirichlet boundary data.};
\node (modhelm) [blackbox, below of=input] {Solve for $U^{I}_{i}$ from \eqref{eq:IMEXmodHelm}, by solving \eqref{eq:ModHelmEq}--\eqref{eq:ModHelmEqBC} as shown in Figure \ref{fig:flowchartSolveModHelm}.};
\node (computeStage) [process, below of=modhelm] {With $U^{I}_{i}$ known, compute implicit stage $k^{I}_{i}$ by \eqref{eq:computeimplicitstage} and then explicit stage $k^{E}_{i} = F^{E}(t_{N}+\delta t_{N+1}c^{E}_{i},U^{I}_{i})$.};
\node (dec1) [decision, below of=computeStage] {Are all stages computed?};
\node (computenextU) [process, below of=dec1] {Compute $U_{N+1}$ and $\hat{U}_{N+1}$  from \eqref{eq:imexsol} and update $\delta t_{N+1}$ as in  \eqref{eq:adaptiveScheme}. Approximate error $r$ \eqref{eq:estError}. };
\node (dec2) [decision, below of=computenextU] {Is $r <$ TOL?};
\node (output) [io, below of=dec2] {Output: Solution $U_{N+1}$, time $t_{N+1} = t_{N} + \delta t_{N+1}$, time step $\delta t_{N+1}$ and $k_{6}^{I}$, which is $k^{I}_{i}$ for the next iteration in time.};
\node (stop) [startstop, below of=output] {Approximate solution $U_{N+1}$ to the diffusion equation at time $t_{N+1}$.};
\draw [arrow] (input) -- node[midway,xshift=1.5cm] {\texttt{for} $i = 1:N_{S}$} (modhelm);
\draw [arrow] (modhelm) -- (computeStage);
\draw [arrow] (computeStage) -- (dec1);
\draw [arrow] (dec1) -- (computenextU);
\draw [arrow] (dec1.west) -- node[midway,yshift=0.3cm] {No ($i<N_{S}$)} (-5,-8) -- (-5,-4) -- (modhelm.west);
\draw [arrow] (dec1.south) -- node[midway,xshift=1.3cm] {Yes ($i=N_{S}$)}  (computenextU.north);
\draw [arrow] (computenextU.south) --   (dec2.north);
\draw [arrow] (dec2.south) -- node[midway,xshift=0.5cm] {Yes}  (output.north);
\draw [arrow] (dec2.west) -- node[midway,yshift=0.3cm] {No} (-6,-12) -- (-6,-2) -- (input.west);
\end{tikzpicture}
\caption{Flowchart over the procedure for updating the approximate solution $U_{N}$ at $t_{N}$ for the heat equation \eqref{eq:HeatEq}--\eqref{eq:HeatEqBC}. Note that the grey block corresponds to the flowchart in Figure \ref{fig:flowchartSolveModHelm}.}\label{fig:flowchartSolveHeat}
\end{figure}
\FloatBarrier

\section{Flowchart over solution procedure}
\nopagebreak
\label{appendix:flowchart}

\begin{figure}
  \label{chart:modhelm}
\centering
\begin{tikzpicture}[node distance=2cm]
\node (start) [startstop] {Solve modified Helmholtz equation   \begin{align*}
  \alpha^{2}u(\mathbf{x}) - \Delta u(\mathbf{x}) = f(\mathbf{x}),\quad \mathbf{x} \in \Omega,\\
  u(\mathbf{x}) = g(\mathbf{x}),\quad \mathbf{x} \in \Gamma.
\end{align*}
Decompose $u = u^{P} + u^{H}$.
};
\node (input) [io, below of=start,yshift = -0.7cm] {Input: $\alpha^2$, $f$ and $g$};
\node (up1) [process, below left of=input,xshift=-4cm,yshift=-1cm] {Construct extension $f^{e}(\x)$ of $f$ with PUX, see Section \ref{sss:discretisation_aptive_modhelm_pux}.};
\node (up2) [process, below  of=up1,yshift=-0.5cm] {Solve $\alpha^{2}u^{P}-\Delta u = f^{e}$ in Fourier space with FFT, Section  \ref{sss:discretisation_adaptive_modhelm_inhomo}.\\ Compute $u^{P}$ in $\Omega$ and $u^{P}|_{\Gamma}$ on $\Gamma$.};
\node (uh1) [process,right of=up2,xshift=8cm] {Solve $
  \alpha^{2}u^{H}-\Delta u^{H} = 0$ in $\Omega$ with $u^{H} = \tilde{g}$ on $\Gamma$ as in Section \ref{sss:discretisation_adaptive_modhelm_homo}.\\Compute $u^{H}$ in $\Omega$};
\node (output) [io, below of=input,yshift=-6cm] {The solution to the modified Helmholtz equation is $u = u^{P} + u^{H}$.\\
Output: $u$};
\node (stop) [startstop, below of=output] {Done};
\draw [arrow] (input) |- (up1);
\draw [arrow] (up1) -- (up2);
\draw [arrow] (up2) -- node[midway,yshift=0.3cm] {$\tilde{g}=g-u^{P}|_{\Gamma}$} (uh1);
\draw [arrow] (uh1.south) -- node[anchor=west,xshift=-0.8cm,yshift=0.1cm] {$u^{H}$} (output.north);
\draw [arrow] (up2.south) -- node[anchor=west,xshift=0.8cm,yshift=0cm] {$u^{P}$} (output.north);
\end{tikzpicture}
\caption{Flowchart over the procedure for solving the modified Helmholtz equation \eqref{eq:ModHelmEq}--\eqref{eq:ModHelmEqBC}.}\label{fig:flowchartSolveModHelm}
\end{figure}
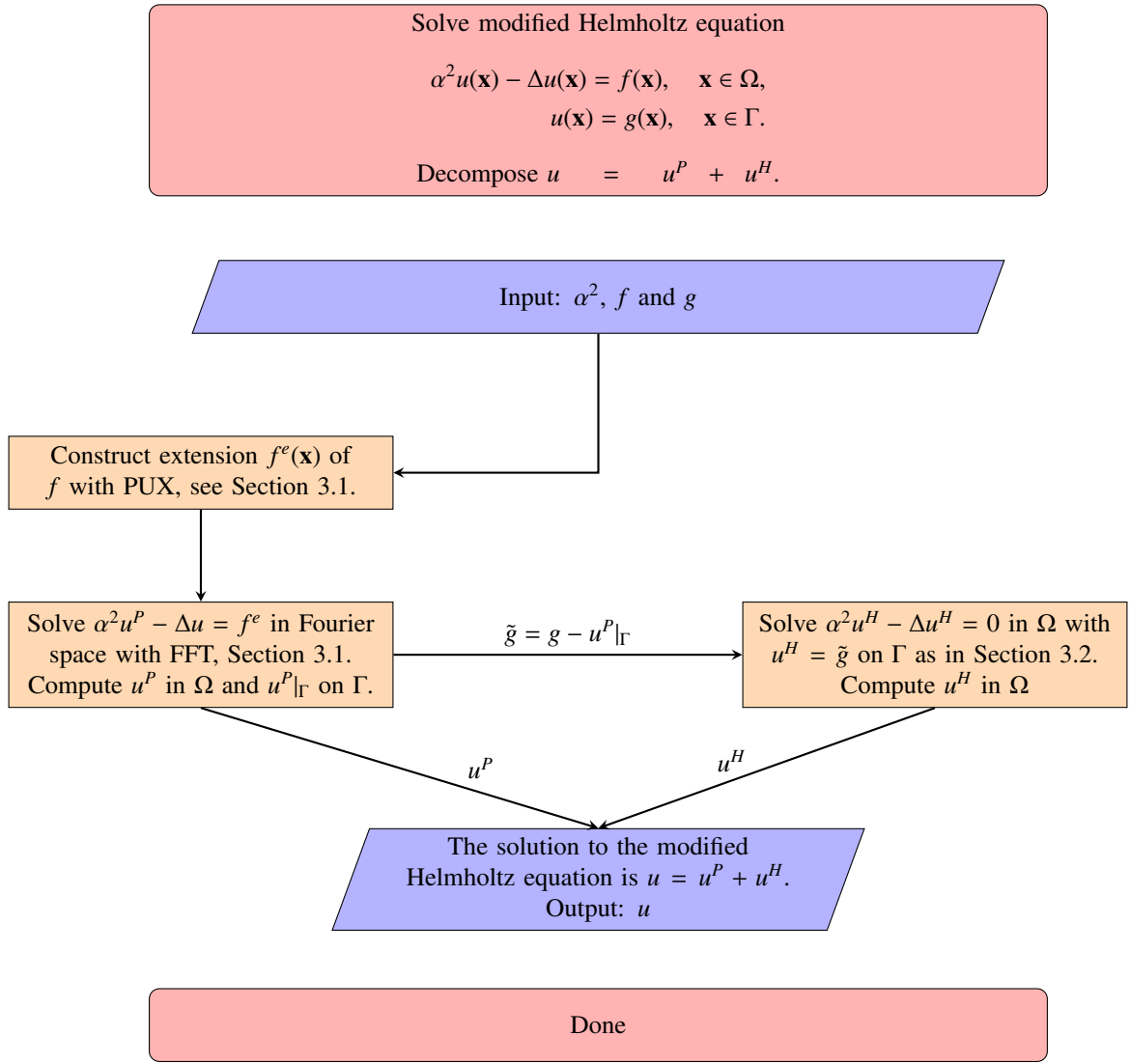
\FloatBarrier
\bibliographystyle{model1-num-names}
\bibliography{sections/references}
\end{document}